\documentclass[preprint,review,12pt]{elsarticle}
\usepackage{graphicx}
\usepackage{amsmath}
\usepackage{amssymb}
\usepackage{amsfonts}	
\usepackage{amsthm}
\usepackage{subfigure}
\usepackage{bm}
\usepackage{url}
\usepackage{float}
\usepackage{mathrsfs}
\usepackage{tikz}
\usepackage{booktabs}
\usepackage{listings}
\usepackage{textcomp}
\usepackage{pifont}
\usepackage{calc}
\usepackage{xspace}
\usepackage{rotating}
\usepackage{pgf,pgfplots}

\usetikzlibrary{fadings,calc,trees,positioning,arrows,chains,shapes.geometric,%
    decorations.pathreplacing,decorations.pathmorphing,shapes,%
    matrix,shapes.symbols,scopes,trees,patterns}
          
\newcommand{\tp}{p}

\newcommand{\Xtp}{\widehat{\bf X}_{{p}}}
\newcommand{\Xtpp}{\widehat{\bf X}_{{p}^+}}
\newcommand{\Xtpm}{\widehat{\bf X}_{{p}^-}}
\newcommand{\BXtp}{\widehat{\bf \cal X}_{{p}}}

\newcommand{\Xp}{{\bf X}_{p}}
\newcommand{\Jtp}{\widehat{\sf J}_{cp}}
\newcommand{\Htp}{{\sf H}_{cp}}
\newcommand{\Atp}{\widehat{A}_{cp}}
\newcommand{\Ftp}{F_{p}}
\newcommand{\figref}{{\sc Fig.}}
\newcommand{\tabref}{{\sc Tab.}}
\newcommand{\VU}{{\bf U}}
\newcommand{\VX}{{\bf X}}
\newcommand{\VN}{{\bf N}}
\newcommand{\Dt}{\dfrac{d}{dt}}
\newcommand{\CR}{{\cal R}}
\newcommand{\CM}{{\sf M}}
\newcommand{\CZ}{{\sf Z}}

\newcommand{\DR}{{\sf R}}
\newcommand{\CF}{{\cal F}}
\newcommand{\CN}{{\cal N}}
\newcommand{\CC}{{\cal C}}
\newcommand{\eps}{\varepsilon}
\newcommand{\FPC}{{\bf F}_{pc}}
\newcommand{\TOm}{\widetilde{\Omega}}
\newcommand{\TX}{\widetilde{\bf X}}
\newcommand{\HVX}{\widehat{\bf X}}

\newcommand{\DT}{{\sf T}}
\newcommand{\bs }{\boldsymbol}
\theoremstyle{definition}

\journal{Journal of Computational Physics}

\usepackage[linecolor=blue!70!white, backgroundcolor=pink!20!white,bordercolor=red]{todonotes}
\let\origintodo\todo
\newcommand{\xtodo}[2][]{\origintodo[#1]{#2}\xspace}
\let\todo\xtodo

\begin{document}

\begin{frontmatter}

\title{A multi-material CCALE-MOF approach in cylindrical geometry}

\author[ecn]{Marie Billaud Friess}
\ead{marie.billaud-friess@ec-nantes.fr}
\author[celia]{J\'er\^ome Breil \corref{cor1}}
\ead{breil@celia.u-bordeaux1.fr}
\author[celia]{St\'ephane Galera}
\author[cea]{Pierre-Henri Maire}
\author[lanl]{Mikhail Shashkov}
\cortext[cor1]{Corresponding author}

\address[ecn]{GeM, UMR CNRS 6183, Ecole Centrale Nantes, 1 rue de la Noë, BP 92101, 44321 NANTES Cedex 3, France}
\address[celia]{Univ. Bordeaux, CEA, CNRS, CELIA, UMR5107, F-33400 Talence, France}
\address[cea]{CEA CESTA, BP 2, 33114 Le Barp Cedex, France }
\address[lanl]{Los Alamos National Laboratory, XCP-4, Los Alamos, NM 87545, USA}

\begin{abstract}
In this paper we present recent developments concerning a Cell-Centered Arbitrary Lagrangian Eulerian (CCALE) strategy using the Moment Of Fluid (MOF) interface reconstruction for the numerical simulation of multi-material compressible fluid flows on general unstructured grids in cylindrical geometries. Especially, our attention is focused here on the following points. First, we propose a new formulation of the scheme used during the Lagrangian phase in the particular case of axisymmetric geometries. Then, the MOF method is considered for multi-interface reconstruction in cylindrical geometry. Subsequently, a method devoted to the rezoning of polar meshes is  detailed. Finally, a generalization of the hybrid remapping to cylindrical geometries is presented. These explorations are validated by mean of several test cases that clearly illustrate the robustness and accuracy of the new method.
\end{abstract}


\end{frontmatter}

\section{Introduction}

\begin{figure}[h!]
\centering
\begin{tikzpicture}[scale=0.4]
\path (0 ,20) node [rectangle , inner sep=0.3cm, rounded corners=0.3cm, draw =black, line width=0.5mm, dashed](init) 
{\footnotesize  \it \sc Initialization};
\path (15 ,13) node [rectangle , inner sep=0.3cm, rounded corners=0.3cm, draw =black, line width=0.5mm, dashed](t) {\it $t^{n+1}=t^n+\Delta t$};
\path (0,14) node [rectangle , inner sep=0.3cm, rounded corners=0.3cm, draw =black, line width=0.5mm](lag)  {\footnotesize {\sc Lagrangian phase}};
\path (0,9) node [rectangle , inner sep=0.3cm, rounded corners=0.3cm, draw =black, line width=0.5mm](mof)  {\footnotesize   \sc
\begin{tabular}{c}
Thermodynamical closure \\
for multi-material flows
\end{tabular}
};
\path (0,4) node [rectangle , inner sep=0.3cm, rounded corners=0.3cm, draw =black,line width=0.5mm](int)  {\footnotesize  \sc Interface reconstruction};
\path (15,4) node [rectangle , inner sep=0.3cm, rounded corners=0.3cm, draw =black,line width=0.5mm](reg)  {\footnotesize \sc  Rezoning phase};
\path (15,8) node [rectangle , inner sep=0.3cm, rounded corners=0.3cm, draw =black,line width=0.5mm](proj)  {\footnotesize \sc  Remapping phase};
\draw[->, line width=0.3mm] (init)--(lag);
\draw[->, line width=0.3mm] (lag)--(mof);
\draw[->, line width=0.3mm] (mof)--(int);
\draw[->, line width=0.3mm] (int)--(reg);
\draw[->, line width=0.3mm] (reg)--(proj);
\draw[->, line width=0.3mm] (proj)--(t);
\draw[- , line width=0.3mm] (15,17)--(t);
\draw[<-, line width=0.3mm] (0,17)--(15,17);
\end{tikzpicture}
\caption{Multi-material CCALE algorithm flowchart.} \label{fig:0}
\end{figure}
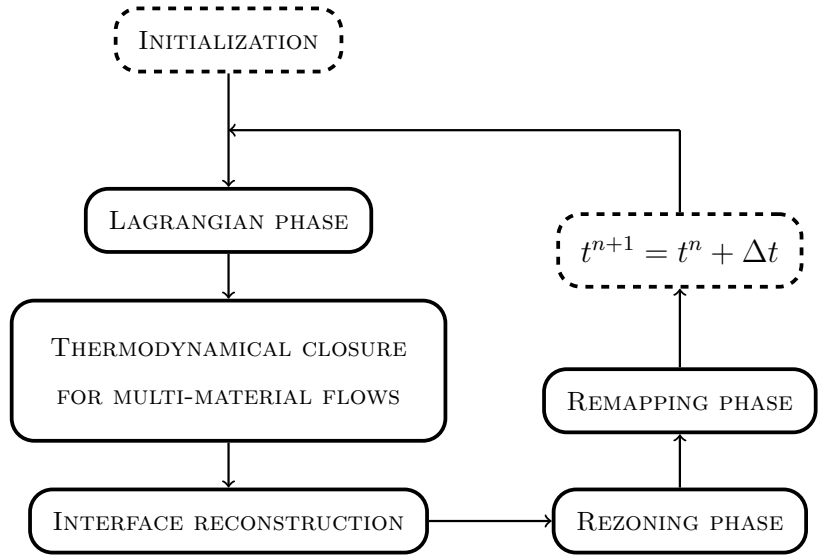

In this work, we consider the simulation of multi-material compressible flows on unstructured meshes in cylindrical geometry. For this, we adopt an ALE description \cite{Hirt1} that has the great advantage to combine the best features of both Eulerian and Lagrangian approaches. Indeed, this choice is not only well adapted to naturally track free surfaces and interfaces between different fluids as purely Lagrangian methods, but also to handle flow distortion as Eulerian methods. Here, a CCALE \cite{Galera1,Galera2} approach is particularly considered whose the main elements are as follow.\\ 
As depicted on figure \figref \ref{fig:0}, the first step of the algorithm relies on an explicit Lagrangian phase in which the physical variables and grid are updated thanks to  a slightly modified version of the Explicite Unstructured Cell-Centered Lagrangian HYDrodynamics (EUCCLHYD) scheme \cite{Maire1,Maire2,Maire3} in cylindrical coordinates. Recently, new investigations have been made about cell-centered Lagrangian schemes \cite{Barlow1,Despres1}. The scheme presented in this paper is a modified version of the area weighted finite volume scheme of \cite{Maire1}.
Then, multi-material flows treatment is done thanks to specific interface capturing method. This choice allows to track the volume fraction of each material used for the thermodynamical closure relying on the equal strain rates assumption. This approach is quite simple to implement and to use and remains sufficient in almost cases \cite{Galera2,Shashkov1}. 
This, leads to constant evolution of the volume fraction during the Lagrangian phase. Such an approach allows to reconstruct with accuracy the interface between each material. In this context, many development have been done for 2D Cartesian geometries. First, a previous version of the CCALE algorithm solving two-material compressible flows using a Volume Of Fluid (VOF) have been proposed in \cite{Breil2,Galera2}. Then an extension to Moment Of Fluid (MOF) approach has been considered to enhance multi-material (more than two components) flows in \cite{Dyadechko1,Galera1}. 
Subsequently, a rezoning phase is realized.  It consists in moving the Lagrangian nodes to improve the geometric quality of the grid \cite{Knupp1}.
Finally,  the physical variables are conservatively interpolated from the Lagrangian grid onto the new rezoned one  during the remapping  phase. 
Here an extension of the hybrib remapping \cite{Berndt1} to cylindrical geometries is introduced. We want to notice that in ALE framework using cell-centered formulation, this phase is straightforward. 
In the lines of these works, the main goal of this paper is to extend the CCALE-MOF algorithm to treat both Cartesian and cylindrical geometry. To this end, several modifications are given to the algorithm previously presented. In a first part, we propose a new formulation of the numerical scheme introduced in \cite{Maire1} for treating axisymmetric geometries during the Lagrangian phase. To build this scheme, an area-weighted formulation of the Lagrangian system of equations is proposed. Then, this system of equations is discretized using a cell-centered finite volume (FV) scheme. Contrary to \cite{Maire1} in which fluxes are directly deduced from the Geometric Conservation Law (GCL) constraint, here a simpler formulation that gives similar results is retained. 
These two main choices lead to a robust first-order scheme conservative for the total energy that has the great advantage to preserve spherical symmetry for one-dimensional flow on uniform angular polar grids. The high order extension has been performed using the Generalized Riemann Problem (GRP) described in \cite{Maire1}.  
To treat interface flows,  a MOF interface reconstruction method is retained in the sequel. Once again, the difficulty here is to propose a natural and consistent adaptation of this approach able to treat axisymmetric interface flows. To this end, formulations of the moments needed to track interface are revisited for cylindrical coordinates as in \cite{Anbarlooei1}. This leads to an accurate and second order interface reconstruction method that allows to treat multi-material (more than two) interfaces in the lines of \cite{Dyadechko1}.
The third part of this study is dedicated to recent enhancement of the rezoning algorithm to improve the mesh quality during computation especially on polar meshes. As it is done in \cite{Galera1,Galera2}, mesh rezoning is based on the Condition Number Smoothing (CNS) \cite{Knupp1} algorithm on unstructured meshes. Moreover, when used for polar meshes, it is well known that CNS algorithm pushes the nodes toward the origin deteriorating the mesh quality. To avoid this drawback, the main idea developed in this paper is to adapt CNS algorithm to polar grids. Then, extension to unstructured grids (Cartesian-polar) is also explored. 
Finally, a generalization of the remapping procedure to cylindrical geometries is proposed. Here, an efficient method adapted to multi-material flows is presented. The main idea is to use an hybrid remapping that combine the main advantages of the swept-face and multi-material cell-intersection remapping as in \cite{Berndt1,Galera1}. Finally, a specific attention is done to polynomial integration that preserves the method efficiency.\\ 

The paper is structured as follows. We detail in the second section a new formulation of the first-order area weighted Lagrangian scheme used for axisymmetric geometries. Further extensions to high-order are notably detailed in \cite{Maire1}.
Afterwards, the extension of the MOF axisymmetric interface reconstruction method is presented for treating  multimaterial flows.
 Then, we describe the General Condition Number Smoothing (GCNS) algorithm for unstructured meshes. Finally, the description of the new hybrid remapping procedure for cylindrical geometry is done. For a complete description of the CCALE-MOF method see \cite{Galera1,Galera2}, except new advances presented in this paper. Then presentation of  numerical experiments is made in Section 4. They demonstrate not only the robustness and the accuracy of the present methodology but also its ability to handle successfully complex two-dimensional multi-material fluid flows notably computed for axisymmetric geometries. Finally concluding remarks and perspectives about future works are given in the last section.

\section{Lagrangian phase in axisymmetric geometry}

 In this part, an extension of the cell-centered Lagrangian scheme \cite{Maire2,Maire3} is presented for the numerical simulation of compressible flows in pseudo-Cartesian geometries for unstructured meshes as in \cite{Maire1}. This choice has the great advantage to treat both axisymmetric and Cartesian geometries. In this paper, a new and simple formulation of the scheme introduced in \cite{Maire1} for first-order approximation is proposed. To this end, an area weighted formulation of classical Lagrangian equations is first introduced. Then these equations are discretized with a node-centered approximate Riemann solver.

\subsection{Governing equations}

During the Lagrangian phase, the rates of change of volume, mass, momentum and total energy are computed assuming that discretized volumes move following the flow. Thus, each arbitrary volume $V(t)$ depending on the time $t>0$ moves satisfying the following system of equations
\begin{align}
\displaystyle &\Dt \int_{V} \rho     dV                                    = 0
\label{eq:lag1}, \\[0.25cm]
\displaystyle & \Dt \int_{V}          dV - \displaystyle \int_{V} \nabla \cdot \VU dV    = 0, 
\label{eq:lag2} \\[0.25cm]
\displaystyle & \Dt \int_{V} \rho \VU dV + \displaystyle \int_{V} \nabla P         dV    = 0, 
\label{eq:lag3} \\[0.25cm]
\displaystyle &\Dt \int_{V} \rho   E dV + \displaystyle \int_{V} \nabla \cdot (P\VU) dV = 0
\label{eq:lag4},
\end{align}
where $\Dt$ is the Lagrangian derivative and $\rho, \VU, P, E$ are respectively the density, velocity, pressure and total energy. In addition, this system is closed thanks to an equation of state (EOS) as
\begin{equation}
p=p(\rho,\eps),
\end{equation}
with the internal energy $\eps$ defined as $ \eps = E - |\VU|^2/2$. At last, we have local kinematic equation 
\begin{equation}
\dfrac{d \VX}{dt} = \VU, \quad \VX(0) = \VX^0,
\end{equation}
with $\VX$ the location of a point of the control volume surface $S(t)$, at time $t>0$ and $\VX^0$ its initial value. This equation is equivalent to \eqref{eq:lag2} also known as geometric conservation of law (GCL). \\

\subsection{Area-weighted formulation}

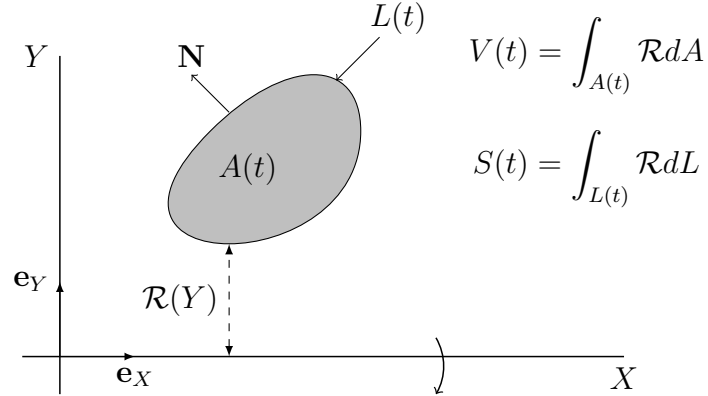
\begin{figure}[h!]
\centering
\begin{tikzpicture}[xscale=0.5,yscale=0.5]
\tikzstyle{fleche1}=[-,>=latex,line width=0.2mm]
\tikzstyle{fleche2}=[<->,>=latex]
\tikzstyle{fleche3}=[->,>=latex]
\draw[fleche1] (-1,0)--(15,0); 
\draw[fleche1] (0,-1)--(0,8);
\draw[fleche3] (0,0)--(0,2); 
\node[left] at (0,2) {${\bf e}_Y$};
\node[below] at (2,0) {${\bf e}_X$};
\node[left] at (0,8) {$Y$};
\node[below] at (15,0) {$X$};
\draw[fleche3] (0,0)--(2,0);
\filldraw[gray!50,draw=black] (4,6) .. controls (6,8) and (8,8) .. (8,6).. controls (8,2) and (0,2) .. (4,6);
\node at (5,5) {$A(t)$};
\node at (9,9) {$L(t)$};
\draw[<->,>=latex,dashed] (4.5,0) -- (4.5,3);
\node[left] at (4.5,1.5) {${\cal R}(Y)$};
\draw[<-,line width=0.2mm] (10,-1) arc (-30:30:1.5);
\draw[->] (4.5,6.5)--(3.5,7.5);
\node at (3.5,8) {$\VN$};
\draw[->] (8.5,8.5)--(7.35,7.35);
\node at (14,8) {$\displaystyle V(t) = \int_{A(t)}  \CR dA$};
\node at (14,5) {$\displaystyle S(t) = \int_{L(t)} \CR dL$};
\end{tikzpicture}  
\caption{Notations related to the pseudo-cartesian grid. \label{fig:5}}
\end{figure}

For defining the differential operators  used in the system of Lagrangian equation \eqref{eq:lag1}-\eqref{eq:lag4} a pseudo-Cartesian reference frame $\{0,X,Y\}$ for the orthonormal basis $({\bf e}_X, {\bf e}_Y)$ is used  (see \figref \ref{fig:5}). Thus each point is localized by means of its positions $X$ and $\CR(Y) = 1-\alpha + \alpha Y$ the pseudo-radius. When $\alpha =0$, the Lagrangian equations for Cartesian geometry are recovered, otherwise for $\alpha = 1$ this corresponds to axisymmetric equations. In this way, axisymmetric geometry is obtained from Cartesian one through a rotational symmetry about the $X$-axis. This implies that the volume $V(t)$ is generated by the rotation of the area $A(t)$ about the $X$-axis. In consequence, the element volume $dV$ writes as $dV = \CR dA$ with $dA = dXdY$ the element area in the pseudo-Cartesian frame. In the same manner, the control surface $S(t)$ delimiting $V(t)$ is obtained through the rotation of $L(t)$ the boundary of $A(t)$ and the surface element is given by $dS = \CR dL$. Note that we have omitted the $2\pi$ factor in the evaluation of the element volume.\\

In a such framework, the velocity divergence and the pressure gradient read as follows
\begin{equation}
\nabla \cdot \VU = \dfrac{1}{\CR} \left[ \dfrac{\partial (\CR u) }{\partial X}
+\dfrac{\partial (\CR v) }{\partial Y}
\right], \text{ where } \VU^t =(u,v) 
\end{equation}
and 
\begin{equation}
\nabla P = \left( \dfrac{\partial P}{\partial X} {\bf e}_X +\dfrac{\partial P }{\partial Y} {\bf e}_Y
\right).
\end{equation}
Using the previous definitions and after some calculations using the Green's formula, it is possible to rewrite \eqref{eq:lag1}-\eqref{eq:lag4} at least in two different ways. The first one, obtained without any approximation is the {\it control volume formulation}. When discretized this formulation leads to a conservative scheme for both equations of energy and momentum, and satisfies the local semi-discrete entropy inequality. However, as shown in \cite{Maire1} it does not preserve symmetries. Consequently, an {\it area-weighted formulation} is adopted here leading to a conservative scheme for energy equation that respect spherical geometries. This formulation is deduced from the control volume one assuming that momentum equation \eqref{eq:lag3} is written in Cartesian geometry. Like this, the area-weighted formulation for the Lagrangian equations reads
\begin{align}
&\displaystyle m \Dt \left< \dfrac{1}{\rho} \right>  -    \displaystyle \int_{L}  \VN \cdot \CR \VU dL  = 0, \label{eq:area1}\\[0.25cm]
&\displaystyle m \Dt \left< \VU \right> + \displaystyle    \overline{\CR}  \int_{L} P \VN   dL=  0, \label{eq:area2} \\[0.25cm]
&\displaystyle m \Dt \left< E \right> + \displaystyle \int_{L} P\VN  \cdot  \CR \VU dL = 0,
\label{eq:area3}
\end{align}
where $m= \int_{V} \rho dV$ represents the mass of the volume $V$. Each physical variable per unit of mass ($E, \VU$) is noted as $\phi$, and has its mass density mean value defined by $\left<  \phi \right>= \frac{1}{m} \int_V \rho \phi dV$. The average $\overline{\CR}$ corresponds to ratio $\overline{\CR} = \frac{V}{A}$. In such case, as $m=\rho V$, the momentum equation is solved in Cartesian geometry. For Cartesian case $V=A$, we recover $\overline{\CR} = 1$. Further details on the derivation of this system are available in \cite{Maire1}.

\subsection{Numerical scheme}

\begin{figure}[h!]
\centering
\begin{tikzpicture}[xscale=0.9,yscale=0.9]
\tikzstyle{fleche1}=[-,>=latex]
\tikzstyle{fleche2}=[<->,>=latex]
\tikzstyle{fleche3}=[->,>=latex]
\draw[fleche1] (2,0)--(10,0); 
\draw[fleche1] (3,-1)--(3,7);
\draw[fleche3] (3,0)--(3,1); 
\node[left] at (3,1) {${\bf e}_Y$};
\node[below] at (4,0) {${\bf e}_X$};
\node[left] at (3,7) {$Y$};
\node[below] at (10,0) {$X$};
\draw[fleche3] (3,0)--(4,0);
\filldraw[fill=gray!30,draw=black] (4.8,3.6)--(5.7,1.8)--(9.5,3.1)--(8.,5.9)--cycle node[black] at (6.6,3.8) {$\Omega_c$}  ;
\node[black] at (4.8,3.6) {$\bullet$};
\node[black] at (5.7,1.8) {$\bullet$};
\node[black] at (9.5,3.1) {$\bullet$};
\node[black] at (8.,5.9) {$\bullet$};
\node[right] at (8.2,5.9)   {$p$};
\node[below right] at (10,3.1)   {$p^{+}$};
\node[above] at  (4.8,3.6)   {$p^{-}$};
\node[black] (ppp) at (8.75,4.5){$\times$}; 
\node[black] (ppm) at (6.4,4.75){$\times$};
\draw[->,>=latex] (8.75,4.5) -- (9.63,4.97);
\draw[->,>=latex] (6.4,4.75) -- (5.8162,5.5622);
\draw[->,>=latex] (8,5.9) -- (8.23,6.88);
\node at (9,5.2) {$\VN^+_{pc}$};
\node at (6.5,5.5) {$\VN^-_{pc}$};
\node at (7.5,6.5) {$\VN_{pc}$};
\draw[fill=blue,draw=blue,<->,>=latex](8.75,4.5) -- (8.,5.9);
\draw[fill=blue,draw=blue,<->,>=latex](6.4,4.75) -- (8.,5.9);
\node[line width=0.5mm] at (8,5) {\textcolor{blue}{$L^+_{pc}$}};
\node[line width=0.5mm] at (7.3,4.8) {\textcolor{blue}{$L^-_{pc}$}};
\node at (15,4) {$L_{pc}\VN_{pc} = L_{pc}^+\VN_{pc}^+ +L_{pc}^-\VN_{pc}^-$};
\node at (15,3) {$\CR_c = V_C/A_C$};
\node at (15,2) {$\CR_p = Y_p$};
\draw[dashed,<->,>=latex](8.,5.9)--(8.,0);
\node[left] at (8,1) {${\cal R}_p$};
\draw[<-,line width=0.2mm] (9,-0.5) arc (-20:20:1.5);
\end{tikzpicture}  
\caption{Notations for the cell-centered scheme. \label{fig:6}}
\end{figure}
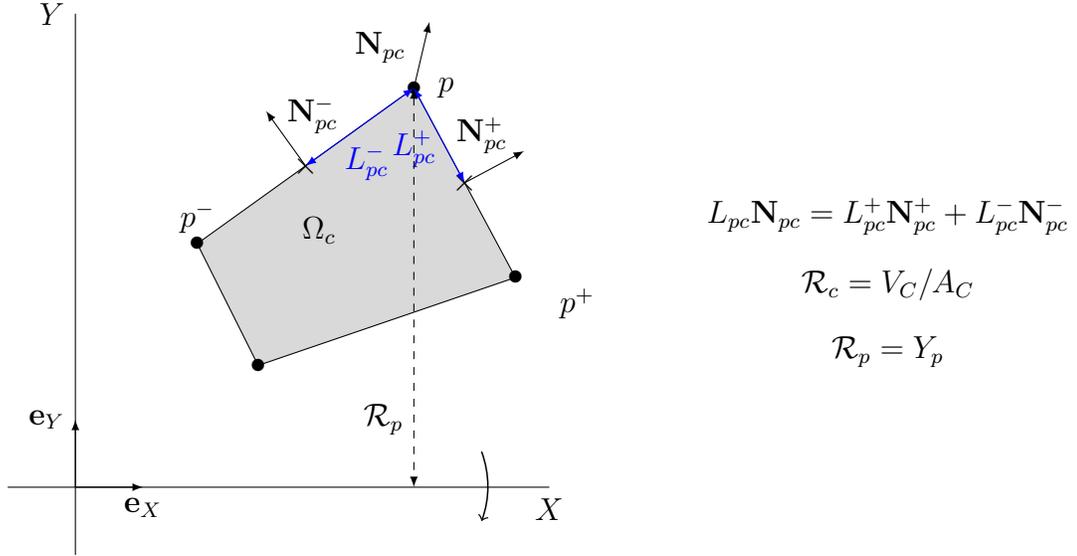

Thereafter, we recall briefly the first order cell-centered Lagrangian scheme introduced in \cite{Maire1}. To this goal, similar notations as \cite{Galera1,Maire1,Maire4} are employed in the sequel. Let us consider a set $\{\Omega_c\}_{c\in\mathbb{N}}$ of non-overlapping polygonal cells that approximates $A(t)$. Each cell noted $\Omega_c$ is assigned a single index $c$. Each vertex of the cell $c$ is labeled with the index $p$ and is localized thanks to its coordinates $\VX_p = (X_p,Y_p)^t$ in the pseudo-Cartesian frame. In addition, we introduce ${\cal P}(c)$ the list of the vertices belonging to the cell $\Omega_c$ and ${\cal C}(p)$  the list of the cells  sharing the vertex $p$. These two sets are counterclockwise ordered. Let us introduce $p^-$ and $p^+$ the previous and the next nodes with respect to p in ${\cal P}(c)$. We denote by $L_{pc}^-,L_{pc}^+$ the half length of the edges $[pp^-], [pp^+]$. Similar notations are used for the normals outward $\VN_{pc}^+$ and $\VN_{pc}^-$. Finally, the corner normal $L_{pc} \VN_{pc}$ is given by $L_{pc}\VN_{pc} = L_{pc}^+\VN_{pc}^+ +L_{pc}^-\VN_{pc}^-$. All these notations have been displayed in \figref \ref{fig:6}.\\  

The first order spatial approximation of \eqref{eq:area1}-\eqref{eq:area3} is obtained considering local integrals  on each cell $\Omega_c$ rotated about the $X$-axis. The mass $m_c$ of the cell $\Omega_c$ is $m_c= \int_{\Omega_c} \rho dV$ and each flow variable $\phi$ (as total energy, velocity) is averaged over each cell through the formula
$$
\phi_c= \frac{1}{m_c} \int_{\Omega_c} \rho \phi dV,
$$
named cell-centered value. Then, we have 
\begin{align}
\displaystyle m_c \Dt \VU_c + \displaystyle    \overline{\CR}_c  \sum_{p \in {\cal P}(c)} \FPC = {\bf 0}, \label{eq:disc1}\\[0.25cm]
\displaystyle m_c \Dt E_c  + \displaystyle \sum_{p \in {\cal P}(c)}  \FPC \cdot \CR_p \VU_p = 0 \label{eq:disc2}.
\end{align}
In addition, the mesh is moved through the local kinematic equation given at each node by
\begin{equation}
\dfrac{d \VX_p}{dt} = \VU_p \text{ for } t>0 \quad \text{ and } \quad \VX_p(0) = \VX_p^0,
\label{eq:CGLdisc}
\end{equation}
with $\VU_p$ and $\VX_p^0$ respectively the velocity and the position of a node $p$ at initial time. In the previous equations, $\FPC$ is the numerical flux at each node $p$ of each cell $c$ defined by 
\begin{equation}
 \FPC = L_{pc} P_c \VN_{pc}-\CM_{pc}(\VU_p-\VU_c),
 \label{eq:FPC}
\end{equation}
with $\VU_p$ the velocity at the point $p$ and $P_c$ the mean value of the pressure in the cell $c$. The $2\times2$ matrices $\CM_{pc}$ and $\CM_{p}$
are defined as
\begin{equation}
 \CM_{pc} = \CZ_c \left( L^-_{pc}\VN^-_{pc} \otimes \VN_{pc}^- + L^+_{pc}\VN^+_{pc} \otimes \VN_{pc}^+ \right), \text{ and } \CM_p = \sum_{c\in {\cal C}(p)} \CM_{pc}.
\end{equation}
Where, we introduce the ``swept mass flux'' \cite{Dukowicz1} associated to the isentropic sound speed $a_c$ that is
\begin{equation}
\CZ_c = \rho_c a_c.
\end{equation}
This is nothing but the acoustic impedance. As it has been demonstrated in \cite{Maire4} the total energy and momentum conservation is equivalent to
\begin{equation}
 \sum_{c \in {\cal C}(p)} \FPC = {\bf 0}.
\label{eq:cons}
\end{equation}
 
Finally using \eqref{eq:FPC}, the nodal velocity $\VU_p$ is  deduced from \eqref{eq:cons} by solving the linear system 
\begin{equation}
 \CM_p \VU_p = \sum_{c \in {\cal C}(p)} (L_{pc} P_c \VN_{pc} +\CM_{pc} \VU_c).
\end{equation}

In \cite{Maire1}, the numerical fluxes used in the discretization of \eqref{eq:area1} and \eqref{eq:area3} are chosen for satisfying the local GCL constraint \eqref{eq:CGLdisc}. Here, we rather adopt a more simple approach that give similar results. Since \eqref{eq:CGLdisc} is explicitly solved for moving the mesh in time, there is no need to solve \eqref{eq:area1}. Thus, each cell volume $V_c$ is directly deduced from \eqref{eq:CGLdisc}. Thereby, it is possible to choose for the numerical flux in \eqref{eq:area3} a simple form as in \eqref{eq:disc2} with $\CR_p = 1-\alpha+\alpha Y_p$. Concerning, the momentum equation the mean value ${\CR}_c$ is equal to the discrete ratio $\overline{\CR}_c = \dfrac{V_c}{A_c}$.

Let us note that this new formulation of the area-weighted discretization relies on a node-centered solver which is exactly the same as the one developed in \cite{Maire2} for two-dimensional Cartesian geometry. However, the present spatial discretization does not satisfy rigorously the GCL compatibility requirement. In what follows, we will assess the discrepancy of our discretization to the GCL by analyzing the corresponding discrete divergence operator. The discrete divergence operator that corresponds to the present scheme writes as
\begin{equation}
\label{divp}
(\nabla \cdot \bs{U})_c=\frac{1}{V_c} \sum_{p \in \mathcal{P}(c)} \mathcal{R}_p (L_{pc}^{-}\bs{N}_{pc}^{-}+L_{pc}^{+}\bs{N}_{pc}^{+})\cdot \bs{U}_p,
\end{equation}
where $\mathcal{R}_p$ denotes the pseudo-radius of vertex $p$. It is shown in \cite{Maire4,Maire1} that the discrete divergence operator deduced from the GCL reads
\begin{equation}
\label{divgcl}
(\nabla \cdot \bs{U})_c^{\scriptstyle{GCL}}=\frac{1}{V_c} \sum_{p \in \mathcal{P}(c)} \frac{1}{3}[(2\mathcal{R}_p+\mathcal{R}_p^{-})L_{pc}^{-}\bs{N}_{pc}^{-}+(2\mathcal{R}_p+\mathcal{R}_p^{+})L_{pc}^{+}\bs{N}_{pc}^{+}]\cdot \bs{U}_p.
\end{equation}
If the time evolution of the position vector, $\bs{X}_{p}$, of vertex $p$ is governed by the trajectory equation [14], then one can prove that the time rate of change of the cell volume, $V_c$, satisfies
$$\frac{1}{V_c}\frac{d V_c}{d t}=(\nabla \cdot \bs{U})_c^{\scriptstyle{GCL}}.$$
Subtracting \eqref{divp} and \eqref{divgcl} leads to
\begin{equation}
\label{subtrac}
(\nabla \cdot \bs{U})_c-(\nabla \cdot \bs{U})_c^{\scriptstyle{GCL}}=\frac{1}{3V_c}\sum_{p \in \mathcal{P}(c)}[(\mathcal{R}_p-\mathcal{R}_p^{-})L_{pc}^{-}\bs{N}_{pc}^{-}+(\mathcal{R}_p-\mathcal{R}_p^{+})L_{pc}^{+}\bs{N}_{pc}^{+}]\cdot \bs{U}_p.
\end{equation}
Knowing that the summation in the previous equation is cyclic, shifting the index in the second term of the right hand-side yields
\begin{equation}
\label{subtrac2}
(\nabla \cdot \bs{U})_c-(\nabla \cdot \bs{U})_c^{\scriptstyle{GCL}}=\frac{1}{3V_c}\sum_{p \in \mathcal{P}(c)}[(\mathcal{R}_p^{+}-\mathcal{R}_p)L_{pc}^{+}\bs{N}_{pc}^{+}]\cdot (\bs{U}_{p^{+}}-\bs{U}_p).
\end{equation}
In case of a one-dimensional spherical flow on an equi-angular polar grid, the right-hand side of the previous equation is equal to zero. To prove this result, let us consider a quadrangular cell of an equi-angular polar grid. The proof proceeds in the following two steps:
\begin{itemize}
\item Either $p$ and $p^{+}$ are located on the same angular sector and thus the nodal velocity $\bs{U}_p$ and $\bs{U}_{p^{+}}$ are colinear to the direction of the angular sector which is orthogonal to the unit outward normal $\bs{N}_{pc}^{+}$. Hence, $(\bs{U}_{p^{+}}-\bs{U}_p)\cdot \bs{N}_{pc}^{+}=0$.
\item Or $p$ and $p^{+}$ are located on the same cercle of radius $R$, then the Cartesian components of their nodal velocities reads as
$$\bs{U}_{p}=U(R) 
\begin{pmatrix} 
\cos \theta \\
\sin \theta 
\end{pmatrix}
,
\quad
\bs{U}_{p^{+}}=U(R) 
\begin{pmatrix} 
\cos (\theta+\Delta \theta) \\
\sin (\theta +\Delta \theta)
\end{pmatrix}
.
$$
Here, $\theta$ denotes the angle of the angular sector, $U(r)$ is the module of the one-dimensional velocity field, and $\Delta \theta$ is size of the angular sector. A straigthforward computation shows that
$$\bs{U}_{p^{+}}-\bs{U}_p=2 U(R) \sin (\frac{\Delta \theta}{2}) \begin{pmatrix} 
-\sin (\theta+\frac{\Delta \theta}{2}) \\
\cos (\theta +\frac{\Delta \theta}{2})
\end{pmatrix}
.$$
Knowing that the unit outward normal is given by
$$\bs{N}_{pc}^{+}=\begin{pmatrix} 
\cos (\theta+\frac{\Delta \theta}{2}) \\
\sin (\theta +\frac{\Delta \theta}{2})
\end{pmatrix}
,
$$
we obtain that $(\bs{U}_{p^{+}}-\bs{U}_p)\cdot \bs{N}_{pc}^{+}=0$.
\end{itemize}
This ends the proof. This result shows that our new area-weighted discretization satisfies rigoroulsy the GCL compatibility requirement for one-dimensional spherical flows on equi-angular polar grids.

\section{MOF multi-material interface reconstruction phase in axisymmetric geometry}

The method used in this work to reconstruct interfaces, is the MOF approach well adapted for treating multi-materials interface problems \cite{Ahn1,Dyadechko1}.  Indeed, such a method enables to capture more accurately interfaces than the classical VOF strategy and allows the treatment of general multi-material flows (more than two materials) \cite{Galera1,Kucharik2}. This method has been recently extended in cylindrical geometries, for a single interface problem \cite{Anbarlooei1}. Here, extension to multi-material interface reconstruction phase to cylindrical coordinates is considered. 

\subsection{Moment of fluid method}

The main idea of MOF is to track each fluid in a cell using the zeroth and first moments \cite{Dyadechko1}.  
Given these two moments, interface is linearly reconstructed  insuring volume conservation. 
To this end, interface update is done minimizing the  discrepancy between the given moments and the reconstructed moments of the polygon behind the interface. 
One should note that no information from neighboring cells is required. This method is exact for linear interfaces and is second order accurate for smoothly curved ones.
In the context of multi-material configurations, one has to face to material ordering when reconstructing interface. The method presented here, allows to automatically determine the order of materials by constructing all the possible combination and choosing the sequence that leads to the configuration where the reconstructed moments are the closest to the given ones. 
The main difference  between cylindrical and planar geometry relies in the definition of the different moments. Since the interface reconstruction is done under volume conservative assumption, the zeroth moment $M_{k,c}^{0}$ of the $k$-th fluid in each cell $c$ is obviously given by
\begin{equation}
M_{k,c}^{0} =\int_{\Omega_{k,c}} \CR dA,
\end{equation} 
from this moment we can deduce the  volume fraction 
\begin{equation}
\alpha_{k,c} = \dfrac{M_{k,c}^{0}}{V_{c}} ,
\end{equation} 
with the cell volume $V_c = \int_{\Omega_{c}} \CR dA$. \\
Contrary to the  zeroth moment, the first moment can be defined without any specific requirement. Thus, it is possible to compute them in the two following different manners. In the one hand we can use the natural extension to axisymmetric geometries
\begin{equation}
M_{k,c}^{1} = \int_{\Omega_{k,c}}\CR \VX dA,
\end{equation} 
and from this moment we deduce the pseudo-centroid
\begin{equation}
\VX_{k,c} = \dfrac{M_{k,c}^{1}}{ V_{k,c}}, \text{ with } V_{k,c} = V_{c} \alpha_{k,c}.
\end{equation}  
This pseudo-centroid for a matter of simplicity will be called here the axisymmetric centroid.\\
On the other hand it can also be done with a planar definition as follows
\begin{equation}
M^{1,pl}_{k,c} =\int_{\Omega_{k,c}}\VX dA,  
\end{equation} 
and thus planar centroid will be obtain from
\begin{equation}
\VX^{pl}_{k,c} = \dfrac{M^{1,pl}_{k,c}}{A_{k,c}},  
\end{equation} 
where $A_{k,c}$ is the area of the $k$-th fluid in the cell $c$. 

Since this interface reconstruction method is coupled to our Lagrangian hydrodynamics scheme it requires to update the  volume fractions and material centroids. 
Using the equal strain assumption, the volume fractions do not evolve during the Lagrangian step (see \cite{Galera2} for more details).  However, the centroid locations are given from the Lagrangian step using a barycentric combination of the new positions of the mesh nodes as done in \cite{Galera1}. 

\subsection{Numerical validation}

\begin{figure}[H]
\centering
\begin{tabular}{cccc} 
{} &
 {Filament} & 
 {T-junction} & 
  {Y-junction} \\
\begin{sideways}{ \small  True partition ($\chi=1$)} \end{sideways} &
\begin{tikzpicture}[xscale=3.8,yscale=3.8]
\tikzstyle{every node}=[font=\scriptsize]
\def\R{1}

\coordinate (A) at (0.45,0.55);
\coordinate (B) at (0.65,0.35);

\coordinate(XI) at ($(A)+({-\R*cos(45)},{\R*sin(45)})$);
\coordinate(XII) at ($(B)+({\R*cos(30)},{-\R*sin(30)})$);
	
\clip (0,0) rectangle (1,1);
\draw[fill=green!20,draw=black] (0,0) rectangle (1,1);
\draw[fill=red!20, line width=0.25mm]  (XI) circle (1cm);
\draw[fill=blue!20,line width=0.25mm] (XII) circle (1cm);

\draw[draw=black,line width=0.5mm] (0,0) rectangle (1,1);

\node[black] at (A) {$\bullet$};
\node[above left] at (A) {$(0.45,0.55)$};
\node[above right] at (A) {~~~$45^{\circ}$};
\node[black] at (B) {$\bullet$};
\node[below left] at (B) {$(0.65,0.35)~$};
\node[above right] at (B) {~~~$60^{\circ}$};

\def\pentea{1}
\def\absa{0.1}
\coordinate (A1) at ($(0.25,\pentea*0.25+\absa)$);
\coordinate (A2) at ($(0.65,\pentea*0.65+\absa)$);
\draw[line width=0.25mm] (A1)--(A2);
\draw[dashed] (0.25,0.55)--(0.65,0.55);
\draw[-] ($(A)+(0.1,0)$) arc (0:45:0.1);

\coordinate (B1) at ($({0.5},{sin(60)/cos(60)*(0.5-0.65)+0.35})$);
\coordinate (B2) at ($({0.8},{sin(60)/cos(60)*(0.8-0.65)+0.35})$);
\draw[line width=0.25mm] (B1)--(B2);
\draw[dashed] (0.5,0.35)--(0.8,0.35);
\draw[-] ($(B)+(0.1,0)$)  arc (0:60:0.1);

\end{tikzpicture}&
\begin{tikzpicture}[xscale=3.8,yscale=3.8]
\tikzstyle{every node}=[font=\scriptsize]
\def\R{1}

     \coordinate (C1) at (0,0);
     \coordinate (C2) at (1,0);
     \coordinate (C3) at (0,1);
     \coordinate (C4) at (1,1);
     
\coordinate (A) at (0.5,0.5);

\coordinate(XI) at ($(A)+({\R},0)$);
\coordinate(XII) at ($(A)+(0,{-\R})$);

\clip (0,0) rectangle (1,1);
\draw[fill=green!20,draw=black,line width=0.25mm] (0,0) rectangle (1,1);
\draw[fill=red!20,draw=black,line width=0.25mm]  (XI) circle (\R );
\draw[fill=blue!20,draw=black,line width=0.25mm] (XII) circle (\R);

\draw[draw=black,line width=0.5mm] (0,0) rectangle (1,1);

\node[black] at (A) {$\bullet$};
\node[below] at (A) {$(0.5,0.5)$};
\node[above right] at (A) {$90^{\circ}$};
\node[above left] at (A) {$90^{\circ}$};


\coordinate (A1) at (0.3,0.5);
\coordinate (A2) at (0.7,0.5);
\coordinate (A3) at (0.5,0.7);
\draw[line width=0.25mm] (A)--(A1);
\draw[line width=0.25mm] (A)--(A2);
\draw[line width=0.25mm] (A)--(A3);
\draw[-] ($(A)+(0.18,0)$) arc (0:90:0.18);
\draw[-] ($(A)+(0,0.18)$) arc (90:180:0.18);

\end{tikzpicture}& 
\begin{tikzpicture}[xscale=3.8,yscale=3.8]
\tikzstyle{every node}=[font=\scriptsize]
\def\R{1}

\coordinate (A) at (0.5,0.5);

\coordinate(XI) at ($(A)+\R*({-cos(60)},{-sin(60)})$);
\coordinate(XII) at ($(A)+\R*({-cos(60)},{sin(60)})$);
\coordinate(XIII) at ($(A)+(\R,0)$);	

\clip (0,0) rectangle (1,1);
\draw[fill=blue!20,draw=black] (0,0) rectangle (1,1);
\draw[fill=blue!20,draw=black]  (XI) circle (1cm);
\draw[fill=green!20,draw=black] (XII) circle (1cm);
\draw[fill=red!20,draw=black] (XIII) circle (1cm);
\begin{scope}
\clip (0.5,0) rectangle (1,0.5);
\draw[fill=blue!20,draw=black]  (XI) circle (1cm);
\end{scope}
\draw[draw=black,line width=0.5mm] (0,0) rectangle (1,1);

\node[black] at (A) {$\bullet$};
\node[below] at (0.5,0.33) {$(0.5,0.5)$};

\def\pentea{1}
\def\absa{0.1}
\coordinate (A1) at (0.5,0.75);
\draw[line width=0.25mm] (A)--(A1);
\draw[-] (A) circle (0.15);

\coordinate (B1) at ($({0.25},{sin(30)/cos(30)*(0.25-0.5)+0.5})$);
\draw[line width=0.25mm] (A)--(B1);

\coordinate (C1) at ($({0.75},{sin(-30)/cos(-30)*(0.75-0.5)+0.5})$);
\draw[line width=0.25mm] (C1)--(A);

\node[above right,fill=red!20] at (0.52,0.5) {$120^{\circ}$};
\node[above left,,fill=green!20] at (0.48,0.5) {$120^{\circ}$};
\node[below,fill=blue!20] at (0.5,0.43) {$120^{\circ}$};

\end{tikzpicture}
\\
\begin{sideways}{~~~~\small Planar centroid } \end{sideways} &
\includegraphics[scale=0.21,clip,trim=5.3cm 2.3cm 4.5cm 0cm]{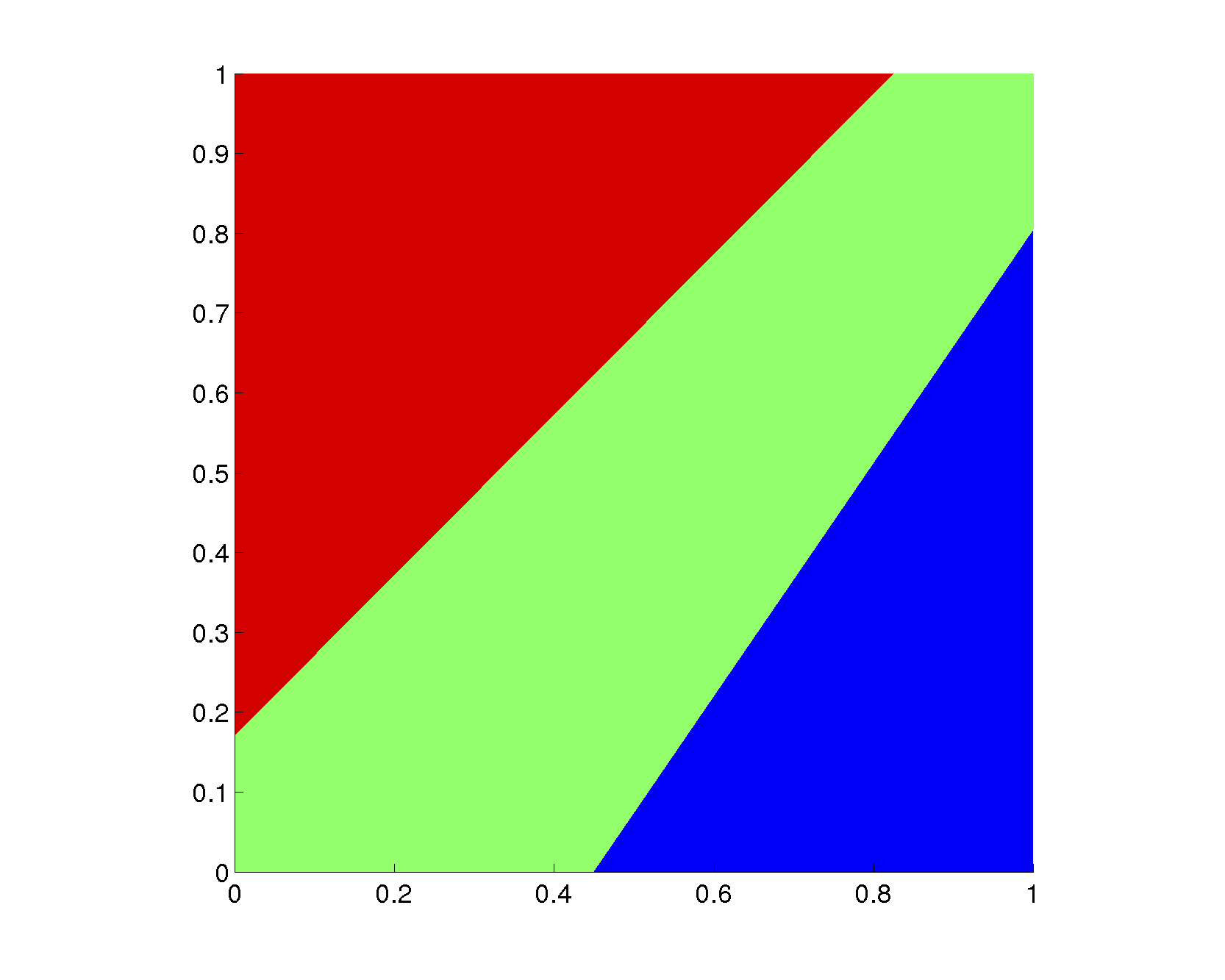} & 
\includegraphics[scale=0.21,clip,trim=5.3cm 2.3cm 4.5cm 0cm]{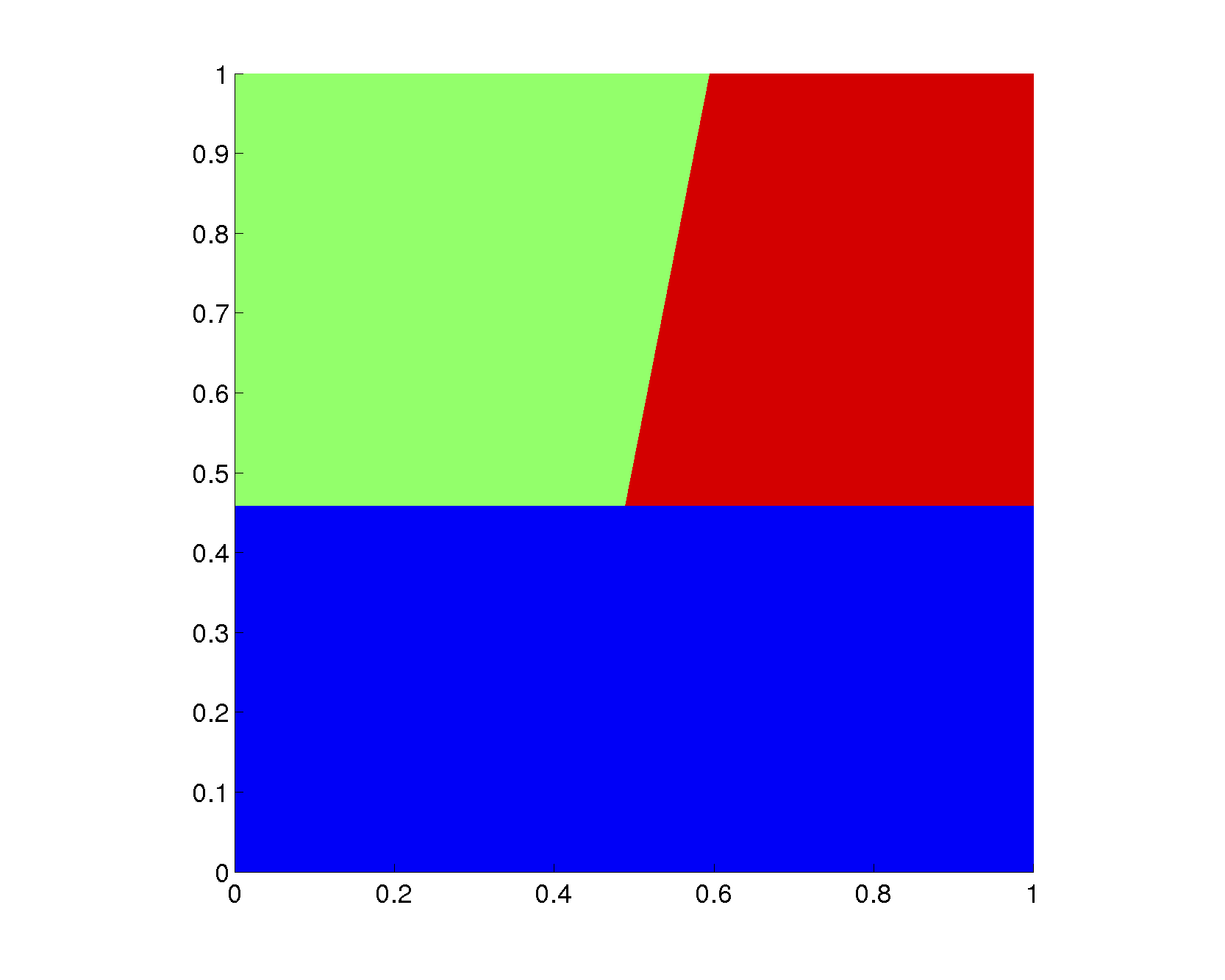}&
\includegraphics[scale=0.21,clip,trim=5.3cm 2.3cm 4.5cm 0cm]{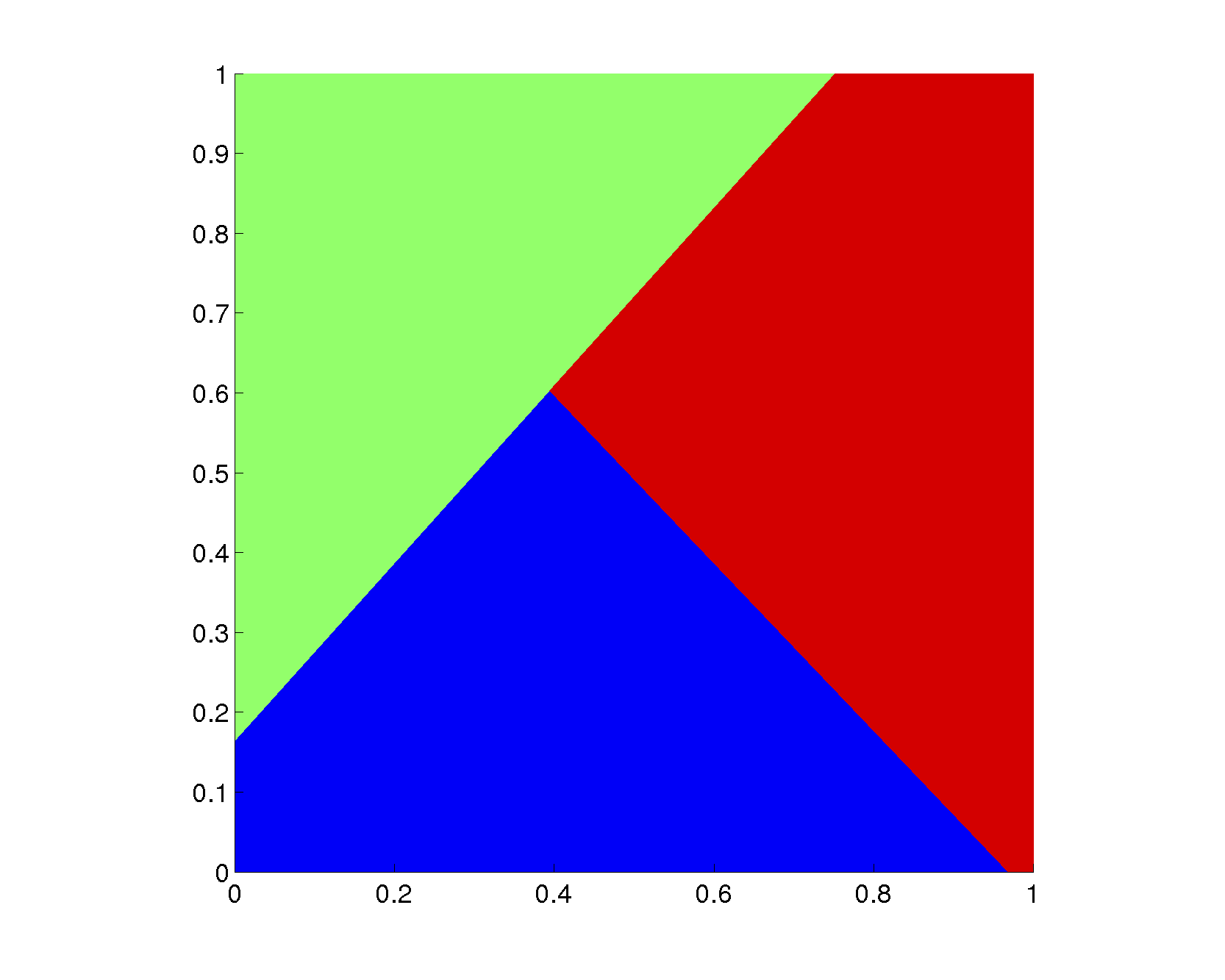} \\
\begin{sideways}{\small Axisymmetric centroid } \end{sideways} &
\includegraphics[scale=0.21,clip,trim=5.3cm 2.3cm 4.5cm 0cm]{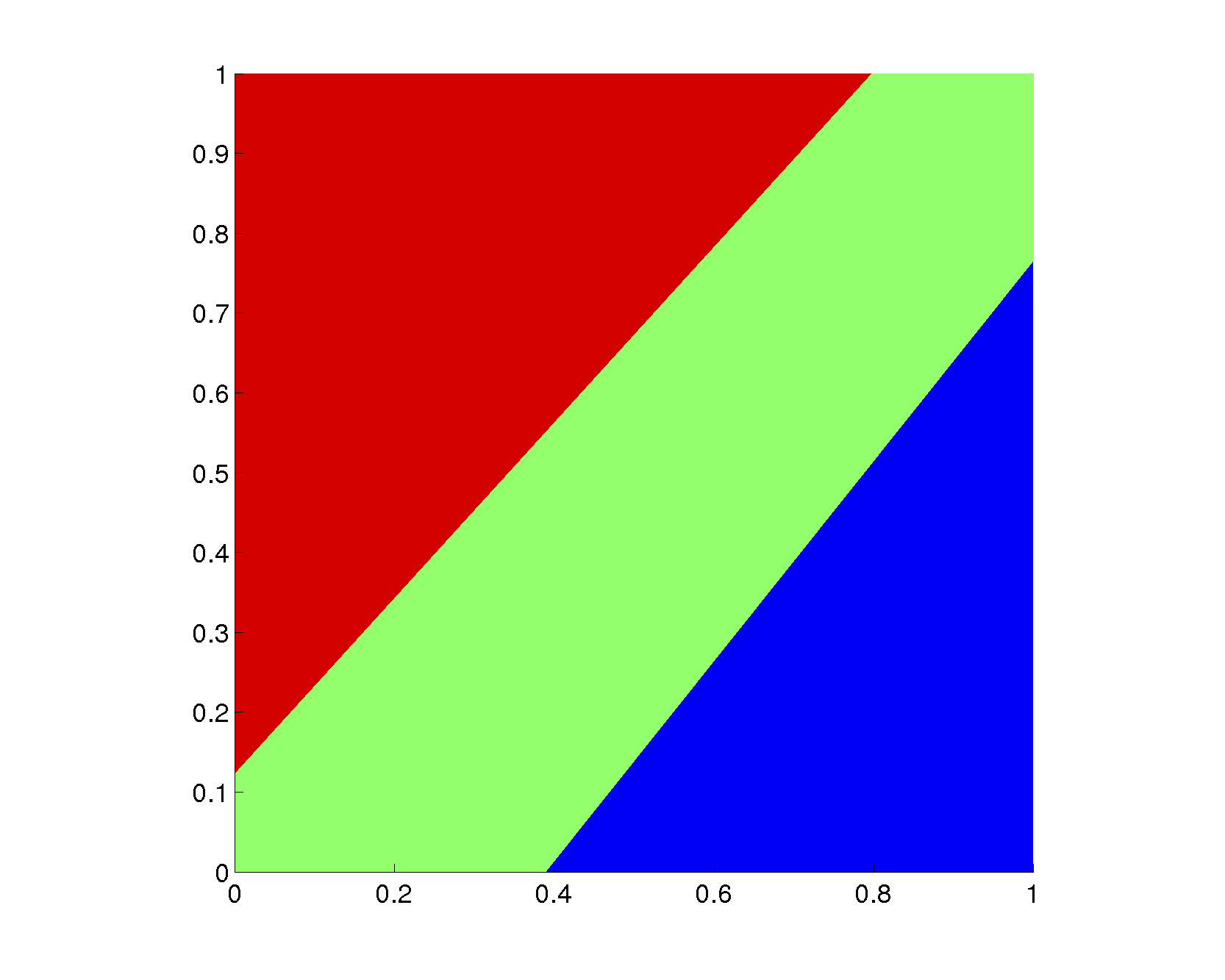} & 
\includegraphics[scale=0.21,clip,trim=5.3cm 2.3cm 4.5cm 0cm]{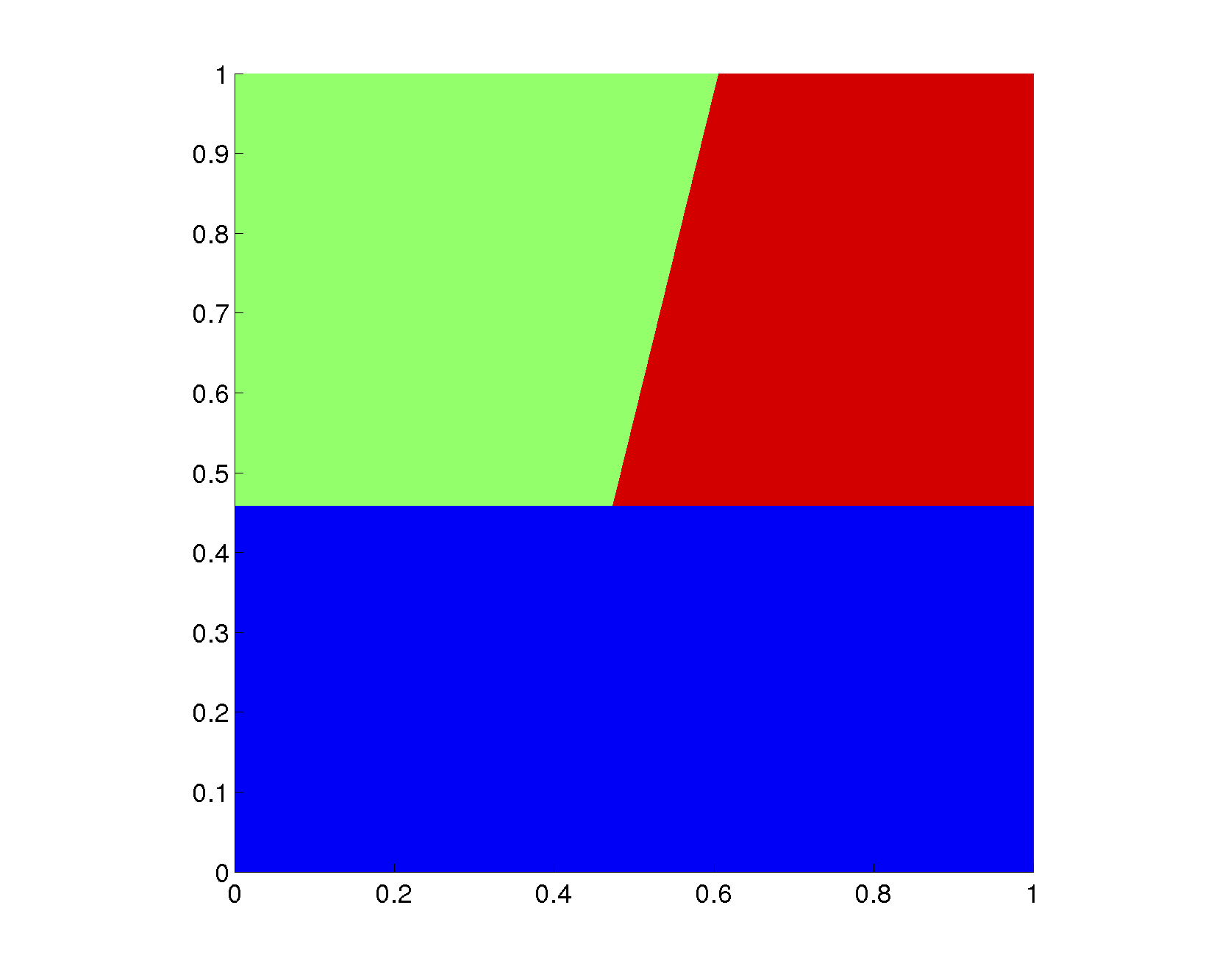}&
\includegraphics[scale=0.21,clip,trim=5.3cm 2.3cm 4.5cm 0cm]{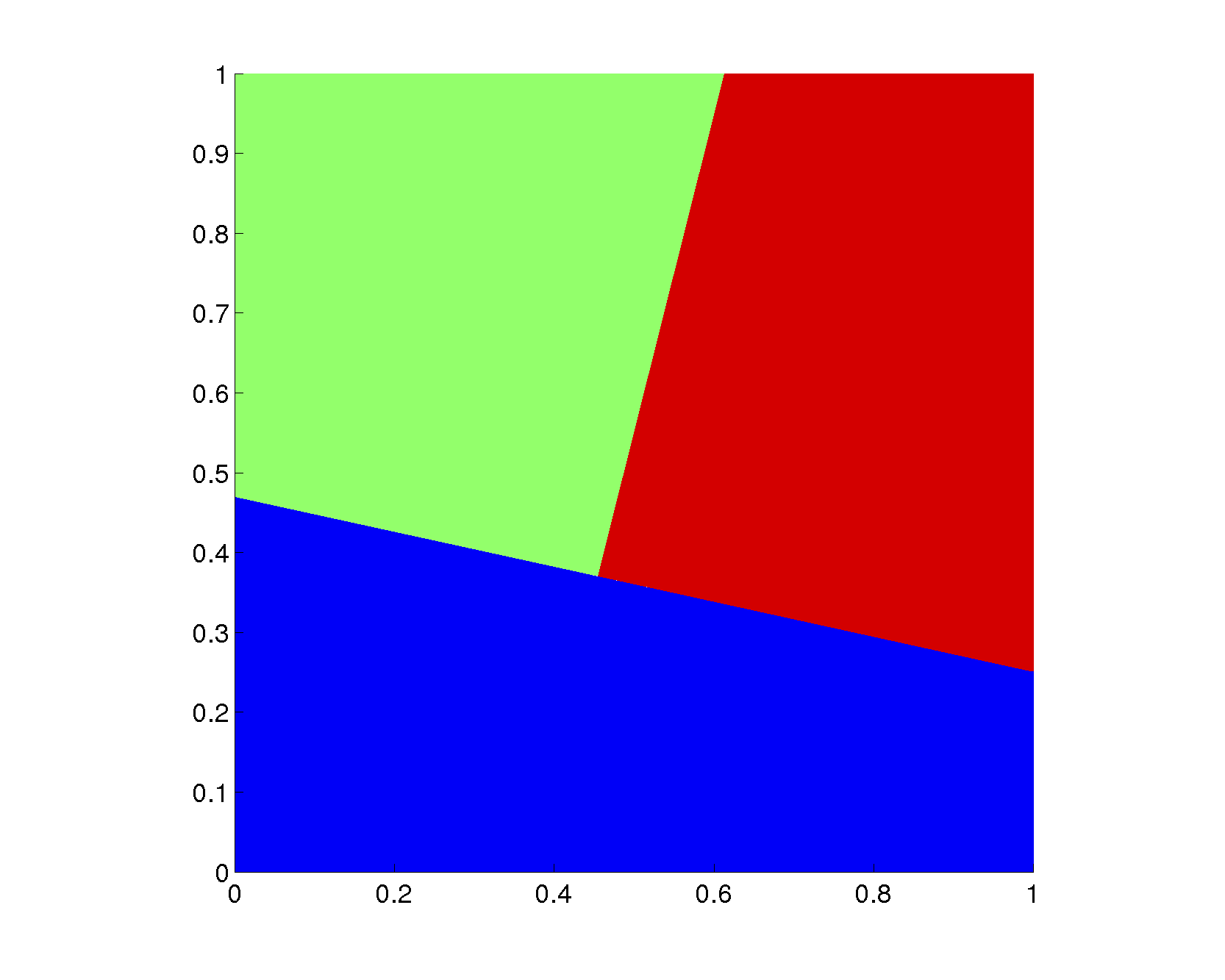} \\
\end{tabular}
\caption{MOF interface reconstruction test for three materials. From the top to the bottom: the true partitions for $\chi=1$ and their MOF reconstructions obtained with planar and axisymmetric centroids. From the left to the right: filament, T-junction and Y-junction configurations. \label{mof:1}}
\end{figure}
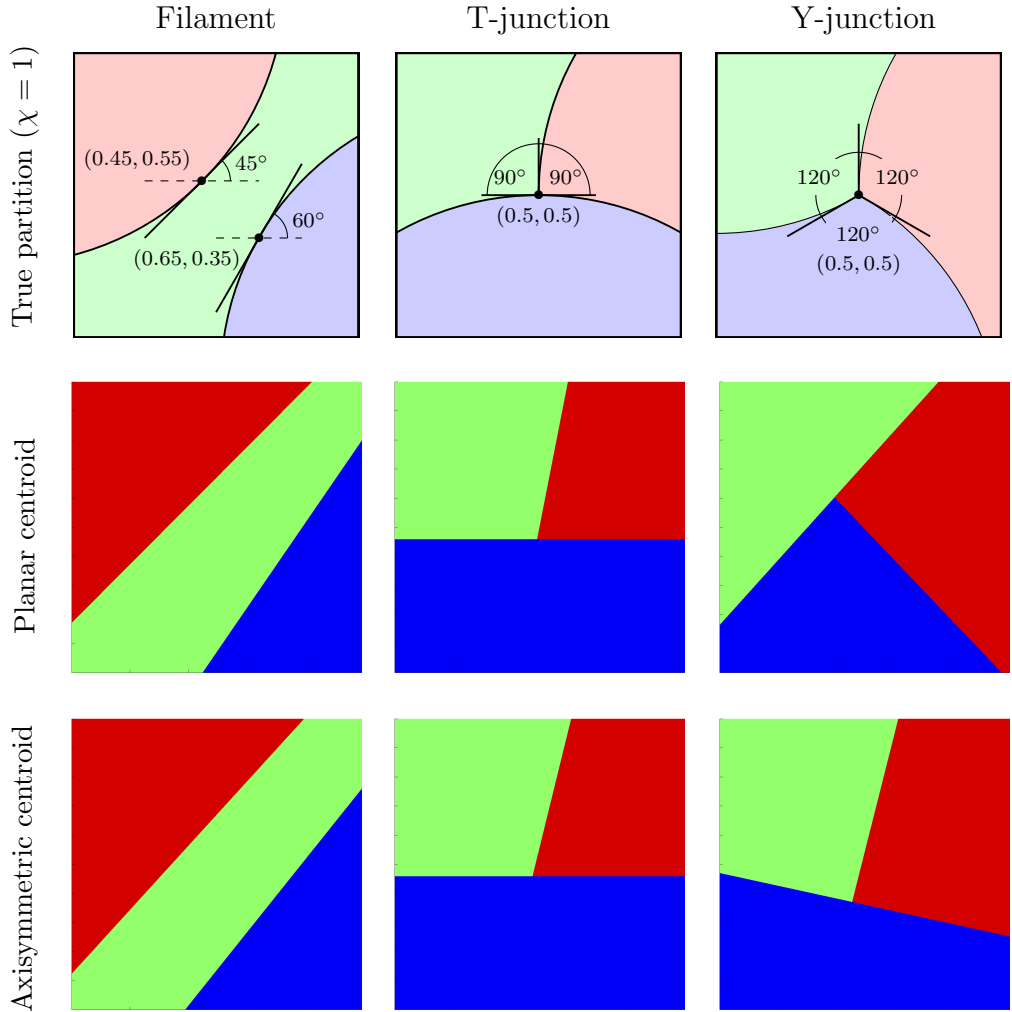

\begin{figure}[H]
\centering
\begin{tabular}{cccc} 
{} &
 {Filament} & 
 {T-junction} & 
  {Y-junction} \\
\begin{sideways}{\small  True partition ($\chi=64$)} \end{sideways} &
\begin{tikzpicture}[xscale=3.8,yscale=3.8]
\tikzstyle{every node}=[font=\scriptsize]
\def\R{1}

\coordinate (A) at (0.45,0.55);
\coordinate (B) at (0.65,0.35);
\def\pentea{1}
\def\absa{0.1}
\coordinate (A1) at ($(0.,\absa)$);
\coordinate (A2) at ($({\pentea*(1-0)-\absa},1)$);
\coordinate (A3) at ($(0,1)$);
\coordinate (B1) at ($({1},{sin(60)/cos(60)*(1-0.65)+0.35})$);
\coordinate (B2) at ($({-(0-0.35)*cos(60)/sin(60)+0.65},{0})$);
\coordinate (B3) at (1,0);

\clip (0,0) rectangle (1,1);
\draw[fill=green!20,draw=black] (0,0) rectangle (1,1);
\draw[fill=red!20, line width=0.25mm] (A1)--(A2)--(A3)--cycle; 
\draw[fill=blue!20, line width=0.25mm] (B1)--(B2)--(B3)--cycle; 
\draw[draw=black,line width=0.5mm] (0,0) rectangle (1,1);

\end{tikzpicture}&
\begin{tikzpicture}[xscale=3.8,yscale=3.8]
\tikzstyle{every node}=[font=\scriptsize]
\def\R{64}

     \coordinate (C1) at (0,0);
     \coordinate (C2) at (1,0);
     \coordinate (C3) at (0,1);
     \coordinate (C4) at (1,1);
     
\coordinate (A) at (0.5,0.5);

\coordinate(XI) at ($(A)+({\R},0)$);
\coordinate(XII) at ($(A)+(0,{-\R})$);

\clip (0,0) rectangle (1,1);
\draw[fill=green!20,draw=black] (0,0) rectangle (1,1);
\draw[fill=red!20,draw=black,line width=0.25mm]  (XI) circle (\R );
\draw[fill=blue!20,draw=black,line width=0.25mm] (XII) circle (\R);

\draw[draw=black,line width=0.5mm] (0,0) rectangle (1,1);

\end{tikzpicture}

%
%
%
&
\begin{tikzpicture}[xscale=3.8,yscale=3.8]
\tikzstyle{every node}=[font=\scriptsize]

\coordinate (A) at (0.5,0.5);
\coordinate (A1) at (0.5,1);
\coordinate (B1) at ($({0},{sin(30)/cos(30)*(0-0.5)+0.5})$);
\coordinate (C1) at ($({1},{sin(-30)/cos(-30)*(1-0.5)+0.5})$);

\clip (0,0) rectangle (1,1);
\draw[fill=green!20,draw=black] (0,0) rectangle (1,1);
\draw[draw=black,line width=0.5mm] (0,0) rectangle (1,1);
\draw[fill=blue!30,line width=0.25mm] (A)--(B1)--(0,0)--(1,0)--(C1)--cycle;
\draw[fill=red!30,line width=0.25mm] (A)--(C1)--(1,1)--(A1)--cycle;
\draw[draw=black,line width=0.5mm] (0,0) rectangle (1,1);
\end{tikzpicture}\\
\begin{sideways}{~~~~\small Planar centroid } \end{sideways} &
\includegraphics[scale=0.21,clip,trim=5.3cm 2.3cm 4.5cm 0cm]{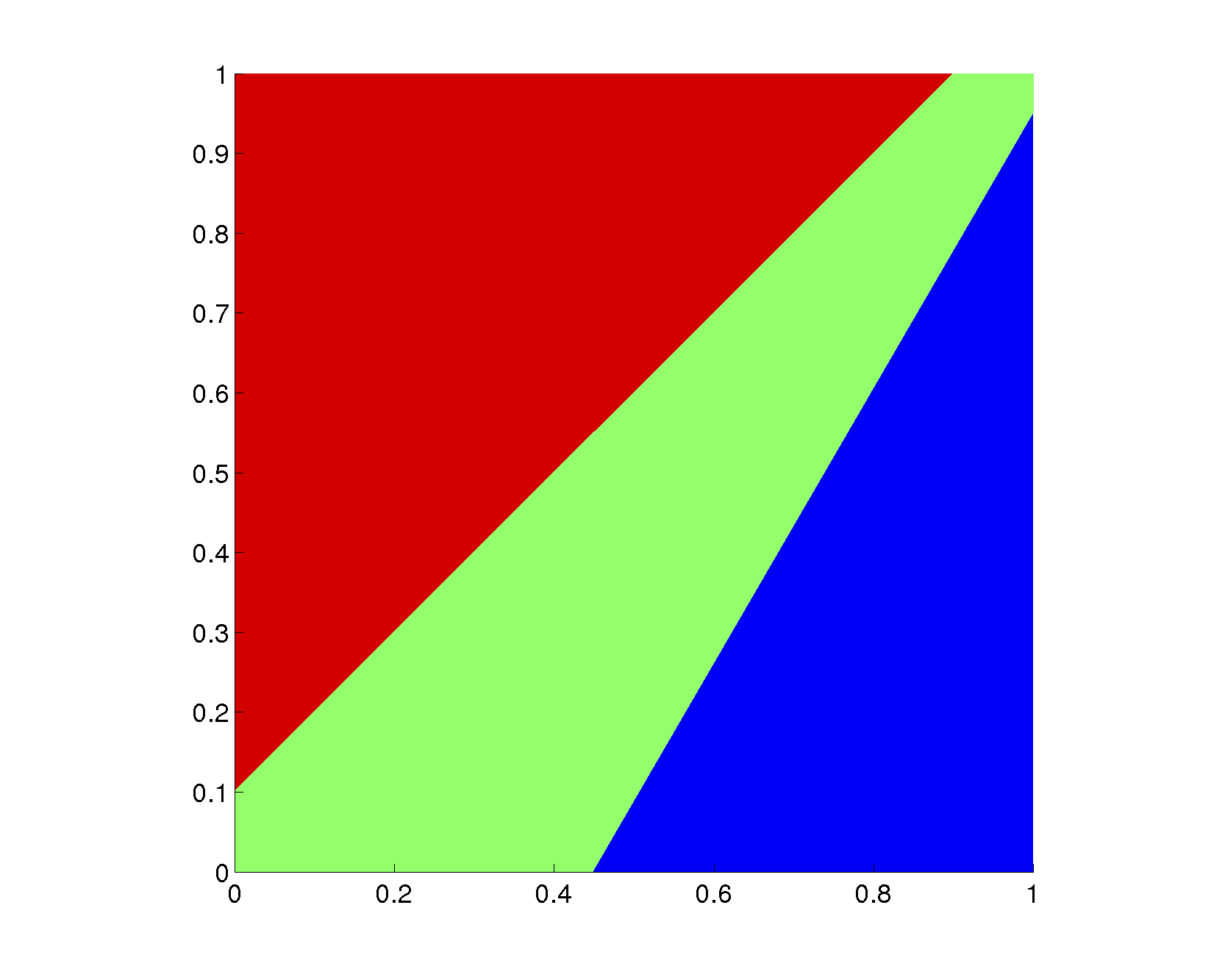} & 
\includegraphics[scale=0.21,clip,trim=5.3cm 2.3cm 4.5cm 0cm]{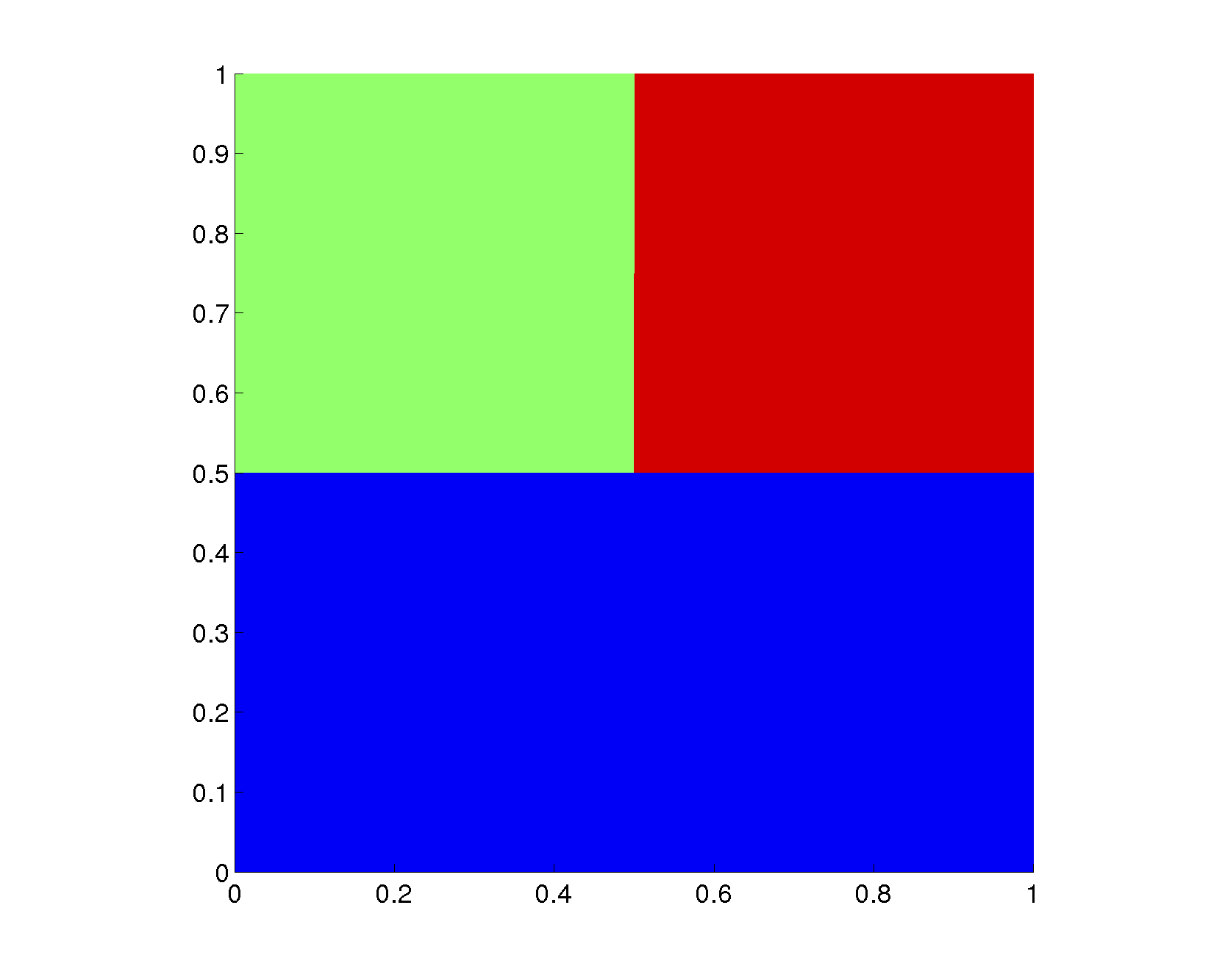}&
\includegraphics[scale=0.21,clip,trim=5.3cm 2.3cm 4.5cm 0cm]{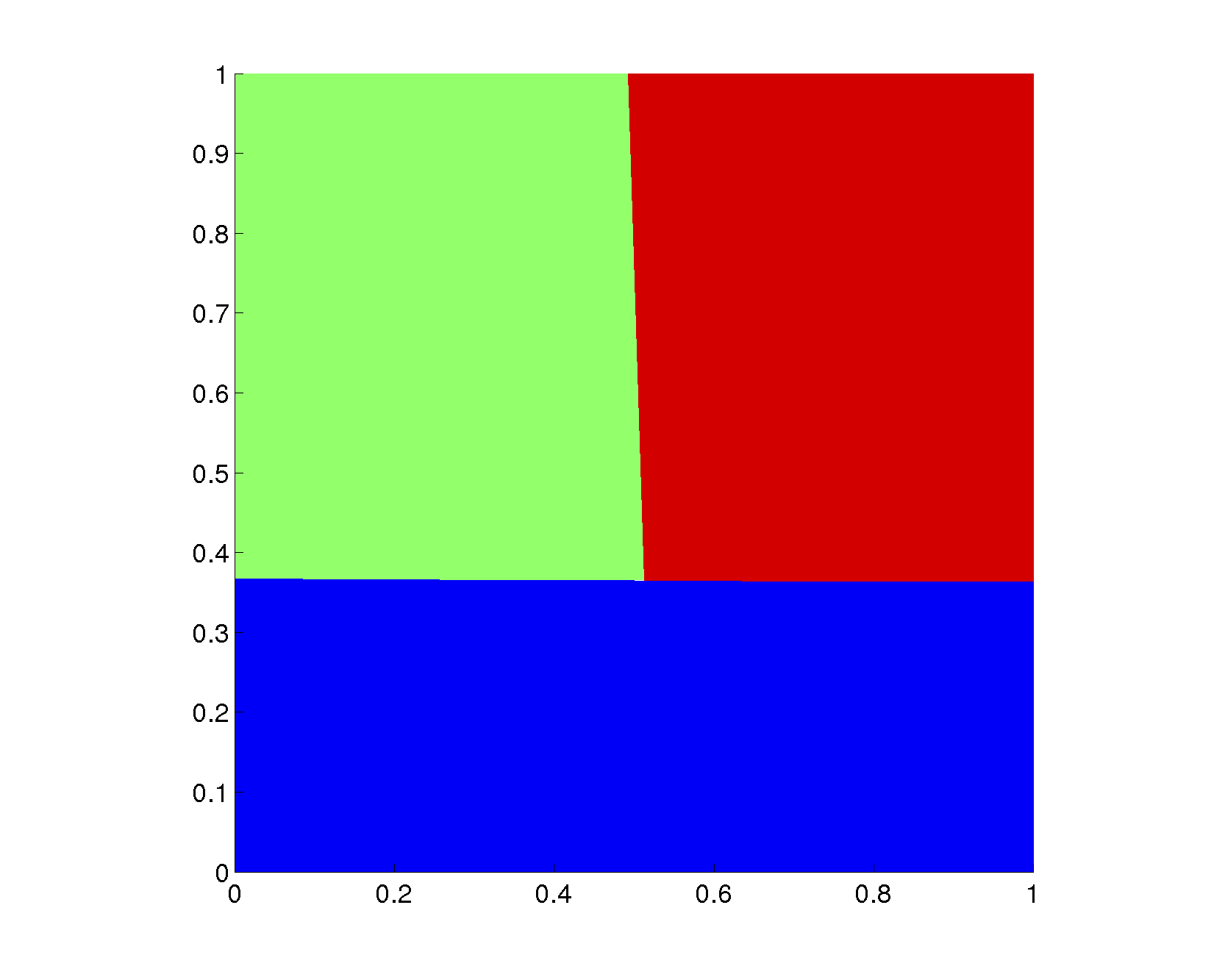} \\
\begin{sideways}{\small Axisymmetric centroid } \end{sideways} &
\includegraphics[scale=0.21,clip,trim=5.3cm 2.3cm 4.5cm 0cm]{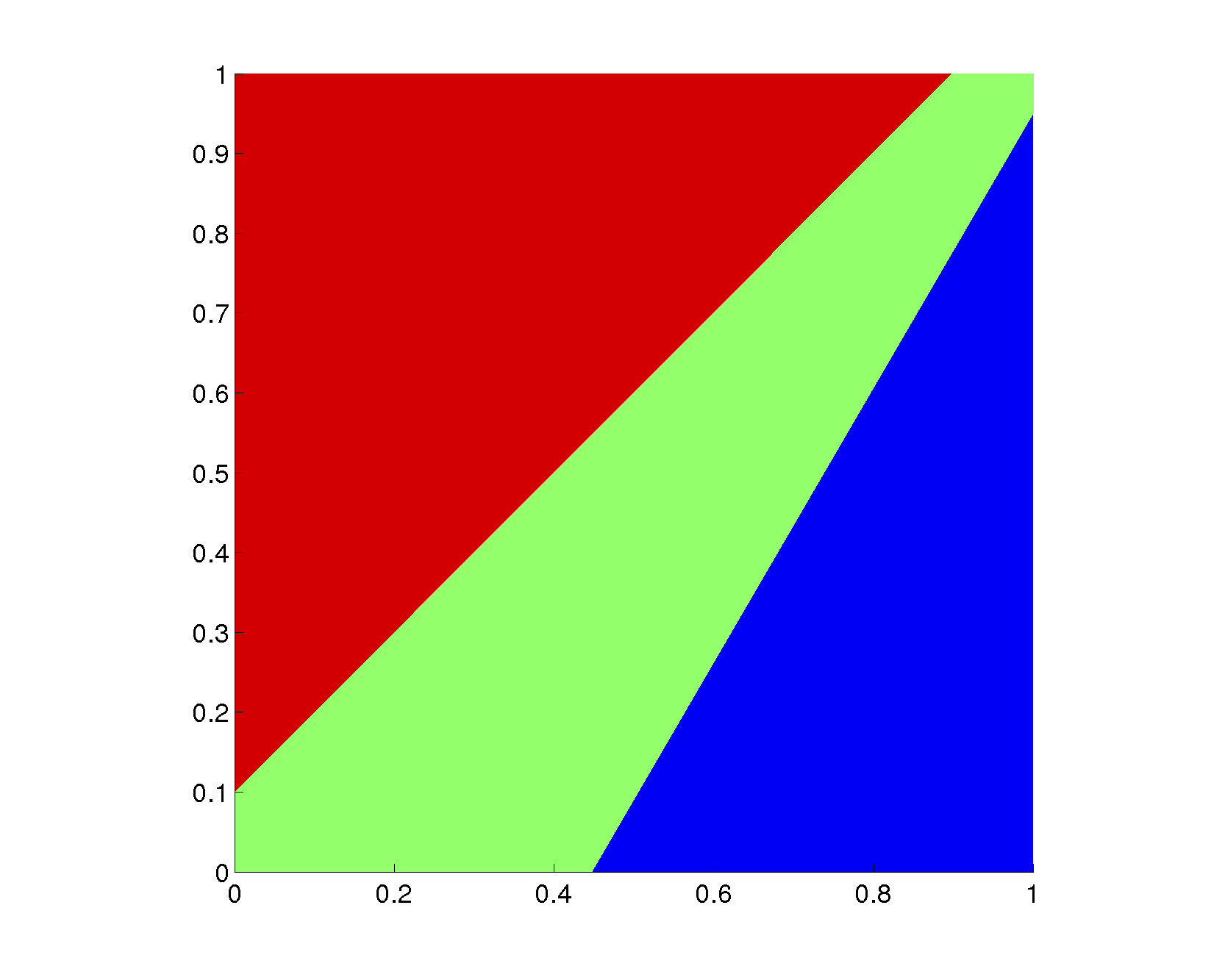} & 
\includegraphics[scale=0.21,clip,trim=5.3cm 2.3cm 4.5cm 0cm]{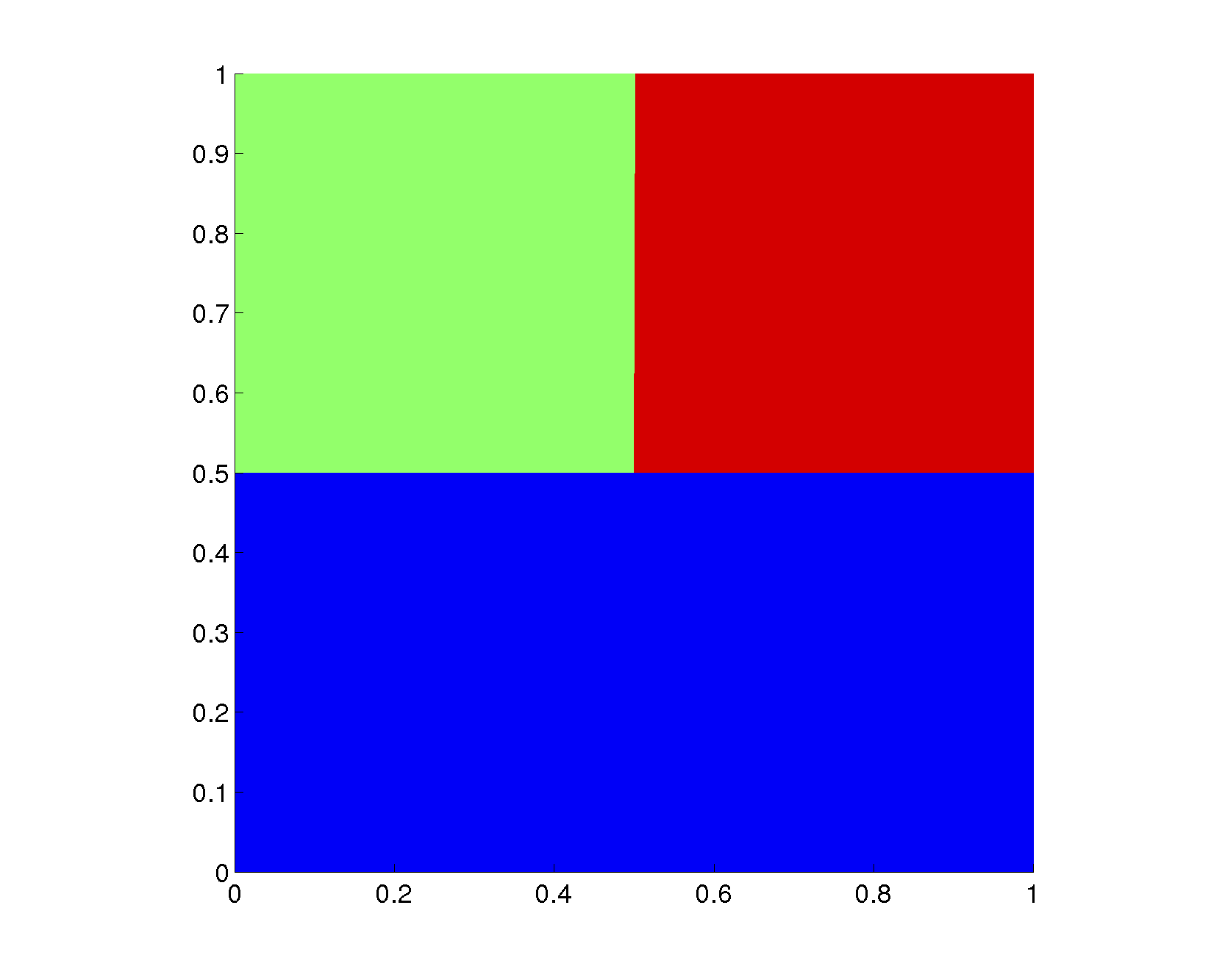}&
\includegraphics[scale=0.21,clip,trim=5.3cm 2.3cm 4.5cm 0cm]{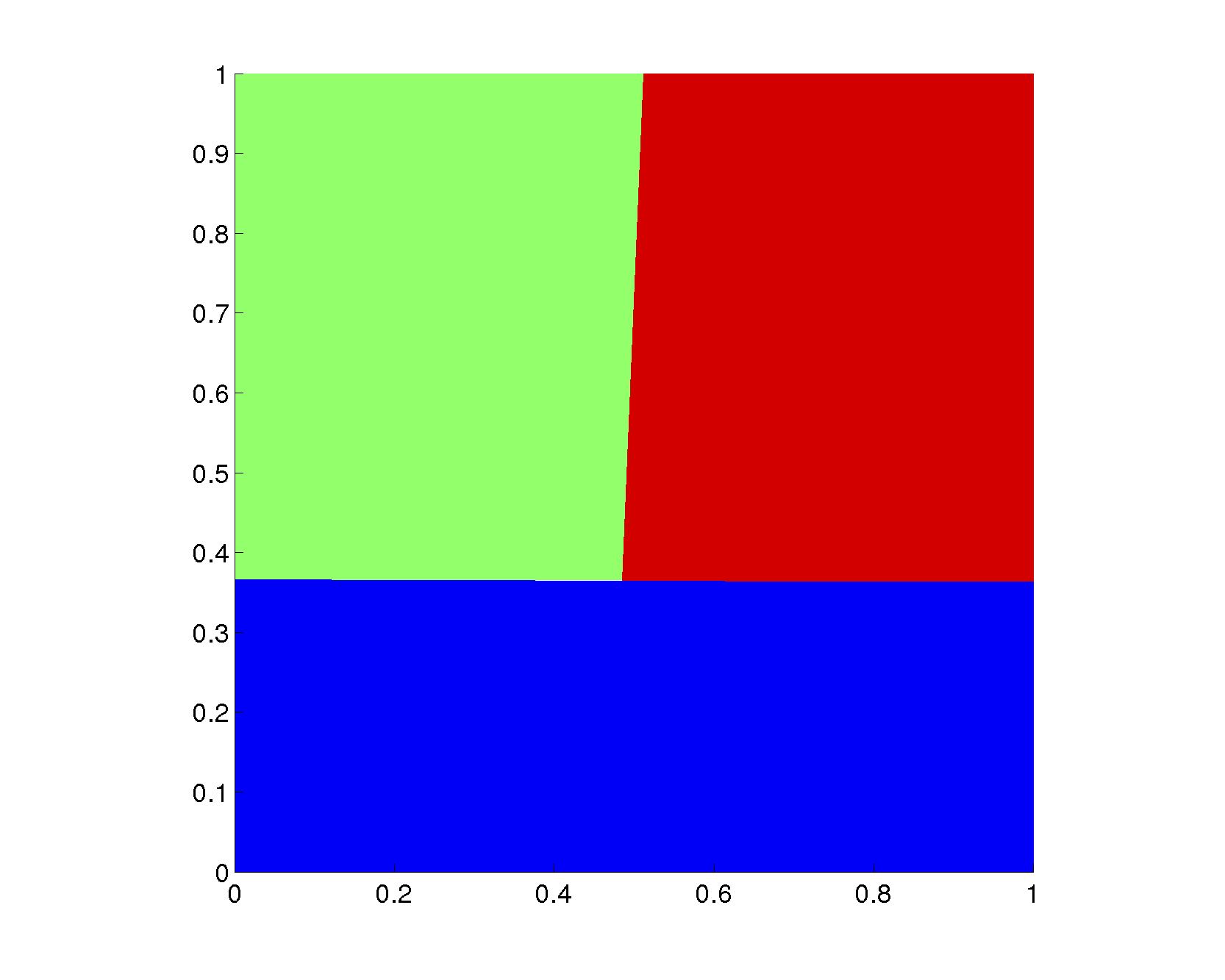} \\
\end{tabular}
\caption{MOF interface reconstruction test for three materials. From the top to the bottom: the true partitions for $\chi=64$ and their MOF reconstructions obtained with planar and axisymmetric centroids. From the left to the right: filament, T-junction and Y-junction configurations. \label{mof:2}}
\end{figure}
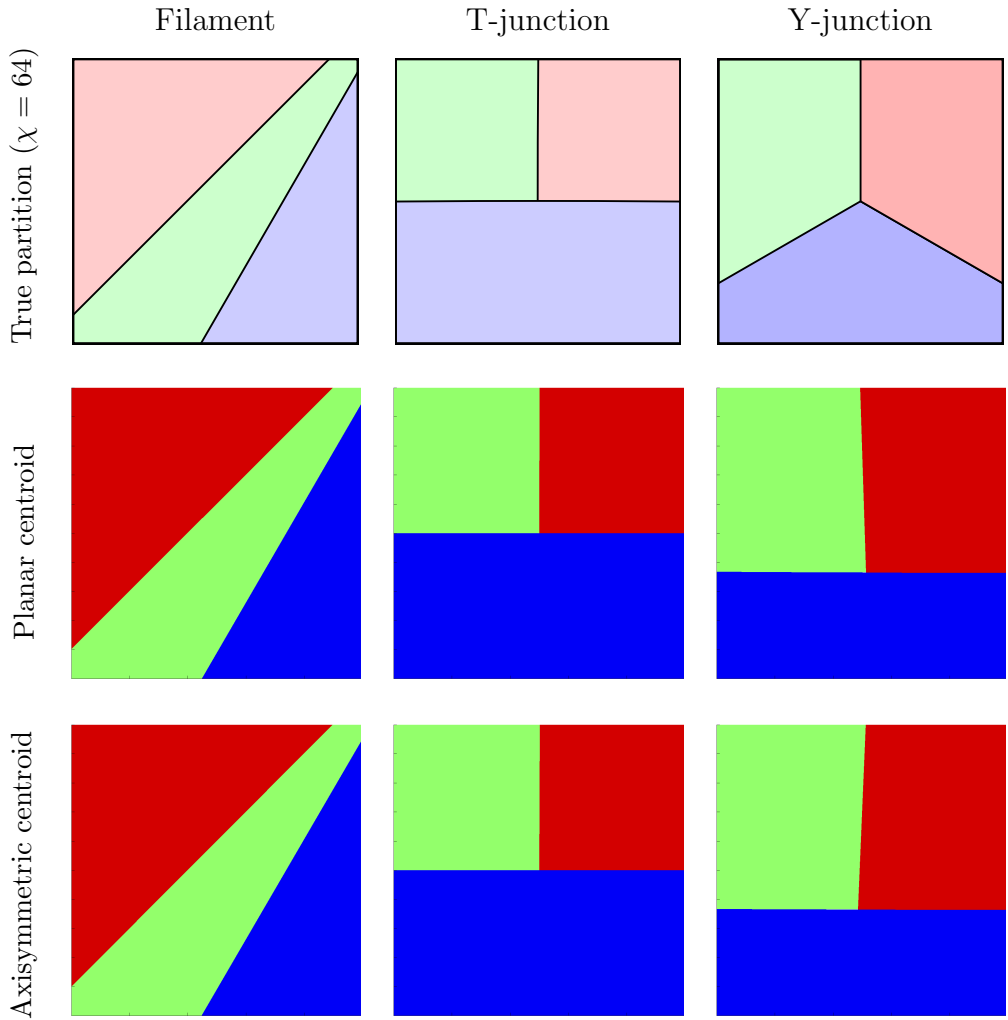

The main goal of  this section is to compare the results given by both axisymmetric and planar formulations of the centroids on several static test cases in one cell. As in \cite{Dyadechko1}, we consider three different mixed-cell layouts that are filament (without junction), T-junction and Y-junction.  The first two configurations correspond to $\mathcal{C}^2$-serial partitions whereas the third is not. In the considered test cases, the parameter $\chi$ corresponds to  the radius of the circles defining the interfaces. Two values are considered  with  $\chi=1$ and $\chi=64$.  In addition, the computation domain is reduce to the cell $[0;1]\times[0;1]$ (see figures \figref \ref{mof:1} and  \figref  \ref{mof:2}).\\

 In the first case, with $\chi=1$, we notice small differences for the filament case, no notable difference on the T-Junction but the Y-junction results for axisymmetric and planar formulations  present distinct interface positions due to a different ordering of the materials. For a large radius $\chi=64$, the curves are reduced to piecewise linear interfaces. Then, the result using both formulations are very close to each other. For the two first cases filament and T-junction, the results are exact. Regarding the Y-junction, it remains a good approximation. These results illustrate the capability of both planar and axisymmetric centroid formulation for MOF to treat accurately multi-material problem. Nevertheless, for consistency with the global cylindrical coordinate formulation, the axisymmetric formulation for the centroids is retained in the sequel.

\section{Rezoning phase improvement for polar meshes}

The rezoning phase introduced in \cite{Galera1,Galera2} consists in moving the Lagrangian grid to improve its geometric quality. The objective of this part is to extend this approach to polar meshes. To this end, the proposed procedure relies on two main steps. The first phase is dedicated to compute the smoothed grid from the Lagrangian one through CNS method. Then the final mesh is deduced from the smoothed one by a relaxation procedure to keep the rezoned grid as close as possible to  the Lagrangian grid in order to insure computation accuracy and avoid unphysical mesh rezoning. In the sequel one should note that rezoning is formulated only for planar geometry in the frame $\{0,X,Y\}$.\\

For the sake of readability, in the rest of the paper the quantities without any accent $a$ are associated to Lagrangian mesh. After the rezoning step we use $a^{rez}$, and finally after relaxation the quantities related to the rezoned mesh are noted with the tilde accent $\tilde{a}$.

\subsection{General condition number smoothing (GCNS)}

As it is pointed out in the introduction, CNS approach is well adapted to rezone Cartesian meshes but it still suffers from drawbacks for polar ones. Indeed, in this case the mesh seems to collapse (like an implosion) to the origin. To circumvent this difficulty, it has been proposed to modify the CNS algorithm using specific weight associated to the mesh geometry  \cite{Vachal1} that controls mesh rezoning with regards to the radius for example. Nevertheless, this approach is  not completely satisfactory. First, it strongly depends on the choice of the weight, that may affect the quality of the mesh which can be shifted in the opposite direction to the origin for example. Furthermore, there is still a residual compression near the origin due to singularity at this point. In conclusion, it does not preserve a uniform polar mesh. For this reason, a different strategy is presented here. 
The main idea developed here is to apply the CNS rezoning algorithm in $(r,\theta)$-coordinate system. 
In fact, a polar mesh initially expressed using a Cartesian coordinates $(X,Y)$ leads to a structured Cartesian mesh in $(r,\theta)$-coordinates.  
Here, a general presentation of the algorithm is made for unstructured meshes. \\
 
Assuming that the resulting mesh from the Lagrangian phase is unfolded (otherwise untangling procedure is used to correct invalid cells \cite{Vachal2}) . Thus, the proposed algorithm consists for polar meshes in three different steps as depicted on \figref\ref{fig:1}. For the sake of simplicity, we consider in the sequel only the case of Cartesian and polar structured meshes. 

\begin{figure}[h!]
 \centering
 
 \begin{tikzpicture}[x=0.2cm,y=0.2cm]
\draw (0,0) -- (0,10);
\draw (0,0) -- (10,0);
\draw (0,0) -- (9.5,3.1);
\draw (0,0) -- (8.,5.9);
\draw (0,0) -- (5.9,8.);
\draw (0,0) -- (3.1,9.5);

\draw (10,0.) -- (9.5,3.1);
\draw (9.5,3.1) -- (8.,5.9);
\draw (8.,5.9) -- (5.9,8.);
\draw (5.9,8.) -- (3.1,9.5);
\draw (3.1,9.5)-- (0,10);

\draw (1,0.) -- (0.95,0.31);
\draw (0.95,0.31) -- (0.8,0.59);
\draw (0.8,0.59) -- (0.59,0.8);
\draw (0.59,0.8) -- (0.31,0.95);
\draw (0.31,0.95)-- (0,1);

\draw (2,0.) -- (1.9,0.62);
\draw (1.9,0.62) -- (1.6,1.2);
\draw (1.6,1.2) -- (1.2,1.6);
\draw (1.2,1.6) -- (0.6,1.9);
\draw (0.6,1.9)-- (0,2);

\draw (4,0.) -- (3.8,1.24);
\draw (3.8,1.24) -- (3.2,2.4);
\draw (3.2,2.4) -- (2.4,3.2);
\draw (2.4,3.2) -- (1.2,3.8);
\draw (1.2,3.8)-- (0,4);

\draw (7,0.) -- (6.7,2.1);
\draw (6.7,2.1) -- (5.6,4.2);
\draw (5.6,4.2) -- (4.2,5.6);
\draw (4.2,5.6) -- (2.1,6.7);
\draw (2.1,6.7)-- (0,7);


\node[below] at (5,0){\footnotesize $(X^{lag},X^{lag})$};
\end{tikzpicture}
\hspace{-0.25cm}
 \begin{tikzpicture}[x=0.2cm,y=0.2cm]
	\tikzstyle{fleche}=[->,>=latex,line width=0.5mm,draw=black,fill=black]
	\node (D) at ( 5,7.5){};
	\node (G) at (0,7.5){};
	\node (0) at (0,0){};
	\node (texte) at (0,4){};
	\draw[fleche] (G)--(D);
	\node[above] at (2.5,7.5) {$\DT$};
\end{tikzpicture} 
\hspace{-0.25cm}
 \begin{tikzpicture}[x=0.2cm,y=0.2cm]
\draw (0,0) -- (0,10);
\draw (1,0) -- (1,10);
\draw (2,0) -- (2,10);
\draw (4,0) -- (4,10);
\draw (7,0) -- (7,10);
\draw (10,0) -- (10,10);
\draw (0,0) -- (10,0);
\draw (0,2) -- (10,2);
\draw (0,4) -- (10,4);
\draw (0,6) -- (10,6);
\draw (0,8) -- (10,8);
\draw (0,10) -- (10,10);
\node[below] at (5,0){\footnotesize $(r^{lag},\theta^{lag})$};
\end{tikzpicture}
\hspace{-0.25cm}
 \begin{tikzpicture}[x=0.2cm,y=0.2cm]
	\tikzstyle{fleche}=[->,>=latex,line width=0.5mm,draw=black,fill=black]
	\node (D) at ( 10,7.5){};
	\node (G) at (0,7.5){};
	\node (0) at (0,0){};
	\node (texte) at (0,4){};
	\draw[fleche] (G)--(D) ;
	\node[above] at (5,7.5) {\small CNS};
	\node[below] at (5,7.5) {\small rezoning};
\end{tikzpicture} 
\hspace{-0.25cm}
\begin{tikzpicture}[x=0.2cm,y=0.2cm]
\draw (0,0) -- (0,10);
\draw (2,0) -- (2,10);
\draw (4,0) -- (4,10);
\draw (6,0) -- (6,10);
\draw (8,0) -- (8,10);
\draw (10,0) -- (10,10);
\draw (0,0) -- (10,0);
\draw (0,2) -- (10,2);
\draw (0,4) -- (10,4);
\draw (0,6) -- (10,6);
\draw (0,8) -- (10,8);
\draw (0,10) -- (10,10);
\node[below] at (5,0){\footnotesize $(r^{rez},\theta^{rez})$};
\end{tikzpicture}
\hspace{-0.25cm}
 \begin{tikzpicture}[x=0.2cm,y=0.2cm]
	\tikzstyle{fleche}=[->,>=latex,line width=0.5mm,draw=black,fill=black]
	\node (D) at ( 5,7.5){};
	\node (G) at (0,7.5){};
	\node (0) at (0,0){};
	\node (texte) at (0,4){};
	\draw[fleche] (G)--(D);
	\node[above] at (2.5,7.5) {$\DT^{-1}$};
\end{tikzpicture} 
\hspace{-0.25cm}
 \begin{tikzpicture}[x=0.2cm,y=0.2cm]
\draw (0,0) -- (0,10);
\draw (0,0) -- (10,0);
\draw (0,0) -- (9.5,3.1);
\draw (0,0) -- (8.,5.9);
\draw (0,0) -- (5.9,8.);
\draw (0,0) -- (3.1,9.5);

\draw (10,0.) -- (9.5,3.1);
\draw (9.5,3.1) -- (8.,5.9);
\draw (8.,5.9) -- (5.9,8.);
\draw (5.9,8.) -- (3.1,9.5);
\draw (3.1,9.5)-- (0,10);

\draw (2,0.) -- (1.9,0.62);
\draw (1.9,0.62) -- (1.6,1.2);
\draw (1.6,1.2) -- (1.2,1.6);
\draw (1.2,1.6) -- (0.6,1.9);
\draw (0.6,1.9)-- (0,2);

\draw (4,0.) -- (3.8,1.24);
\draw (3.8,1.24) -- (3.2,2.4);
\draw (3.2,2.4) -- (2.4,3.2);
\draw (2.4,3.2) -- (1.2,3.8);
\draw (1.2,3.8)-- (0,4);

\draw (8,0.) -- (7.6,2.48);
\draw (7.6,2.48) -- (6.4,4.8);
\draw (6.4,4.8) -- (4.8,6.4);
\draw (4.8,6.4) -- (2.4,7.6);
\draw (2.4,7.6)-- (0,8);

\draw (6,0.) -- (5.7,1.8);
\draw (5.7,1.8) -- (4.8,3.6);
\draw (4.8,3.6) -- (3.6,4.8);
\draw (3.6,4.8) -- (1.8,5.7);
\draw (1.8,5.7)-- (0,6);

\node[below] at (5,0){\footnotesize $(X^{rez},Y^{rez})$};
\end{tikzpicture}
\caption{Rule representation for GCNS algorithm for a polar mesh. \label{fig:1}}
\end{figure}
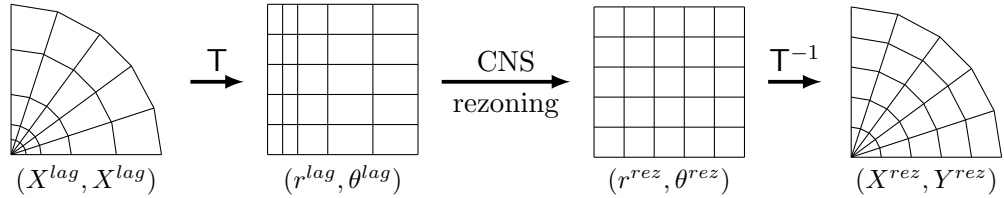
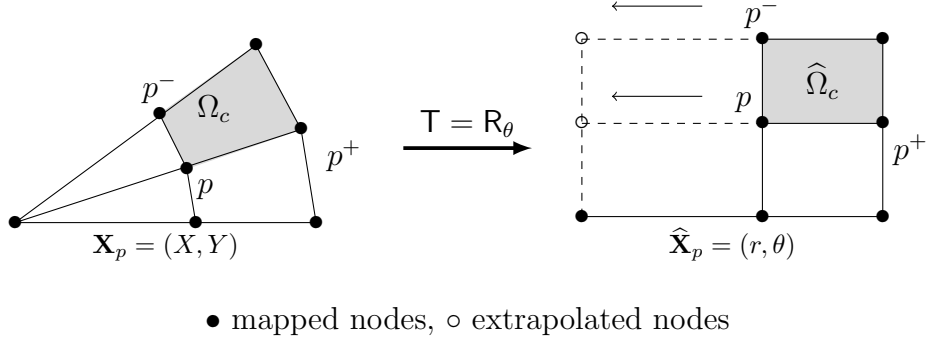
\begin{figure}[h!]
 \centering

\begin{tikzpicture}[x=0.4cm,y=0.4cm]
\filldraw[gray!30] (4.8,3.6)--(5.7,1.8)--(9.5,3.1)--(8.,5.9)--cycle node[black] at (6.6,3.8) {$\Omega_c$};
\draw (0,0) -- (10,0);
\draw (0,0) -- (9.5,3.1);
\draw (0,0) -- (8.,5.9);
\draw (10,0.) -- (9.5,3.1);
\draw (9.5,3.1) -- (8.,5.9);
\draw (6,0.) -- (5.7,1.8);
\draw (5.7,1.8) -- (4.8,3.6);
\node[below] at (5,0){\footnotesize $\Xp=(X,Y)$};
\node at (5.7,1.8)  {$\bullet$};
\node at (0,0)  {$\bullet$};
\node at (8.,5.9)  {$\bullet$};
\node at (9.5,3.1)  {$\bullet$};
\node at (6,0.)  {$\bullet$};
\node at (4.8,3.6)  {$\bullet$};
\node at (10,0)  {$\bullet$};
\node[below right] at (5.7,1.8)   {$p$};
\node[below right] at (10,3.1)   {$p^{+}$};
\node[above] at  (4.8,3.6)   {$p^{-}$};
\end{tikzpicture}
 \begin{tikzpicture}[x=0.4cm,y=0.2cm]
\tikzstyle{fleche}=[->,>=latex,line width=0.5mm,draw=black,fill=black]
\node (D) at ( 5,7.5){};
\node (G) at (0,7.5){};
\node (0) at (0,0){};
\node (texte) at (0,4){};
\draw[fleche] (G)--(D);
\node[above] at (2.5,7.5) {$\DT = \DR_{\theta}$};
\end{tikzpicture} 
\begin{tikzpicture}[x=0.4cm,y=0.4cm]
\filldraw[gray!30] (6,3.1) rectangle (10,5.9) node[black] at (8,4.5) {$\widehat{\Omega}_c$};
\draw (0,0) -- (10,0);
\draw[dashed] (0,0) -- (0,5.9);
\draw (6,0) -- (6,5.9);
\draw (10,0) -- (10,5.9);
\draw[dashed] (0.,5.9) -- (6,5.9);
\draw[<-](1,7) -- (4,7);
\draw (6,5.9) -- (10,5.9);
\draw[dashed] (0.,3.1) -- (6,3.1);
\draw[<-] (1,4) -- (4,4);
\draw (6.,3.1) -- (10,3.1);
\node at (6,3.1)  {$\bullet$};
\node at (0,0)    {$\bullet$};
\node at (0,5.9)  {$\circ$};
\node at (0,3.1)  {$\circ$};
\node at (10.,5.9){$\bullet$};
\node at (10,3.1) {$\bullet$};
\node at (6,0.)   {$\bullet$};
\node at (6,5.9)  {$\bullet$};
\node at (10,0)   {$\bullet$};
\node[above left] at (6,3.1)   {$p$};
\node[below right] at (10,3.1)   {$p^{+}$};
\node[above] at  (6,5.9)   {$p^{-}$};
\node[below] at (5,0){\footnotesize $\Xtp=(r,\theta)$};
\end{tikzpicture}
\\$\bullet$ mapped nodes, $\circ$ extrapolated nodes
\caption{Notations and mapping between cartesian and polar coordinates. \label{fig:GCNS}}
\end{figure}

The first step, is dedicated to the mapping between Cartesian and polar coordinates. To this end, consider $c$ a given cell of the Lagrangian grid for $(X,Y)$-coordinates, ${p} \in {\cal P}({c})$ a node of this cell. Notation used in the sequel are depicted on \figref\ref{fig:GCNS}. The mapping between a point $p\in c$ of Cartesian coordinates $\Xp = (X,Y)^t$ to $\HVX_p=(r_p,\theta_p)^t$ in polar ones is done through the following linear transformation
\begin{equation}
\HVX_p = \DT \VX_p,
\label{eq:remap}
\end{equation}
where $\DT(\VX_p) \stackrel{def}{=} I_d - \beta I_d + \beta \DR(\VX_p)$ with $\beta \in \{0,1\}$. For $\beta = 1$ then $\DT(\VX_p)~=~\DR(\VX_p)$ with the rotation matrix 
$$
\DR(\VX_p) = 
\left(
\begin{array}{cc}
\cos(\theta_p) & -\sin(\theta_p) \\
\sin(\theta_p) & \cos(\theta_p)
\end{array}
\right),
$$
using the definition $\theta_p = \arctan\left(\dfrac{Y_p}{X_p}\right)$ and $r_p = \sqrt{X_p^2+Y_p^2}$. In the case $\beta = 0$, this formula leads to Cartesian rezoning with the identical transformation $\HVX_p = \VX_p$.
When mapping $(X,Y)$ to $(r,\theta)$, the origin node has to be specifically treated. Indeed the transformation \eqref{eq:remap} is not defined for this point. 
Then as it is needed in the rezoning algorithm in the $(r,\theta)$ frame, the origin node is defined by a mapping of the first row on $r=0$ axis (see \figref \ref{fig:GCNS}). Note that these nodes are not used for the final backward mapping.

The second step is the GCNS algorithm. It is based on a minimization problem of a local functional that controls the quality of the mesh. As done in \cite{Galera1,Galera2}, one has to distinguish boundary nodes and internal node for which the smoothing procedure is different.\\
For internal nodes, let us introduce as in \cite{Knupp1} the condition number for $(r,\theta)$-coordinates that writes
\begin{equation}
 \kappa(\Jtp)= \dfrac{||\widehat{\VX}_{pp^+}||^2+||\widehat{\VX}_{pp^-}||^2}{\Atp},
\end{equation}
where $\widehat{\VX}_{pp^{\pm}}=\Xtp-\widehat{\VX}_{p^{\pm}}$, and  $\Atp=\det(\Jtp)$ is the area of the triangle delimited by $\{\tp,\tp^+,\tp^-\}$ in the rezoned grid and $\Jtp=[\widehat{\VX}_{pp^+}|-\widehat{\VX}_{pp^-}]$ the $2\times2$ Jacobian matrix associated to each corner at vertex $\tp$ of cell $c$. Thanks to this condition number we define the local function associated to the node~$\tp$ 
\begin{equation}
 \Ftp(\Xtp)=\sum_{c \in{\cal C}(\tp)} \kappa(\Jtp),
\end{equation}
Finally, the new position $\Xtp^{rez}$ is obtained by the minimization of the local function $\Ftp$ using the first step of a Newton algorithm. This leads to the formula
\begin{equation}
 \Xtp^{rez} =  \Xtp-\Htp^{-1}(\Xtp) \nabla \Ftp(\Xtp),
\end{equation}
where $\Htp^{-1}$ and $\nabla \Ftp$ are respectively the Cartesian $2\times2$ Hessian matrix and gradient  related to the local functional $\Ftp$.\\
For boundary nodes, the rezoned position $\Xtp^{rez}$ of $\tp$ is computed in consistent way with the GCNS algorithm. To this end, $\Xtp^{rez}$ is given thanks to a second-order interpolation B\'ezier curve \cite{Galera2} leading to
\begin{equation}
\Xtp^{rez}= \BXtp(s^{rez}) = (1-(s^{rez})^2)\Xtpm+2(1-s^{rez})\widehat{\bf X}_i+(s^{rez})^2\Xtpp,
\end{equation}
where $\widehat{\bf X}_i$ such that $\BXtp(1/2) = \Xtp$. Furthermore, the parameter $s^{rez}$ is computed to minimize $\Ftp(\BXtp(s))$ (for more details on this procedure see  \cite{Galera2}).\\

Finally, the third step consists in backward mapping between $\Xtp^{rez}$ and $\Xp^{rez}$ using \eqref{eq:remap}, where the inverse of the transformation matrix is taken equal to $\DT^{-1}(\VX_p) \stackrel{def}{=} I_d - \beta I_d + \beta \left[\DR(\VX_p)\right]^{-1}$ with $\beta \in \{0,1\}$ and $\left[\DR(\VX_p)\right]^{-1}$ the inverse rotation matrix . 

\subsection{Relaxation algorithm}

The relaxation algorithm consists in making a convex combination between rezoned grid obtained from GCNS step and its location after Lagrangian step. This reads for each mesh node $p$ by: 
$$
\widetilde{\VX}_p = \Xp +\omega_p(\Xp^{rez}-\Xp), \text{ with } \omega_p \in [0,1],
$$
where $\widetilde{\VX}_p $ is the new mesh node position after the complete rezoning phase. The coefficient $\omega_p$ is computed as a function of the right Cauchy-Green tensor associated to the Lagrange grid deformation over a time step (for details see \cite{Galera2, Loubere1}).

\subsection{Numerical validation}
In this section, we compare results obtained by the GCNS algorithm to those obtained for classical CNS for the rezoning of uniform polar and unstructured meshes.

\begin{figure}[h!]
\begin{tabular}{ccc}
\includegraphics[scale=0.2]{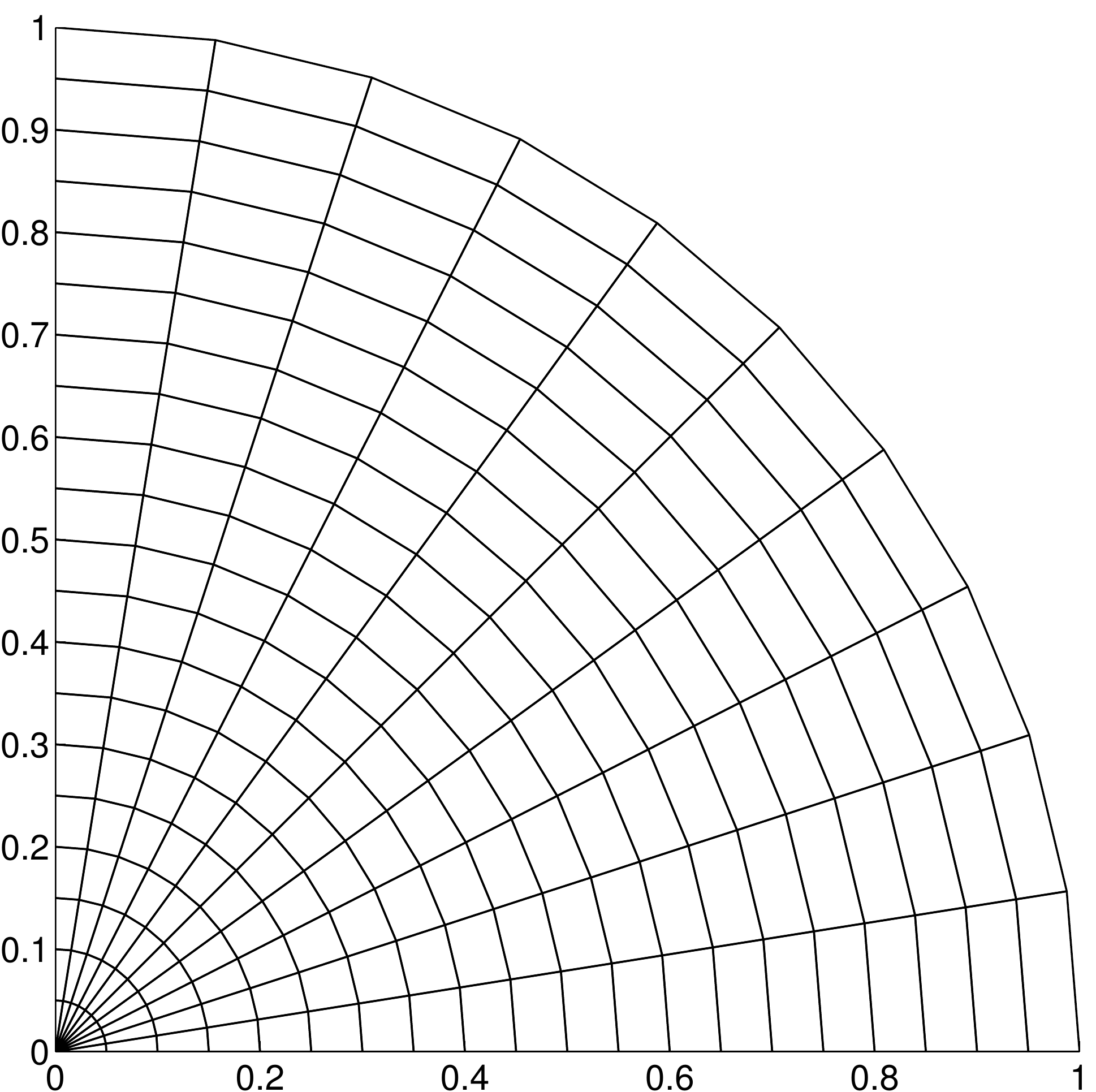} & 
\includegraphics[scale=0.2]{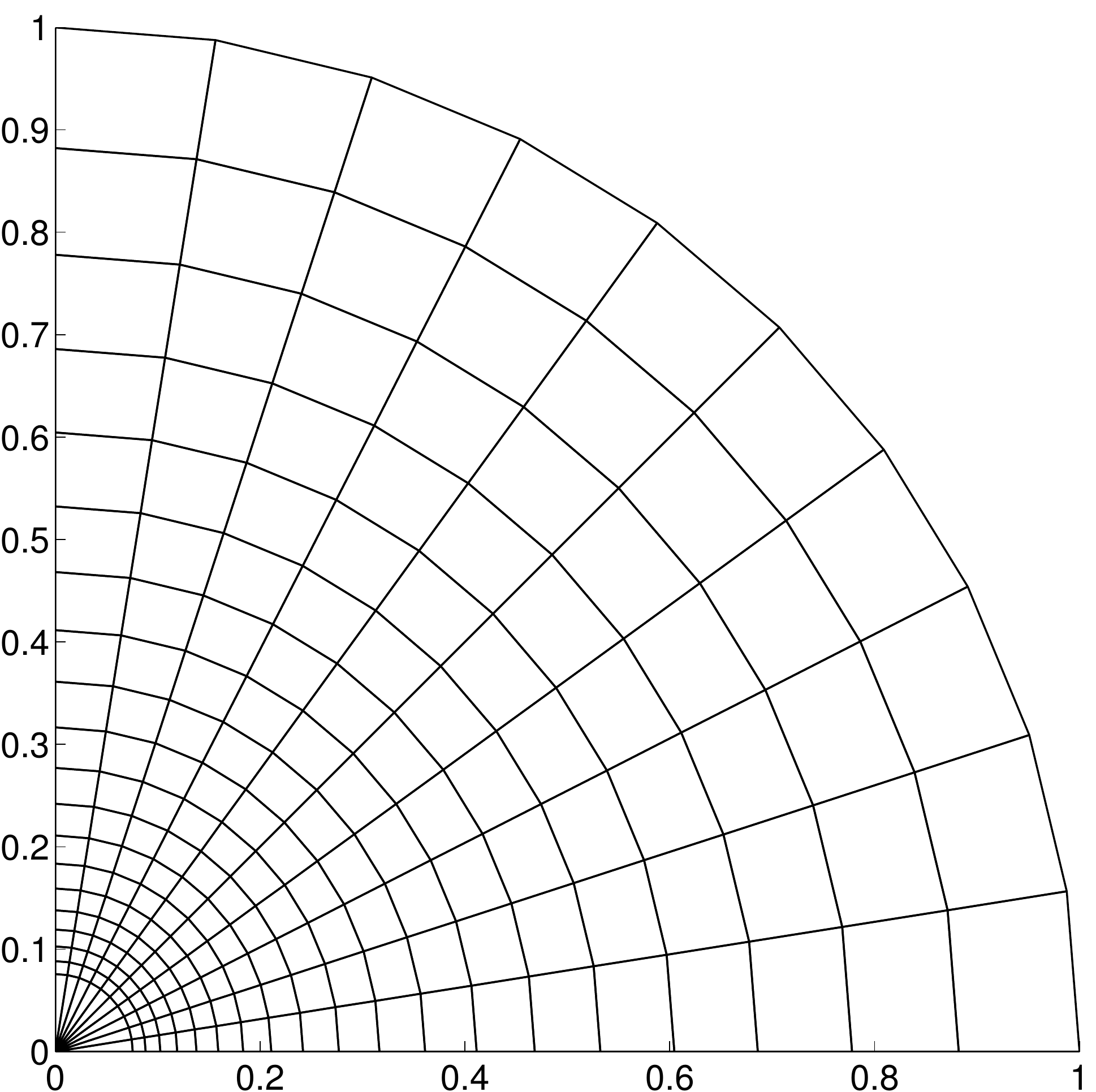} & 
\includegraphics[scale=0.2]{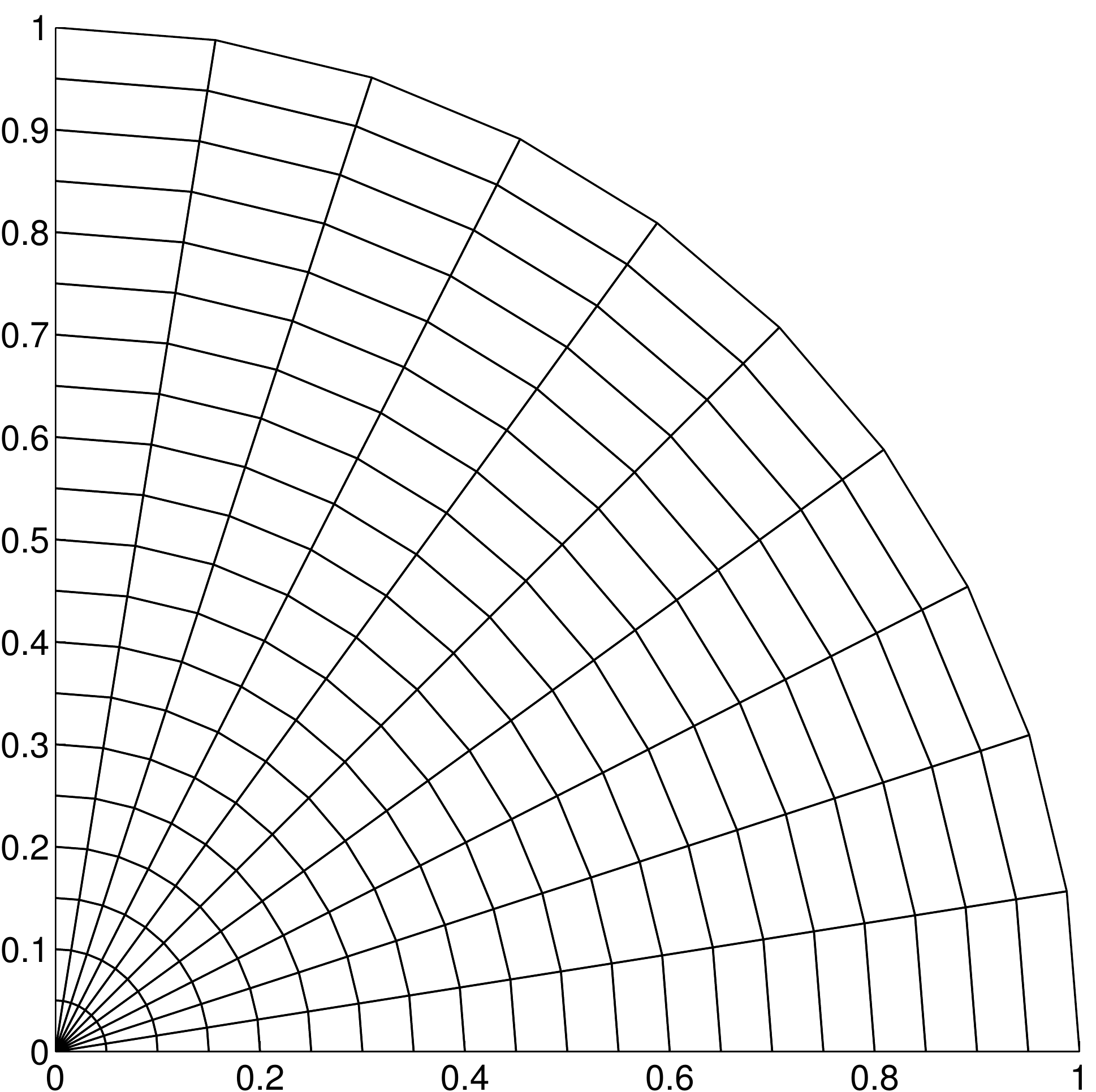} \\
(a) & (b) & (c) 
\end{tabular}
\caption{Smoothing of a static polar grid $16\times10$: (a) initial grid; Smoothed grids after 100 iterations: (b) CNS, (c) GCNS. \label{fig:2}}
\end{figure}

\paragraph{\bf Uniform mesh} First, we consider an uniform polar mesh made of $20\times10$ elements see \figref \ref{fig:2}-(a). Results obtained after 100 iterations for the classical and general smoothing are presented on \figref  \ref{fig:2}. For each method the relaxation coefficient $\omega_p$ is taken equal to 1. As already mentioned, the classical smoothing does not converge on polar mesh and implies the collapse of cell layers to the origins (see \figref \ref{fig:2}-(b)). However, for the GCNS, the result obtained (see \figref \ref{fig:2}-(c)) is converged. The mesh initially uniform, is not modified at the end of the computation. This clearly illustrates the good behavior of our smoothing algorithm.

\begin{figure}[h!]
\begin{tabular}{ccc}
\includegraphics[scale=0.2]{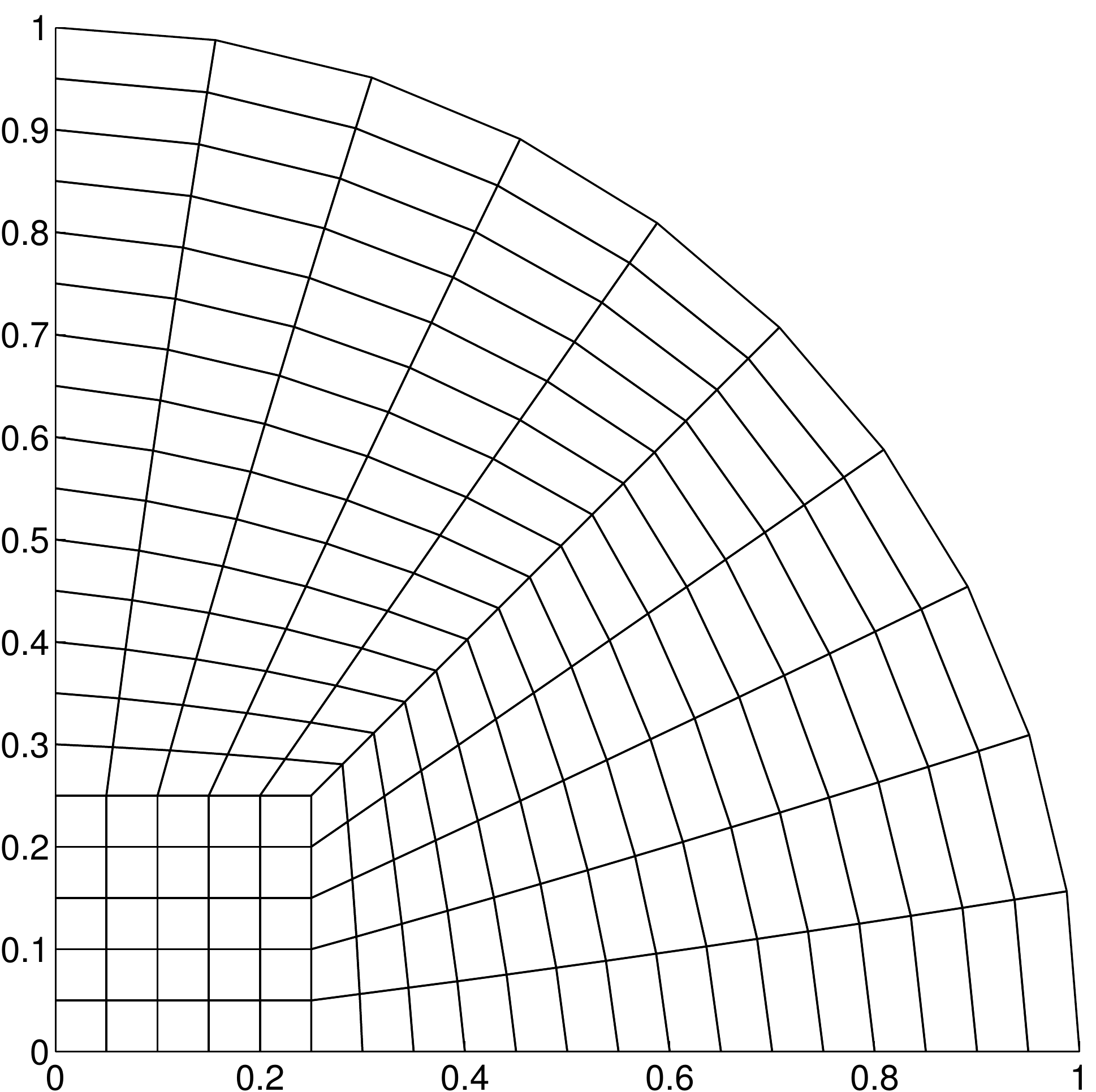} & 
\includegraphics[scale=0.2]{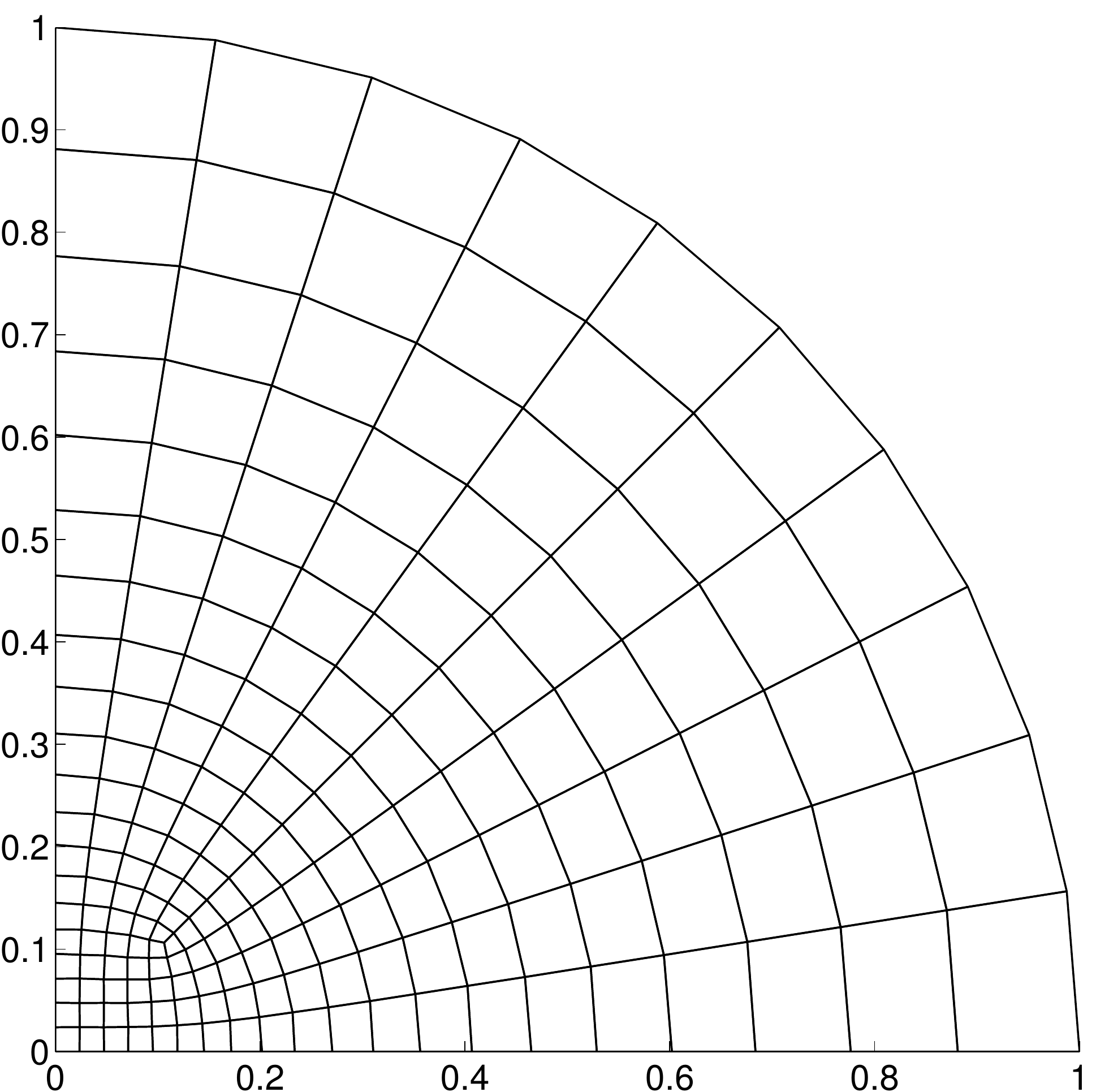} & 
\includegraphics[scale=0.2]{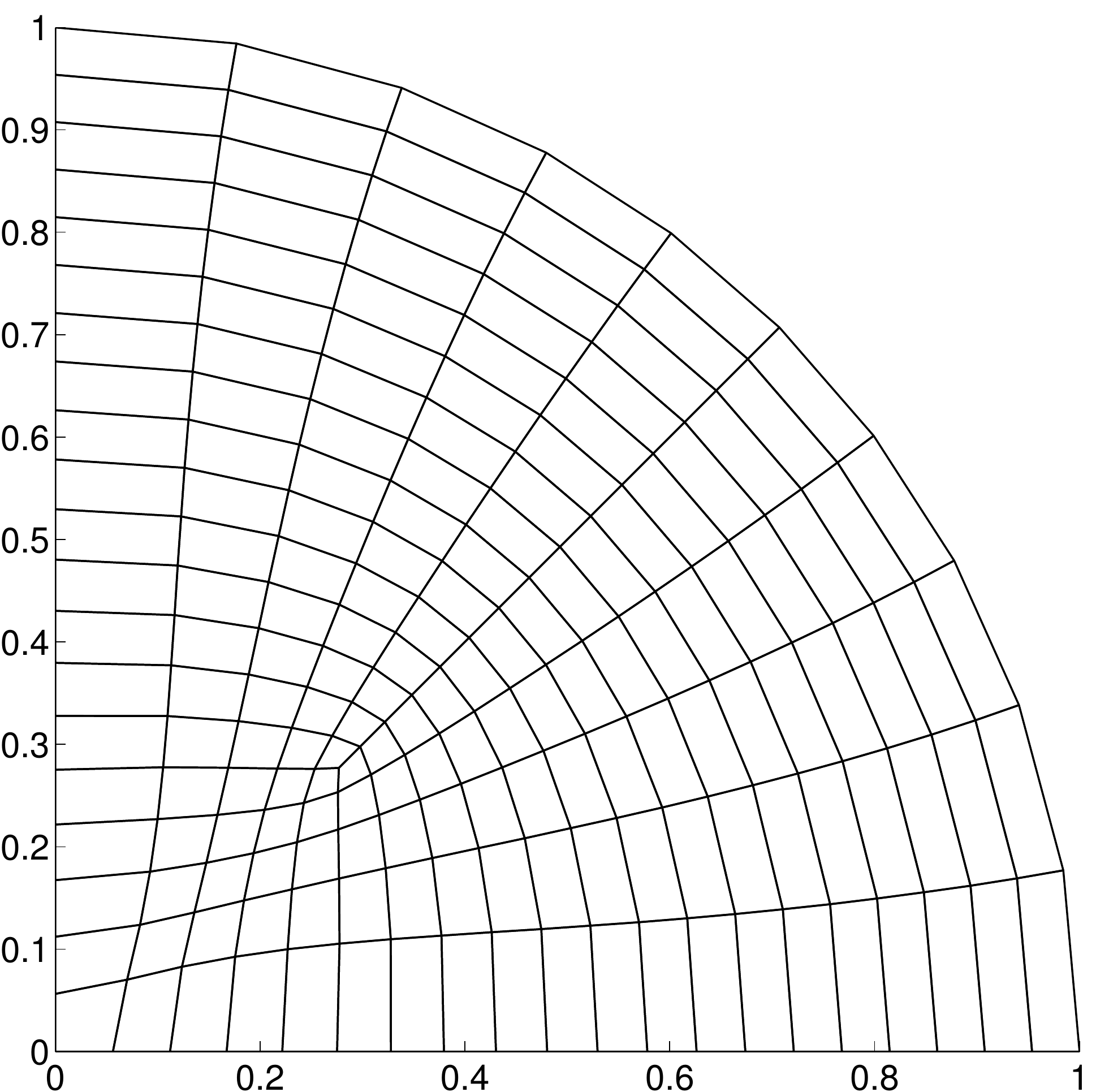} \\
(a) & (b) & (c) \\
\end{tabular}
\caption{Smoothing of a static unstructured grid: (a) initial grid; Smoothed grids after 100 iterations (b) CNS, (c) GCNS after 100 iterations. \label{fig:3}}
\end{figure}

\begin{figure}[h!]
\begin{tabular}{ccc}
\includegraphics[scale=0.2]{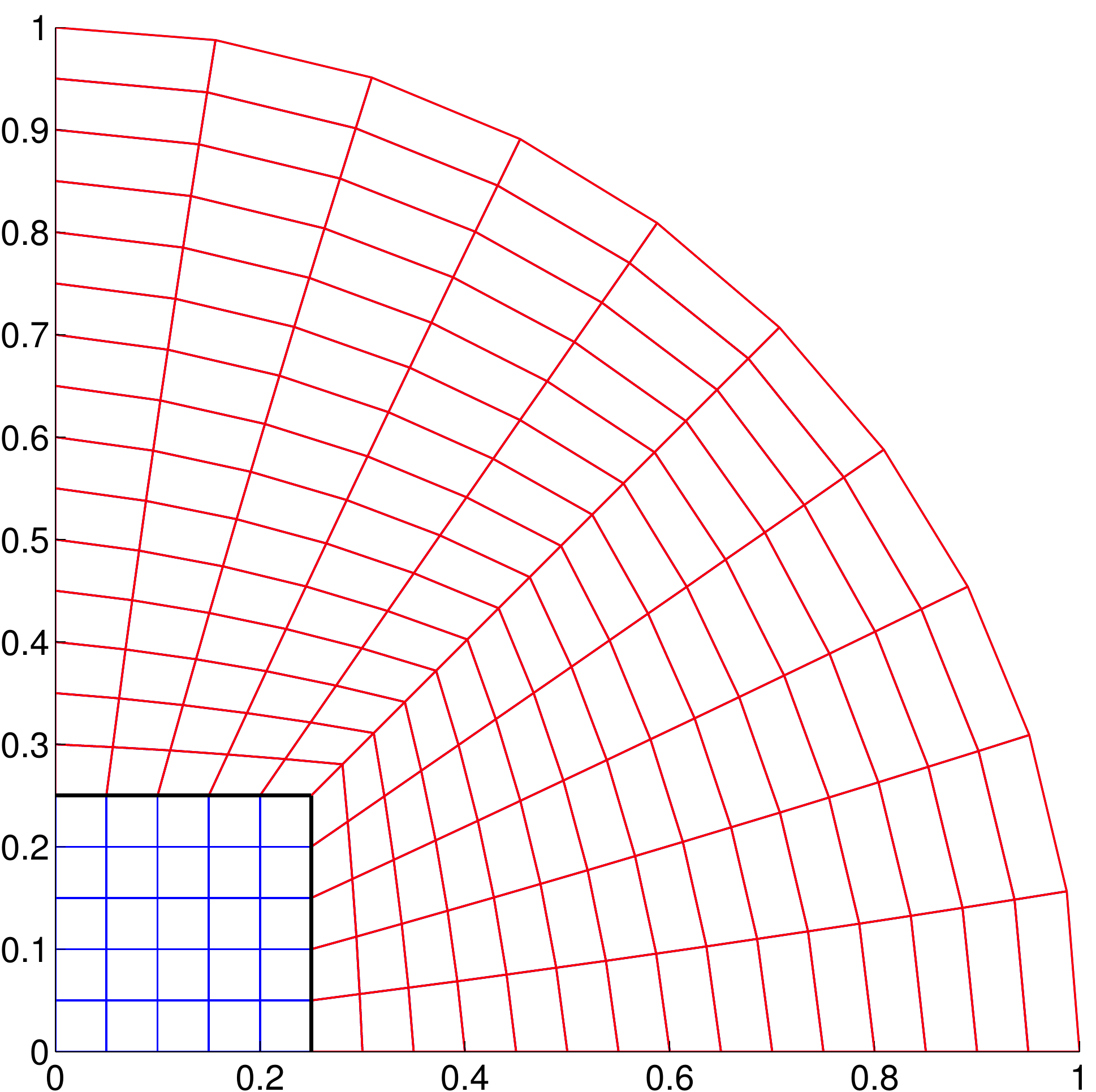} & 
\includegraphics[scale=0.2]{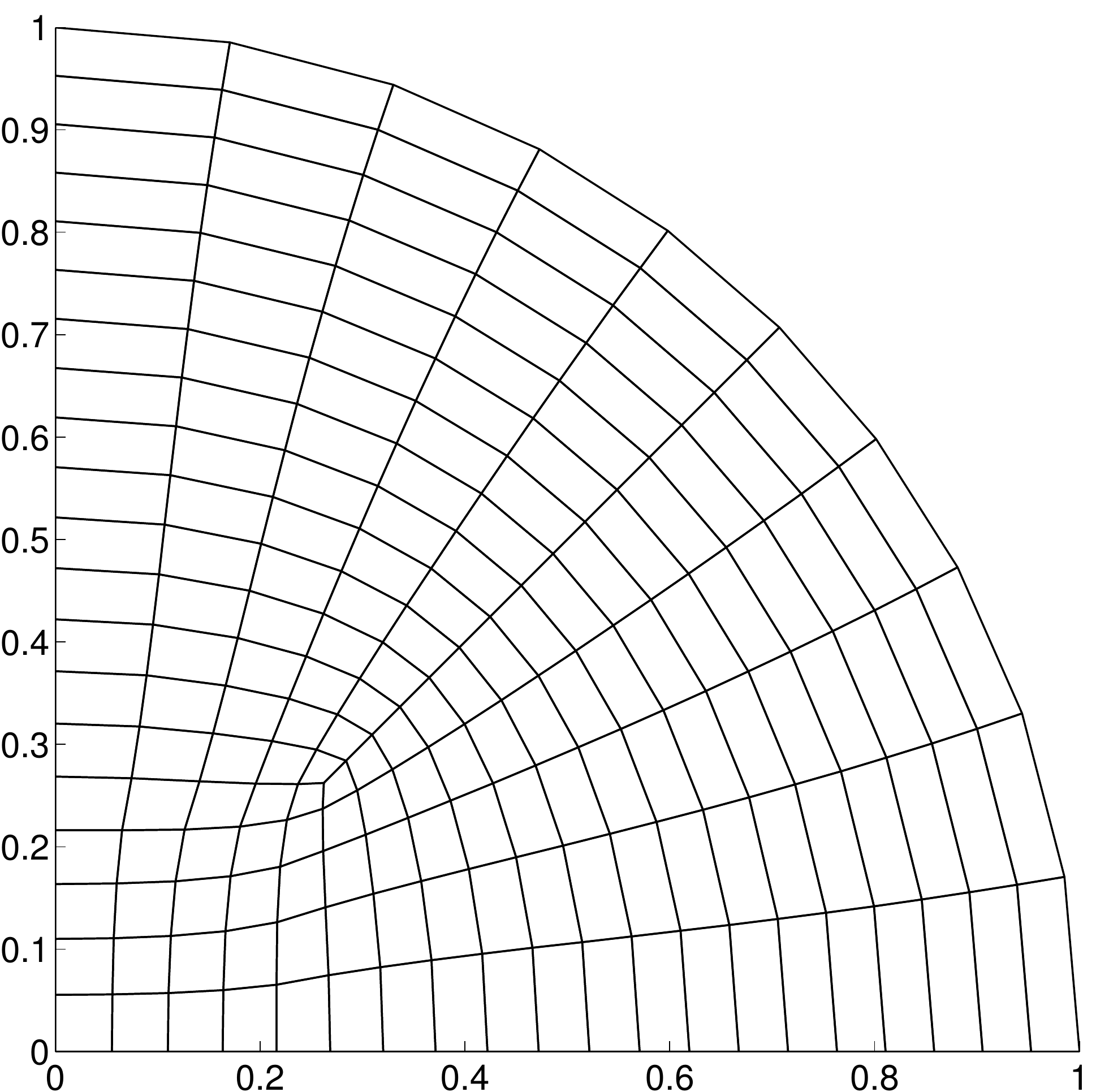} & 
\includegraphics[scale=0.2]{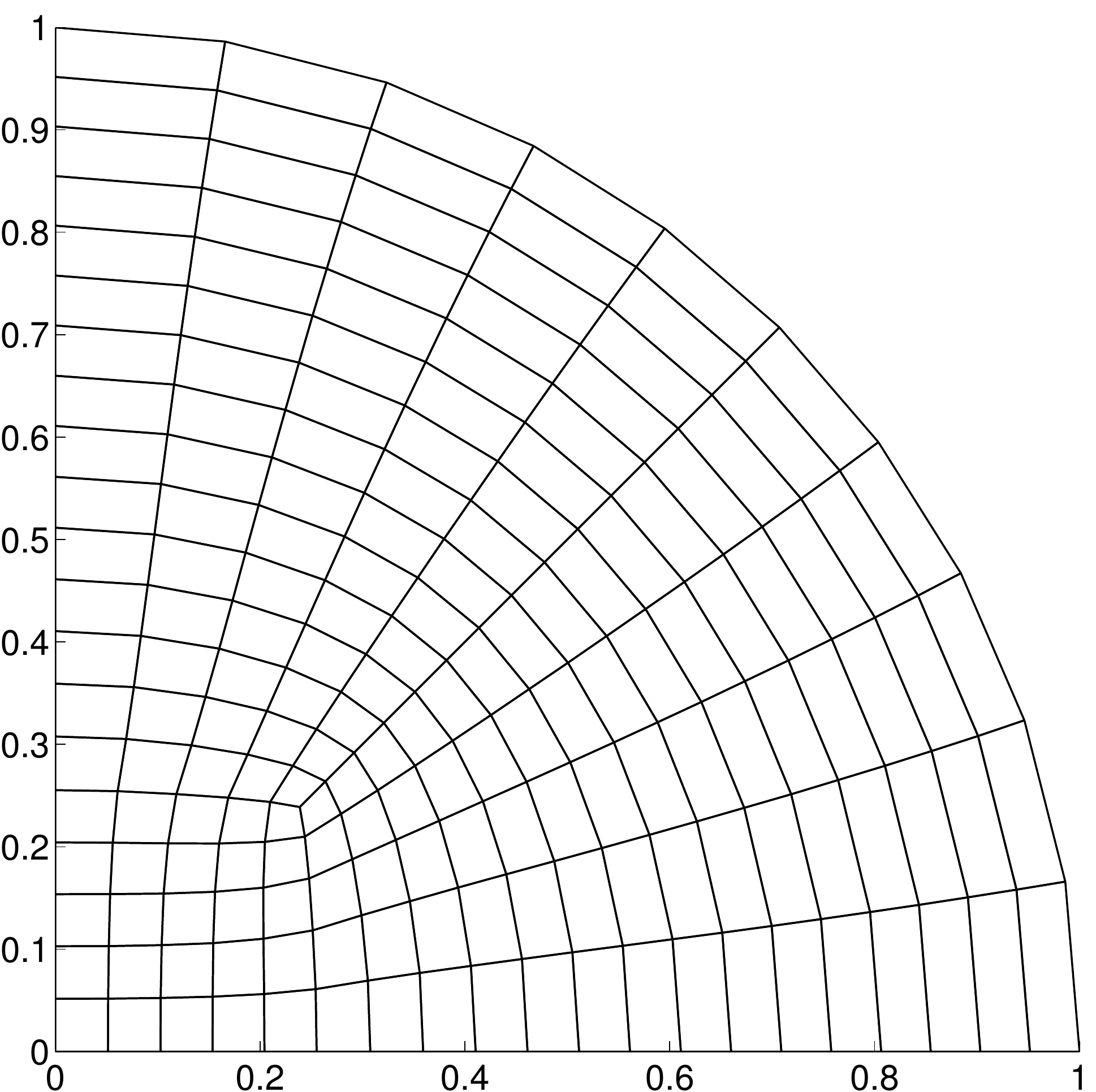} \\
(a) & (b) & (c)\\
\end{tabular}
\caption{Smoothing of a static unstructured grid: (a) initial grid with Cartesian (blue) and polar (red) rezoning regions; Smoothed grids after 100 iterations (b) GCNS with interfacial polar rezoning, (c) GCNS with interfacial Cartesian rezoning. \label{fig:smooth}}
\end{figure}

\paragraph{\bf Unstructured mesh} Now, rezoning for an unstructured mesh is studied. Let us consider a mesh made of $175$ quadrangular cells as depicted on \figref \ref{fig:3}-(a).
When applying the full Cartesian rezoning to the mesh, similar observations as previously can be made. It suffers from an implosion of central cells to the origin and does not converge (see  \figref\ref{fig:3}-(b)).  For the full GCNS algorithm, one can see  after convergence, the formation of mesh distortion on the square region and a polar mesh far from the center (see \figref\ref{fig:3}-(c)). 
Nevertheless, it is possible to improve this rezoning. Thus, the main idea developed in the sequel is to apply the GCNS rezoning algorithm differently for a node belonging initially to a Cartesian or polar region of the mesh. To this end, the transformation $T$ between $(X,Y)$ and $(r,\theta)$ coordinates is modified in the following way
\begin{equation}
 \DT(\VX_p) \stackrel{def}{=} I_d - \beta(\VX_p) I_d + \beta(\VX_p) \DR(\VX_p), 
\end{equation}
with 
$$
\beta(\VX_p) = 
\left\{
\begin{array}{rl}
 1 & \text{ if } \VX^0_p \in {\cal P}^{pol}, \\
 0 & \text{ if } \VX^0_p \in {\cal P}^{car},
\end{array}
\right.
$$
where ${\cal P}^{car}$ and ${\cal P}^{pol}$ are the sets of nodes that belong to the Cartesian and respectively polar region of the mesh at the initial time. These regions are represented thanks to red and blue color (see \figref \ref{fig:smooth}-(a)) for the considered mesh. Nodes localized at the frontier between the polar and Cartesian meshes (black nodes on \figref \ref{fig:smooth}-(a)) can be considered either polar, or Cartesian. As represented on \figref \ref{fig:smooth}-(b,c), both possibilities are tested. The obtained results illustrate that the Cartesian choice remains better contrary to the polar one that introduce mesh distortion.

\section{Hybrid remapping in axisymmetric geometry}

During the remapping phase, the physical unknowns (density, velocity, total energy) computed thanks to the Lagrangian step are conservatively remapped from the Lagrangian mesh to the rezoned one. To this end, an extension of the {\it Hybrid Remapping Algorithm} for multi-material flows \cite{Galera1,Kucharik1,Berndt1} to cylindrical geometry is proposed here. This strategy consists in the following two steps. First a {\it swept-faced remapping} is used to  treat cells and nodes localized far from the interface. Then, a {\it cell-intersection-based} method \cite{Galera2} is applied to the cells and nodes in the neighborhood of the interface. In this way, this approach combines the ability of the cell-intersection method  to remap the interface and the efficiency of the swept flux approach for the other cells that significantly reduce the global computing cost of the method. As done previously, in the perspective of general use of the method, a global formulation including both Cartesian and axisymmetric framework is presented.\\

We assume in the sequel, that there is no topology change of the mesh, the cells of the Lagrangian and rezoned grids are respectively designed by $\Omega_c$ and $\TOm_c$. 

\subsection{Multi-material cell-intersection-based (MCIB) remapping}

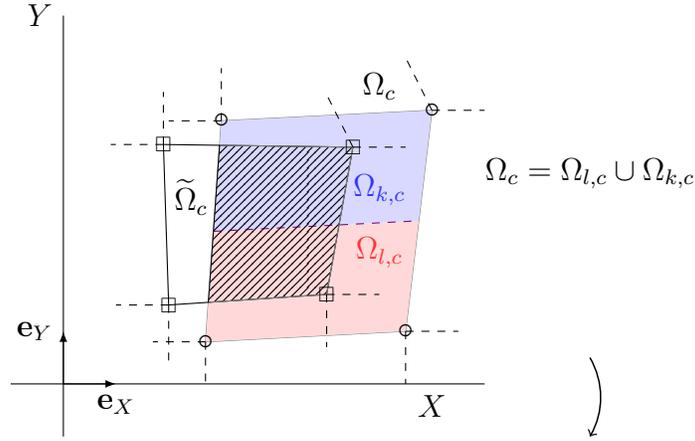
\begin{figure}[h!]
\centering
\begin{tikzpicture}[xscale=0.7,yscale=0.7]
\tikzstyle{fleche1}=[-,>=latex]
\tikzstyle{fleche2}=[<->,>=latex]
\tikzstyle{fleche3}=[->,>=latex]
\newcommand{\quadA}{(6,5)--(5.7,0.8)--(9.5,1)--(10,5.2)--cycle}
\newcommand{\quadC}{(5.85,2.9)--(5.7,0.8)--(9.5,1)--(9.75,3.1)--cycle}
\newcommand{\quadCC}{(6,5)--(5.85,2.9)--(9.75,3.1) --(10,5.2)--cycle}
\newcommand{\quadB}{(4.9,4.55)--(5,1.5)--(8,1.7)--(8.5,4.5)--cycle}
\draw[color=gray] \quadA;
\fill[fill=red!15] \quadC;
\fill[fill=blue!15] \quadCC;
\draw[color=black] \quadB;
\draw[violet,dashed] (5.85,2.9)--(9.75,3.1);
\begin{scope}
\clip \quadA;
\draw[color=gray,fill=gray!20,pattern=north east lines] \quadB;
\end{scope}
\begin{scope}
\clip \quadB;
\draw[color=black] \quadA;
\end{scope}
\draw[fleche1] (2,0)--(11,0); 
\draw[fleche1] (3,-1)--(3,7);
\node[left] at (3,1) {${\bf e}_Y$};
\node[below] at (4,0) {${\bf e}_X$};
\node[left] at (3,7) {$Y$};
\node[below] at (10,0) {$X$};
\draw[fleche3] (3,0)--(3,1); 
\draw[fleche3] (3,0)--(4,0);
\node at (6,5) {$\circ$};
\node at (5.7,0.8) {$\circ$};
\node at (9.5,1) {$\circ$};
\node at (10,5.2) {$\circ$};

\node[black] at (4.9,4.55) {\scriptsize $\square$};
\node[black] at (5,1.5) {\scriptsize $\square$};
\node[black] at (8,1.7) {\scriptsize $\square$};
\node[black] at (8.5,4.5) {\scriptsize $\square$};
\draw[dashed,black](4.9,4.55)--(4.9,5.55);
\draw[dashed,black](4.9,4.55)--(3.9,4.55);
\draw[dashed,black](5,1.5)--(5,0.4);
\draw[dashed,black](5,1.5)--(4,1.5);
\draw[dashed,black](8,1.7)--(8.,0.7);
\draw[dashed,black](8,1.7)--(9,1.7);
\draw[dashed,black](8.5,4.5)--(9.5,4.5);
\draw[dashed,black](8.5,4.5)--(8,5.5);

\draw[dashed](6,5)--(6,6);
\draw[dashed](6,5)--(5,5);
\draw[dashed](5.7,0.8)--(5.7,0);
\draw[dashed](5.7,0.8)--(4.7,0.8);
\draw[dashed](9.5,1)--(9.5,0);
\draw[dashed](9.5,1)--(10.5,1);
\draw[dashed](10,5.2)--(9.5,6.2);
\draw[dashed](10,5.2)--(11,5.2);
\node[above] at (9,5.2) {$\Omega_{c}$};
\node[blue!80,above] at (9,3.2) {$\Omega_{k,c}$};
\node[red!80,above] at (9,2) {$\Omega_{l,c}$};
\node[below right] at (4.9,4.1) {$\TOm_c$};
\draw[<-,line width=0.2mm] (13,-1) arc (-30:30:1.5);
\node  at (13,4) {$\Omega_{c}=\Omega_{l,c}\cup\Omega_{k,c} $};
\end{tikzpicture} 
\caption{Notations for MCIB method. \label{fig:7a}}
\end{figure}

The main goal of remapping is as follows. Given the piecewise constant representation of the physical variables per unit of volume ($\rho,\rho \VU, \rho E$) noted $\psi_c=\rho_c \phi_c$ in each cell of the Lagrangian grid, we want to compute its equivalent $\widetilde{\psi}_c$ in each cell of the rezoned grid given as 
\begin{equation}
\widetilde{\psi}_c = \dfrac{1}{\widetilde{V}_c}\int_{\TOm_c} \widetilde{\rho \phi} \CR dA,
\end{equation}
with $\widetilde{V}_c$ the volume of the cell $\TOm_c$. 
Contrary to single fluid approach, here the rezoned values $\widetilde{\psi}_c$ can not be computed directly in each cell $c$. In fact, one has to take into account multi-material aspects.\\ 

First of all, let us introduce some notations. Each material of the flow noted $k$ occupies the polygon $\Omega_{k,c} \subset \Omega_c$, within the MOF framework, such that $\displaystyle\Omega_c= \bigcup_{k}\Omega_{k,c}$ and is characterized by its partial mass, density, pressure, internal energy and variables per unit of volume (total energy, momentum) whose averaged values in each sub-cell are respectively $m_{k,c}, \rho_{k,c}, P_{k,c}, \varepsilon_{k,c}$ and $\psi_{k,c} = \rho_{k,c} \phi_{k,c}$ with $\phi_{k,c}$ the partial velocity or energy per unit of mass.\\

Thus, for multi-material flow, the main idea of remapping is not to directly compute the global rezoned quantities $\widetilde{\psi}_c$ but the partial rezoned ones noted $\widetilde{\psi}_{k,c}$. This is particularly true for the MCIB method that is dedicated to treat cell in the interface neighborhood. To this end, we first propose a second order reconstruction $\Psi_{k,c}(\VX)$ of $\psi_{k,c}$ over each Lagrangian cell $c$ through  the piecewise linear function
\begin{equation}
\Psi_{k,c}(\VX) = \psi_{k,c} + (\nabla \Psi_k)_c(\VX-\VX_{k,c}),
\end{equation}
where  $ (\nabla \Psi_k)_c$ denotes the constant gradient of $\Psi_{k,c}$ within cell $c$ computed thanks to a least-squares approach. Finally $\VX_{k,c}$ is the centroid related to the $k$-th fluid in the cell $c$ given by
\begin{equation}
 \VX_{k,c} = \dfrac{1}{ V_{k,c}} \int_{\Omega_{k,c}} \CR \VX dA.
\end{equation}
Thanks to these notations, the remapped value for MCIB is given by
\begin{equation}
 \widetilde{\psi}_{k,c} = \dfrac{1}{\widetilde{ V}_{k,c}} \sum_{d\in \CC(c)} \int_{\Omega_{k,d}\cap\widetilde{\Omega}_c} \CR \Psi_{k,c} dA,
 \label{eq:rezpsik}
\end{equation}
where the intersection polygons $\Omega_{k,d}\cap\widetilde{\Omega}_c$ are computed thanks to a specific triangulation of the mesh.  The procedure is detailed in \cite{Galera2}. The set $ \CC(c)$ contains the cells  including $c$ that share at least one node with the cell $c$. At last, the partial volume defined on the rezoned cell is $\displaystyle \widetilde{V}_{k,c} = \sum_{d\in  \CC(c)} \int_{\Omega_{k,d}\cap\widetilde{\Omega}_c} \CR dA $.\\ 

In the context of MOF reconstruction, one has to define additional quantities as the partial remapped mass corresponding to material $k$. It is computed as $\widetilde{m}_{k,c}= \widetilde{\rho}_{k,c} \widetilde{\alpha}_{k,c} \widetilde{V}_{k,c}$ with the volume fraction
\begin{equation}
 \widetilde{\alpha}_{k,c} = \dfrac{1}{ \widetilde{V}_{c}} \sum_{d\in \CC(c)} \int_{\Omega_{k,d}\cap\widetilde{\Omega}_c} \CR dA,
\end{equation} 
thus the partial volume can be also expressed as $ \widetilde{V}_{k,c} = \widetilde{V}_{c} \widetilde{\alpha}_{k,c}$.
In addition, each material centroid position is defined thanks to
\begin{equation}
 \widetilde{\VX}_{k,c} = \dfrac{1}{ \widetilde{V}_{k,c}} \sum_{d\in  \CC(c)} \int_{\Omega_{k,d}\cap\widetilde{\Omega}_c} \CR \VX dA.
\end{equation}

\subsection{Pure cell swept-face (PCSF)  remapping }

\begin{figure}[h!]
\centering
\begin{tikzpicture}[xscale=0.7,yscale=0.7]
\tikzstyle{fleche1}=[-,>=latex]
\tikzstyle{fleche2}=[<->,>=latex]
\tikzstyle{fleche3}=[->,>=latex]
\newcommand{\quadA}{(6,5)--(5.7,0.8)--(9.5,1)--(10,5.2)--cycle}
\newcommand{\quadB}{(4.9,4.55)--(5,1.5)--(8,1.7)--(8.5,4.5)--cycle}
\newcommand{\quadC}{(8,1.7)--(8.5,4.5)--(10,5.2)--(9.5,1)--cycle}
\draw[color=black,fill=blue!10] \quadA;
\draw[color=black] \quadB;
\filldraw[color=gray,fill=gray!20,pattern= horizontal lines] \quadC;
\draw[fleche1] (3,-1)--(3,7);
\draw[fleche1] (2,0)--(11,0); 
\node[below] at (11,0) {$X$};
\node[left] at (3,1) {${\bf e}_Y$};
\node[below] at (4,0) {${\bf e}_X$};
\node[left] at (3,7) {$Y$};
\draw[fleche3] (3,0)--(3,1); 
\draw[fleche3] (3,0)--(4,0);
\node at (6,5) {$\bullet$};
\node at (5.7,0.8) {$\bullet$};
\node at (9.5,1) {$\bullet$};
\node at (10,5.2) {$\bullet$};
\node[above right] at (10.1,5.2) {$p^+$};
\node[right] at (9.6,0.5) {$p$};

\node at (4.9,4.55) {\scriptsize $\blacksquare$};
\node at (5,1.5) {\scriptsize $\blacksquare$};
\node at (8,1.7) {\scriptsize $\blacksquare$};
\node at (8.5,4.5) {\scriptsize $\blacksquare$};
\node[left] at (7.9,2.) {$p$};
\node[left] at  (8.3,4.2) {$p^+$};

\draw[dashed](4.9,4.55)--(4.9,5.55);
\draw[dashed](4.9,4.55)--(3.9,4.55);
\draw[dashed](5,1.5)--(5,0.4);
\draw[dashed](5,1.5)--(4,1.5);
\draw[dashed](8,1.7)--(8.,0.7);
\draw[dashed](8,1.7)--(9,1.7);
\draw[dashed](8.5,4.5)--(9.5,4.5);
\draw[dashed](8.5,4.5)--(8,5.5);

\draw[dashed](6,5)--(6,6);
\draw[dashed](6,5)--(5,5);
\draw[dashed](5.7,0.8)--(5.7,0);
\draw[dashed](5.7,0.8)--(4.7,0.8);
\draw[dashed](9.5,1)--(9.5,0);
\draw[dashed](9.5,1)--(10.5,1);
\draw[dashed](10,5.2)--(9.5,6.2);
\draw[dashed](10,5.2)--(11,5.2);
\node[above] at (9,5.2) {$\Omega_c$};
\node at (5.5,3.8) {$\TOm_c$};
\node[fill=white] at (9,3.2) {$A_{f}$};
\node at (11,3.2) {$\TOm_{c^+}$};
\draw[<-,line width=0.2mm] (13,-1) arc (-30:30:1.5);
\node  at (13.5,4) {$\Omega_{c}=\Omega_{k,c}$};
\node  at (13.5,3) {$\Omega_{c^+}=\Omega_{k,c^+}$};
\end{tikzpicture} 
\caption{Notations for swept face-based method. \label{fig:7b}}
\end{figure}
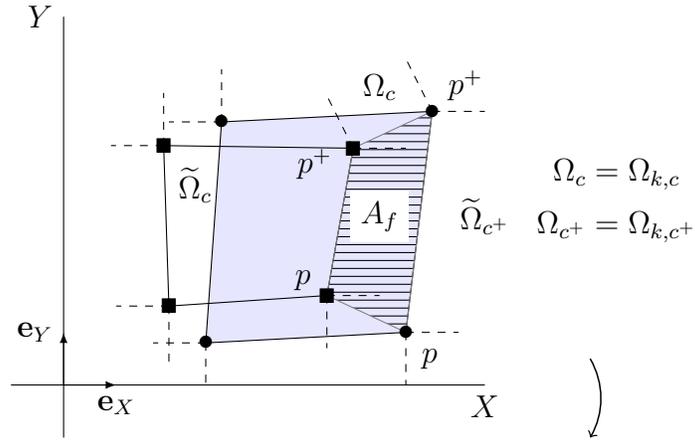

As explained before, the PCSF remapping is used only to treat single fluid cells.  In this context, one should remark that $\Omega_c = \Omega_{k,c}$, thus the mean value $\widetilde{\psi}_{c,k}$ is given through
\begin{equation}
 \widetilde{\psi}_{k,c} =  \psi_{k,c} + \sum_{f\in \CF(c)} \int_{A_f} \CR \Psi_{k,f} dA, 
 \label{eq:rezpsi}
\end{equation}
with $A_f$ the quadrangular signed area swept by the face $f$ of a cell $c$ between the Lagrangian grid and the rezoned grid delimited by the ordered nodes of coordinates $\{\VX_p,\TX_p,\TX_{p^+},\VX_{p^+}\}$ (refer to fig \figref \ref{fig:7b}). We note $\CF(c)$ the set of the faces $f$ of a cell $c$. In addition, $\Psi_{k,f}$ is the upwind value given by
\begin{equation}
 \Psi_{k,f} = 
 \left\{
 \begin{array}{ll}
  \Psi_{k,c^+} & \text{ if } A_f >0 \\ 
  \Psi_{k,c}   & \text{ otherwise.} \\ 
 \end{array}
 \right.
\end{equation}
with $c^+$ the neighbor cell of $c$ through the face $f$. During this step the volume fractions $\widetilde{\alpha}_{k,c}= \alpha_{k,c}$ do not change as we consider single fluid cells and the material centroid can be updated directly from the geometry $\widetilde{\VX}_{k,c}= \widetilde{\VX}_{c}$ where $\widetilde{\VX}_{c}$ is the centroid of the cell $\widetilde{\Omega}_c$.

\subsection{Integration strategy}
For both PCSF and MCIB remapping, one has to compute several surface integrals, on polygons where the integrand is a polynomial function of $(X,Y)$. This can be done using a triangulation of these areas. Nevertheless, this is expensive. Here, we rather adopt a more efficient method as in \cite{Margolin1}. In this context, integrals are simplified using Taylor decomposition of the polynomial integrand and Green's formula leading to compute circular integrals over the edges of the polygons defining the integration areas.  For further details on integral computations see \cite{Margolin1}.

\subsection{Hybrid remapping algorithm}

\begin{figure}[h!]
\centering
\begin{tikzpicture}[xscale=0.7,yscale=1]
\tikzstyle{every node}=[font=\scriptsize]
 \fill[blue!10] (-1,3) rectangle (6.25,4);
 \fill[red!10] (6.25,3) rectangle (13.5,4);
 \draw[line width=1pt] (-1,4)   to[-] (13.5,4);  
 \draw[line width=1pt] (-1,3)   to[-] (13.5,3);
 \node at (0  ,4)  {\Large $\bullet$};
 \node at (2.5,4)  {\Large $\bullet$};
 \node at (5  ,4)  {\Large $\circ$};
 \node at (7.5,4)  {\Large $\circ$};
 \node at (10 ,4)  {\Large $\bullet$};
 \node at (12.5,4) {\Large $\bullet$};
  \node at (0  ,3)  {\Large $\bullet$};
 \node at (2.5,3)  {\Large $\bullet$};
 \node at (5  ,3)  {\Large $\circ$};
 \node at (7.5,3)  {\Large $\circ$};
 \node at (10 ,3)  {\Large $\bullet$};
 \node at (12.5,3) {\Large $\bullet$};
\draw[line width=1pt] (0   ,4) to  (0   ,3);
 \draw[line width=1pt] (2.5 ,4) to  (2.5 ,3);
 \draw[line width=1pt] (5   ,4) to  (5   ,3);
 \draw[line width=1pt] (7.5 ,4) to  (7.5 ,3);
 \draw[line width=1pt] (10  ,4) to  (10  ,3);
 \draw[line width=1pt] (12.5,4) to  (12.5,3);
 
 \fill[blue!10] (-1,1) rectangle (6.25,2);
 \fill[red!10] (6.25,1) rectangle (13.5,2);
 \draw[line width=1pt] (-1,2)   to[-] (13.5,2);
 \draw[line width=1pt] (-1,1)   to[-] (13.5,1); 
 \node at (0.5  ,2)  {$\blacksquare$};
 \node at (3.,2)  { $\blacksquare$};
 \node at (5  ,2)  {\Large $\circ$};
 \node at (7.5,2)  {\Large $\circ$};
 \node at (10.5 ,2)  { $\blacksquare$};
 \node at (13 ,2)  { $\blacksquare$};
 \node at (0.5 ,1)  { $\blacksquare$};
 \node at (3.  ,1)  { $\blacksquare$};
 \node at (5   ,1)  {\Large $\circ$};
 \node at (7.5 ,1)  {\Large $\circ$};
 \node at (10.5,1)  { $\blacksquare$};
 \node at (13  ,1)  { $\blacksquare$};
 \draw[line width=1pt] (0.5   ,2) to  (0.5   ,1);
 \draw[line width=1pt] (3 ,2)     to  (3 ,1);
 \draw[line width=1pt] (5   ,2) to  (5   ,1);
 \draw[line width=1pt] (7.5 ,2) to  (7.5 ,1);
 \draw[line width=1pt] (10.5  ,2) to  (10.5  ,1);
 \draw[line width=1pt] (13,2)     to  (13,1);
 
 \fill[blue!10] (-1,-1) rectangle (6.25,0);
 \fill[red!10] (6.25,-1) rectangle (13.5,0);
 \draw[line width=1pt] (-1,0)   to[-] (13.5,0);
 \draw[line width=1pt] (-1,-1)   to[-] (13.5,-1); 
 \node at (0.5  ,0)  {$\blacksquare$};
 \node at (3,0)      {$\blacksquare$};
 \node at (5.5  ,0)  {$\square$};
 \node at (8,0)      {$\square$};
 \node at (10.5 ,0)  {$\blacksquare$};
 \node at (13 ,0)    {$\blacksquare$};
 \node at (0.5  ,-1)  {$\blacksquare$};
 \node at (3,-1)      {$\blacksquare$};
 \node at (5.5  ,-1)  {$\square$};
 \node at (8,-1)      {$\square$};
 \node at (10.5 ,-1)  {$\blacksquare$};
\node at (13 ,-1)    {$\blacksquare$};
 \draw[line width=1pt] (0.5   ,0) to  (0.5   ,-1);
 \draw[line width=1pt] (3 ,0)     to  (3 ,-1);
 \draw[line width=1pt] (5.5   ,0) to  (5.5   ,-1);
 \draw[line width=1pt] (8 ,0)     to  (8 ,-1);
 \draw[line width=1pt] (10.5  ,0) to  (10.5  ,-1);
 \draw[line width=1pt] (13,0)     to  (13,-1);

\draw[color=gray,fill=gray!20,pattern= horizontal lines] (0.  ,3)--(0.  ,4)--(0.5,2)--(0.5,1)--cycle;
\draw[color=gray,fill=gray!20,pattern= horizontal lines] (2.5  ,3)--(2.5  ,4)--(3,2)--(3,1)--cycle;
\draw[color=gray,fill=gray!20,pattern= horizontal lines] (10.  ,3)--(10.  ,4)--(10.5,2)--(10.5,1)--cycle;
\draw[color=gray,fill=gray!20,pattern= horizontal lines] (12.5  ,3)--(12.5  ,4)--(13,2)--(13,1)--cycle;
\draw[color=gray,fill=gray!20,pattern=north east lines]  (3  ,0)--(3,1)--(10.5,1)--(10.5,0)--cycle;
\node at (14.6,3.5)  {\begin{tabular}{c} Lagrangian \\ grid \\ $\VX_p, \psi_{k,c}$ \end{tabular}};
\node at (14.6,-0.5)  {\begin{tabular}{c} Rezoned    \\ grid \\ $\widetilde{\VX}_p, \widetilde{\psi}_{k,c}$\end{tabular}};
\draw[->,line width=0.2mm]  (-1.5,3.5) arc (140:220:1.5) node[fill=white] at (-1.75,2.5) { Swept-face step};   \draw[->,line width=0.2mm]  (-1.5,1.4) arc (140:220:1.5) node[fill=white] at (-1.75,0.5) { Exact-intersection step};
\draw[dashed,line width=1pt,violet] (6.25,2) to  (6.25,1);
\draw[dashed,line width=1pt,violet] (6.25,0) to  (6.25,-1);
\draw[dashed,line width=1pt,violet] (6.25,3) to  (6.25,4);
\node at (2   ,-2) {{\large $\bullet$}: lagrangian pure nodes };
\node at (2.2 ,-2.5) {{\large $\circ$}: lagrangian mixed nodes};
\node at (8   ,-2) {{ $\blacksquare$}: rezoned pure nodes};
\node at (8.2 ,-2.5) {{$\square$}: rezoned mixed nodes};
\draw[color=gray,fill=gray!20,pattern= horizontal lines] (-0.5,-3.1) rectangle (0,-2.9); 
\node at (1.5,-3) {: swept region }; 
\draw[color=gray,fill=gray!20,pattern=north east lines] (5.5,-3.1) rectangle (6,-2.9); 
\node at (8.5,-3.) {: exact intersection region };
\draw[dashed,line width=1pt,violet] (12,-3.) -- (12.5,-3.);
\node at (13.5,-3.) {: interface }; 
\draw (-1,4.5)   to[->] (2.5,4.5);
\node[fill=white] at (0.75,4.5) {pure cells}; 
\draw (2.5,4.5)   to[<->] (10,4.5);
\node[fill=white] at (6.25,4.5) {mixed cells}; 
\draw (10,4.5)   to[<-] (13.5,4.5);
\node[fill=white] at (11.75,4.5) {pure cells}; 
\end{tikzpicture} 
\caption{Hybrid remapping principle in one-dimension case. \label{fig:8}}
\end{figure}
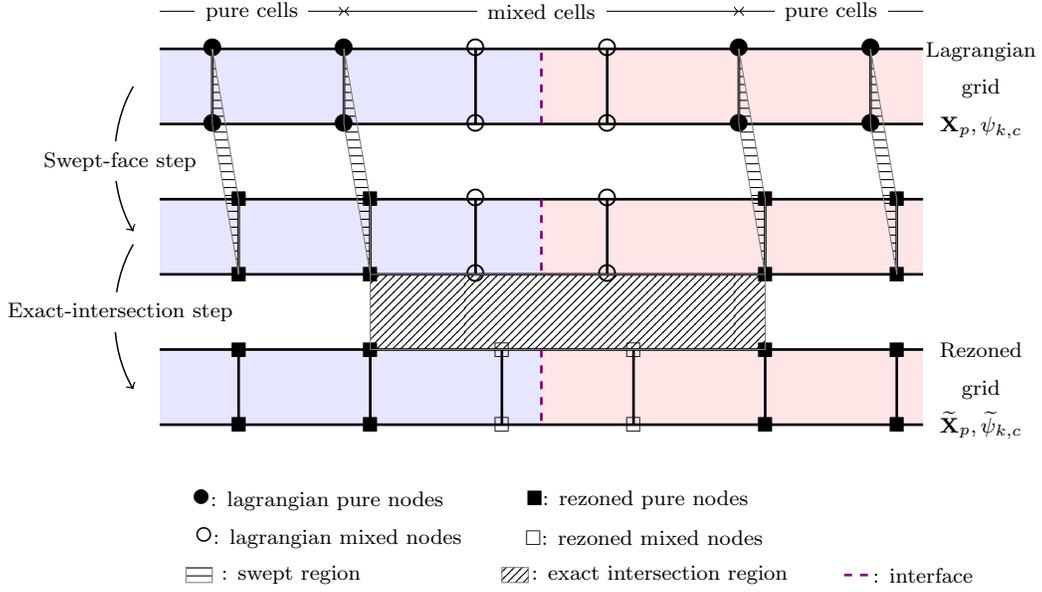

In this part, we detail the hybrid remapping algorithm that is summarized on \figref \ref{fig:8}. To this end, let us introduce $\CN^P$ and $\CN^M$ the sets of nodes and in the same manner $\CC^P$ and $\CC^M$  the sets of cells  respectively used for PCSF and MCIB remapping. Here $\CN^M$ collects mixed nodes belonging to cells that contain the interface or are on this interface (white nodes on \figref \ref{fig:8}) despite $\CN^P$ contains the pure ones (black nodes on \figref \ref{fig:8}). In addition, $\CC^M$ is the set of mixed cell that include cells intersected by the interfaces and their neighbors by nodes. Finally, $\CC^P$ contains the cells that have at least one node in $\CN^P$.\\
The hybrid remapping procedure consists in performing the following steps.
\begin{enumerate}
\item {\bf PCSF step.} In this step 
 we first move the pure nodes included in $\CN^P$ and we remap the quantities in cell $c$ belonging to $\CC^P$. Thus, we have $\widetilde{\psi}_{k,c} = (\widetilde{\rho}_{k,c}, \widetilde{\rho E}_{k,c}, \widetilde{\rho \VU}_{k,c})$ using relation \eqref{eq:rezpsi} and $\widetilde{m}_{k,c}, \widetilde{\alpha}_{k,c}, \widetilde{\VX}_{k,c}$ for each cell $c \in \CC^p$ .
\item {\bf MCIB step.} Now, the mixed nodes in $\CN^M$ are moved and the $\widetilde{\psi}_{k,c} = (\widetilde{\rho}_{k,c}, \widetilde{\rho E}_{k,c}, \widetilde{\rho \VU}_{k,c})$ are remapped thanks to  \eqref{eq:rezpsik} and $\widetilde{m}_{k,c}, \widetilde{\alpha}_{k,c}, \widetilde{\VX}_{k,c}$ are computed for cells $c \in \CC^M$. 
\end{enumerate}
Since $\CC^M\cap\CC^P \neq \{ \emptyset \}$,  one should note that cell included in this intersection are remapped at each step of the algorithm.\\ 

At the end of remapping, only the partial values of the physical variables per unit of volume are known. A this step, a first point is to compute the physical variables per unit of mass. The remapped partial total energy is given using $ \widetilde{E}_{k,c} = \widetilde{(\rho E)}_{k,c}/\widetilde{\rho}_{k,c}$. However, this is different for the remapped partial velocity $\widetilde{\VU}_{k,c}$. Indeed, as explained in the second part of this paper, the Lagrangian computation of the velocity is done in Cartesian geometry. For this reason, the remapped velocity is deduced from the $\widetilde{(\rho \VU)}_{k,c}$ through $\widetilde{\VU}_{k,c} = \widetilde{(\rho \VU)}^{pl}_{k,c}/\widetilde{\rho}^{pl}_{k,c}$ using the {\it planar} remapped density and momentum given through \eqref{eq:rezpsik} and \eqref{eq:rezpsi} with $\CR =1$. The second point is dedicated to the reconstruction of the global values required for the next Lagrangian step. To this end, a classical procedure is to use specific averages 
\begin{equation}
 \widetilde{\phi}_c = \dfrac{1}{ \widetilde{m}_c} \sum_k \widetilde{m}_{k,c} \widetilde{\phi}_{k,c},
\end{equation}
with the global mass and density deduced from
\begin{equation}
\widetilde{m}_c =\sum_k \widetilde{m}_{k,c} \widetilde{ \alpha}_{k,c} \text{ and } 
\widetilde{\rho}_c =\sum_k \widetilde{\rho}_{k,c} \widetilde{ \alpha}_{k,c}.
\end{equation}
At last, thermodynamical variables as pressure $P$ and internal energy $\varepsilon$ are obtained thanks to specific thermodynamical closures as done in \cite{Galera2}.

\section{Numerical results}

We present in this section several numerical test cases performed using the CCALE-MOF computing procedure detailed in \cite{Galera1,Galera2} including the various development proposed in this paper. In the sequel, all the materials are governed by perfect gas equation of state $p=\rho e(\gamma-1)$, where $\gamma$ stands for the polytropic index of gas.

\subsection{Sedov problem}

\begin{figure}[h!]
\centering
\includegraphics[scale=0.45,clip,trim=0cm 0cm 1.6cm 0cm]{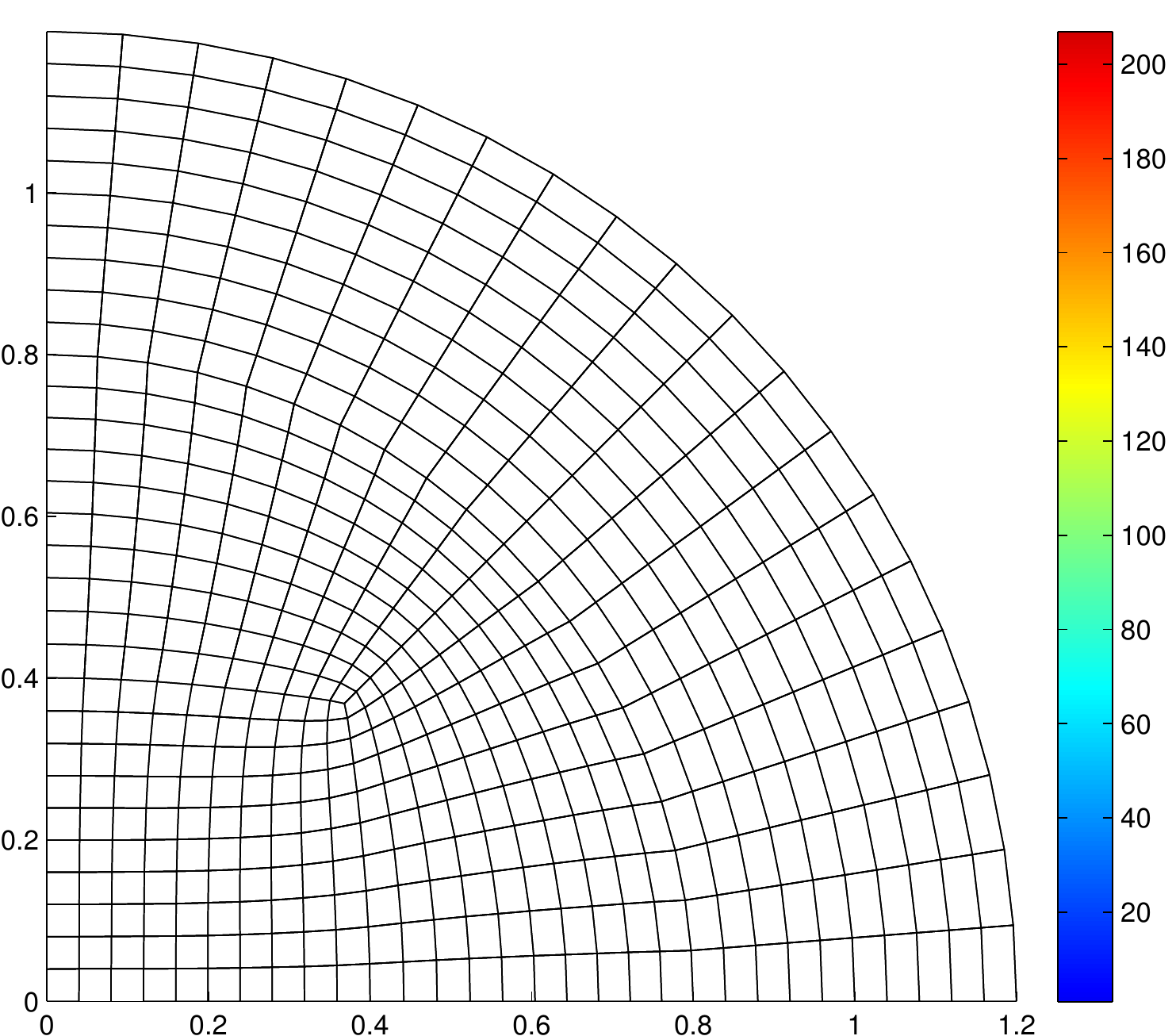}
\includegraphics[scale=0.3]{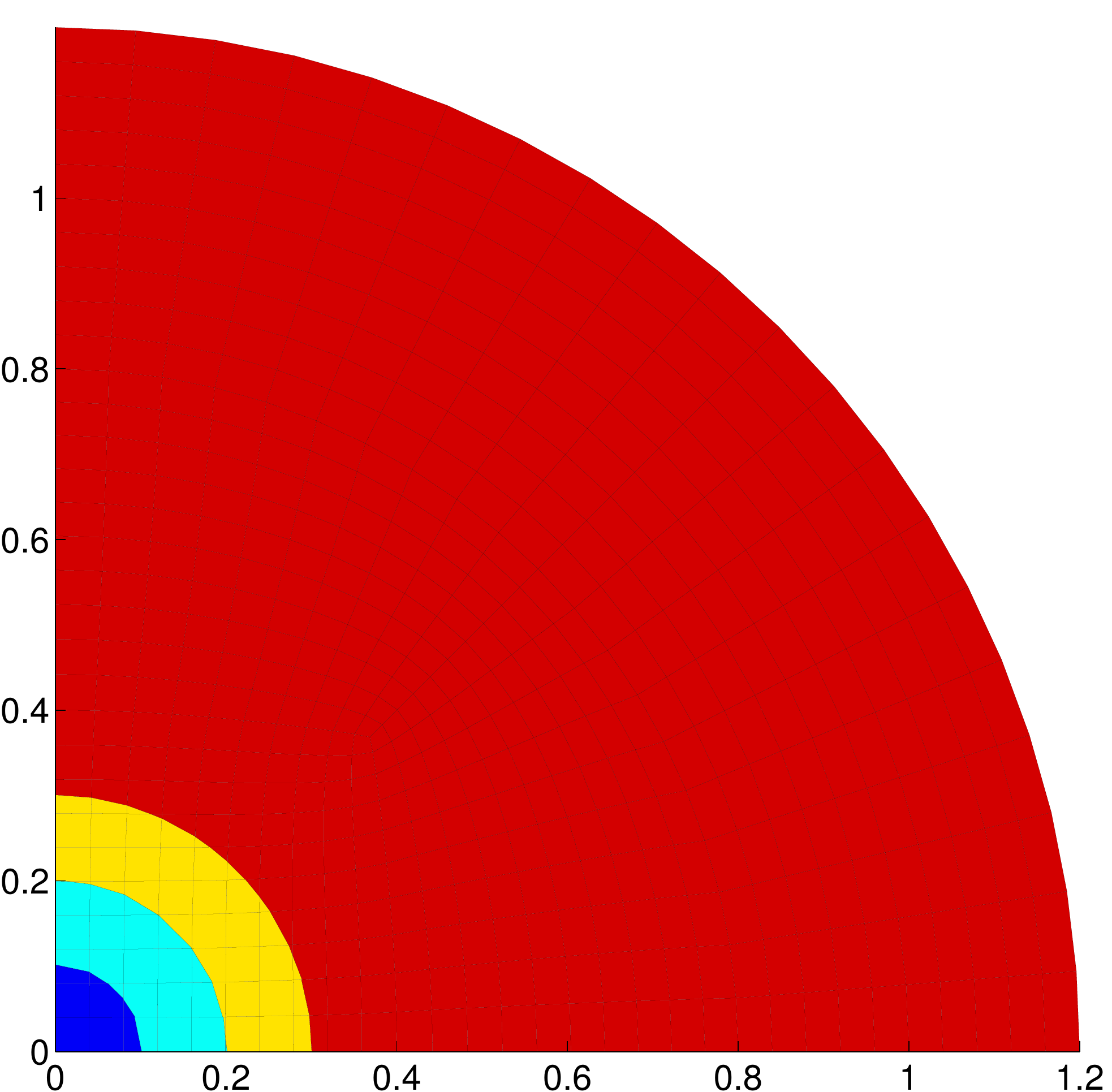}
\caption{Initial grid and material positions for the Sedov problem. \label{sedov:1}}
\end{figure}

\begin{figure}[H]
\centering
\begin{tabular}{cll} 
{} &
 {~~~~New Lagrangian scheme} & 
  {~~~Original EUCCLHYD scheme } \\
\begin{sideways}{~~~~~ Interface positions } \end{sideways} &
\includegraphics[scale=0.25]{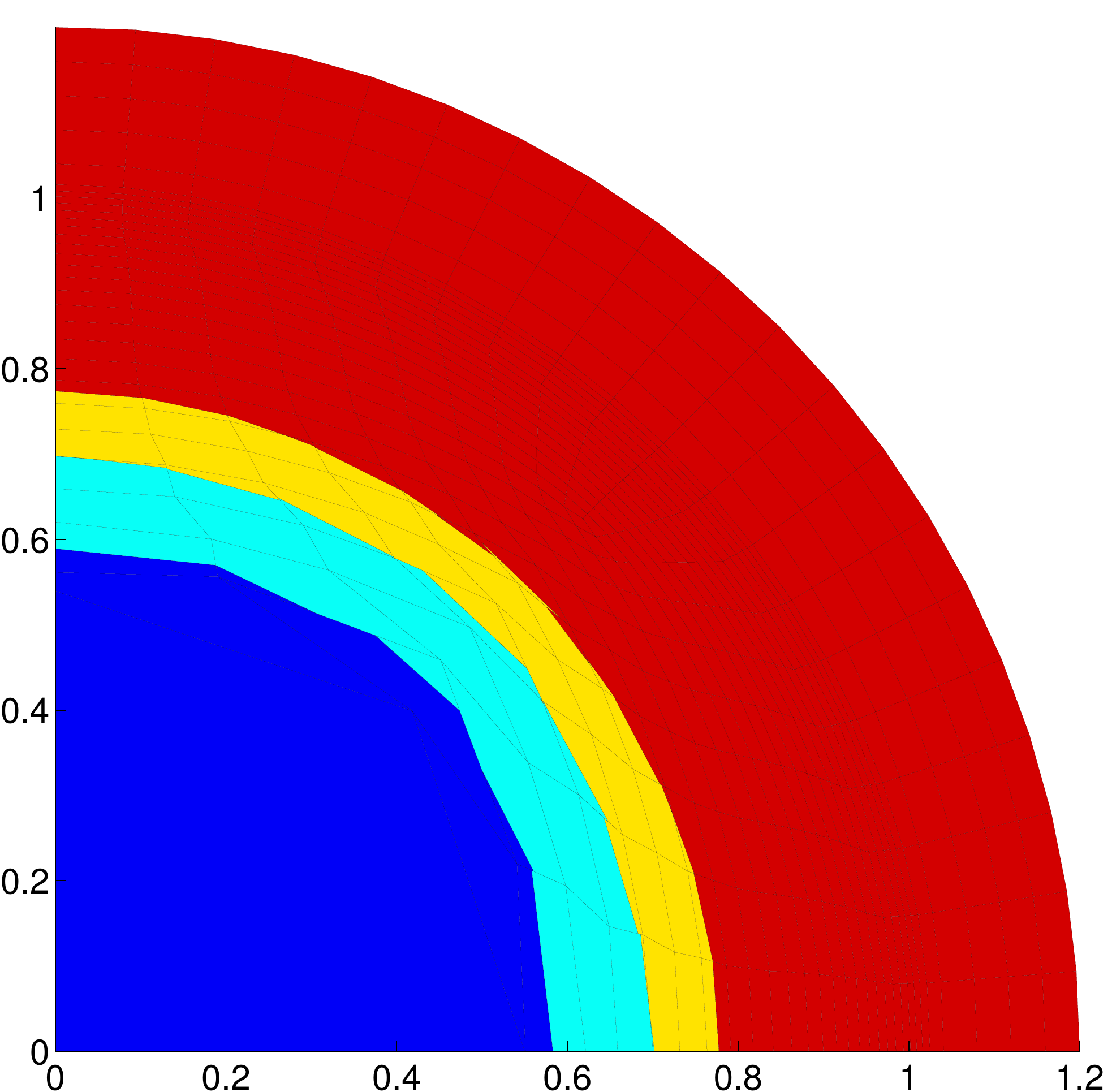} & 
\includegraphics[scale=0.25]{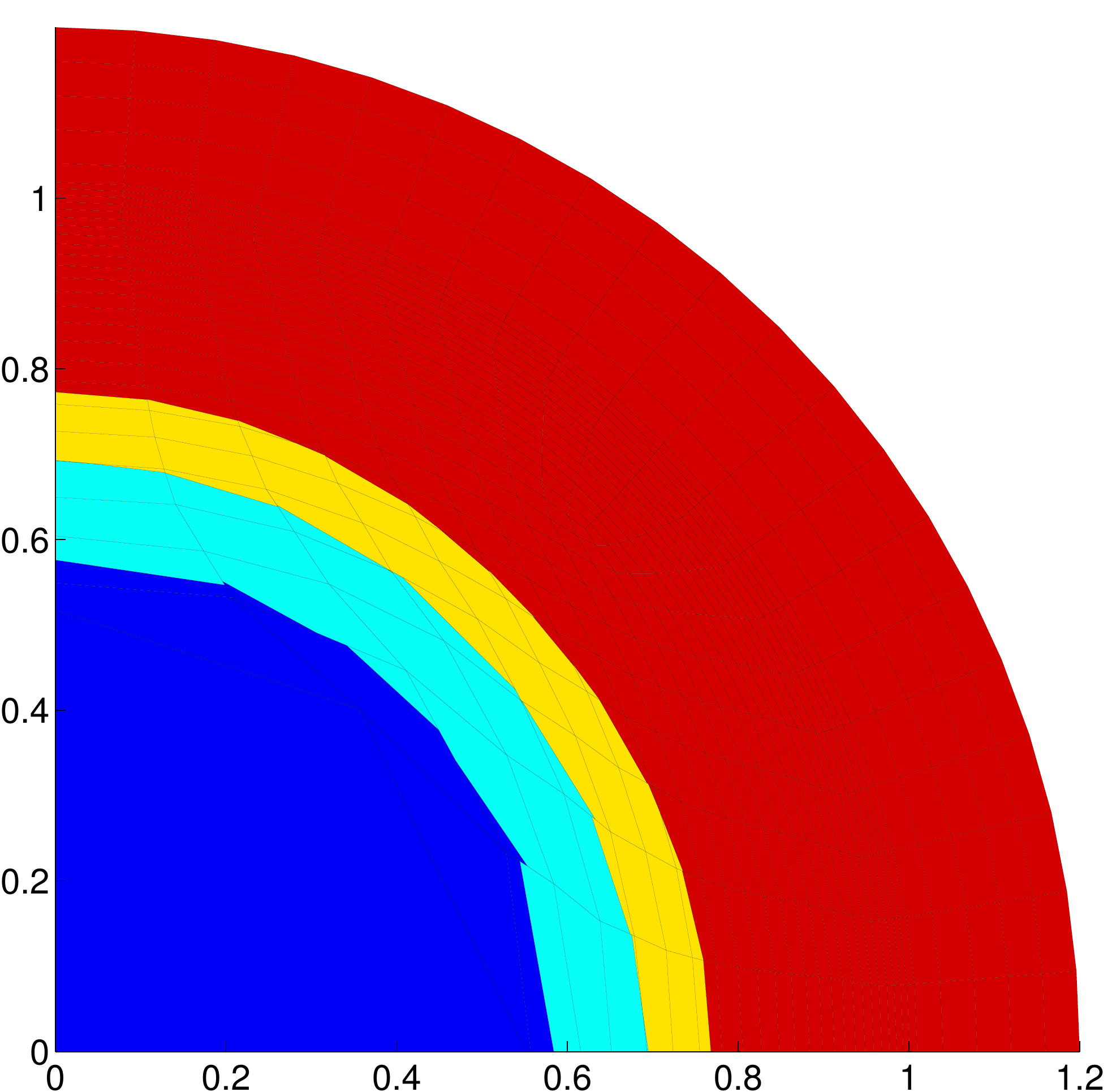} \\
\begin{sideways}{~~~~~~~~ Density} \end{sideways} &
\includegraphics[scale=0.38]{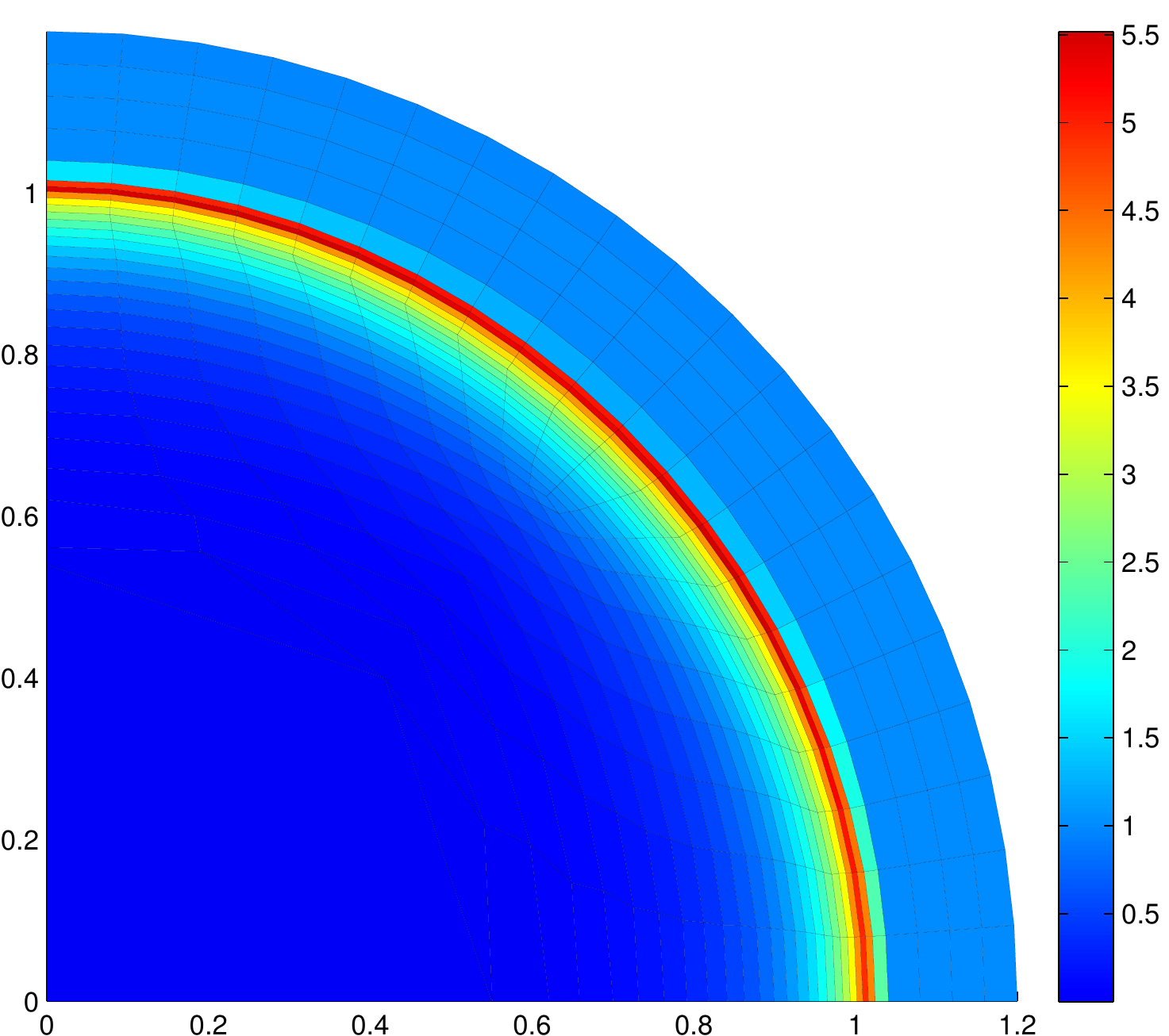} & 
\includegraphics[scale=0.38]{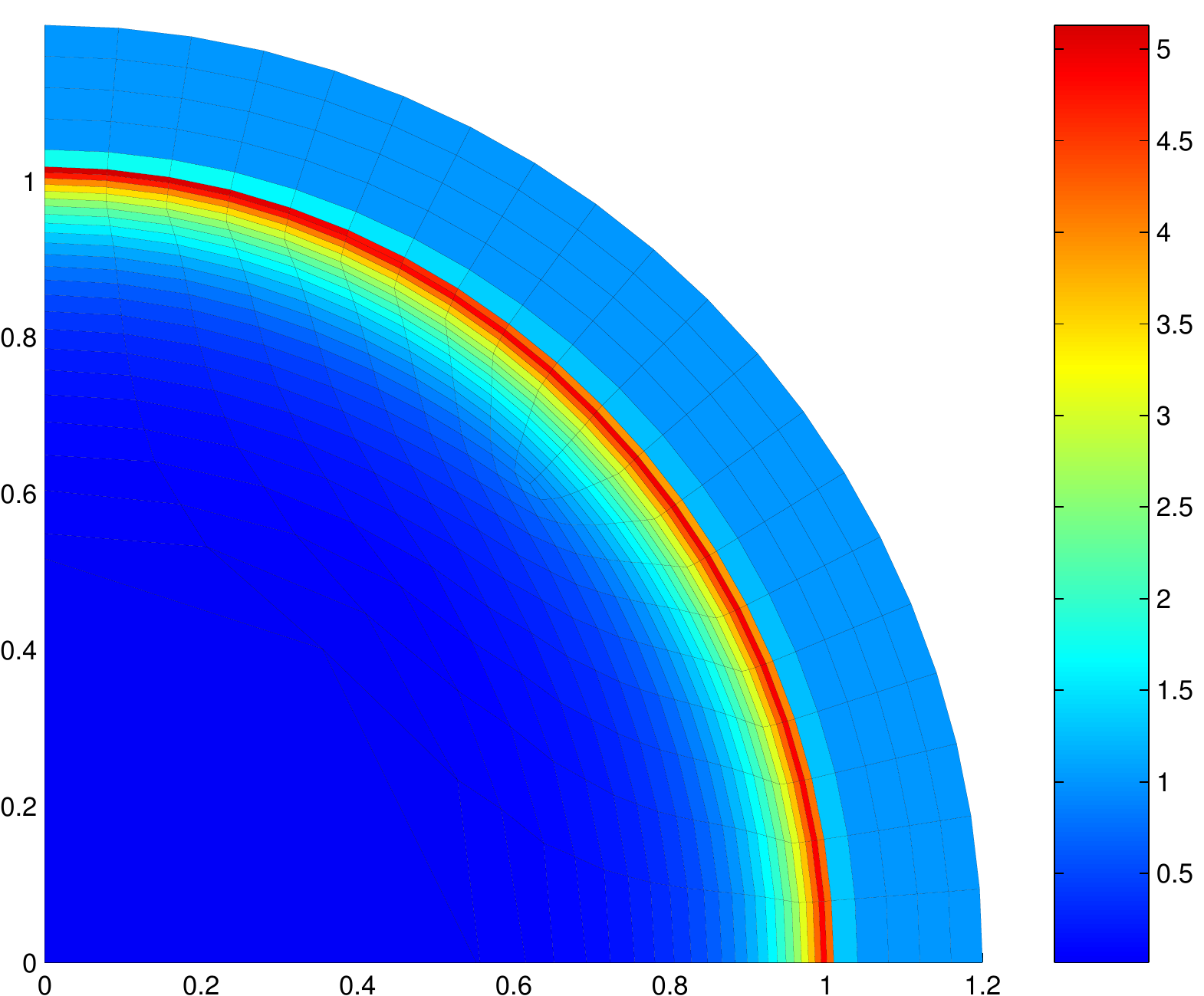} \\
\begin{sideways}{~~~ Density profile} \end{sideways} &
\includegraphics[scale=0.35]{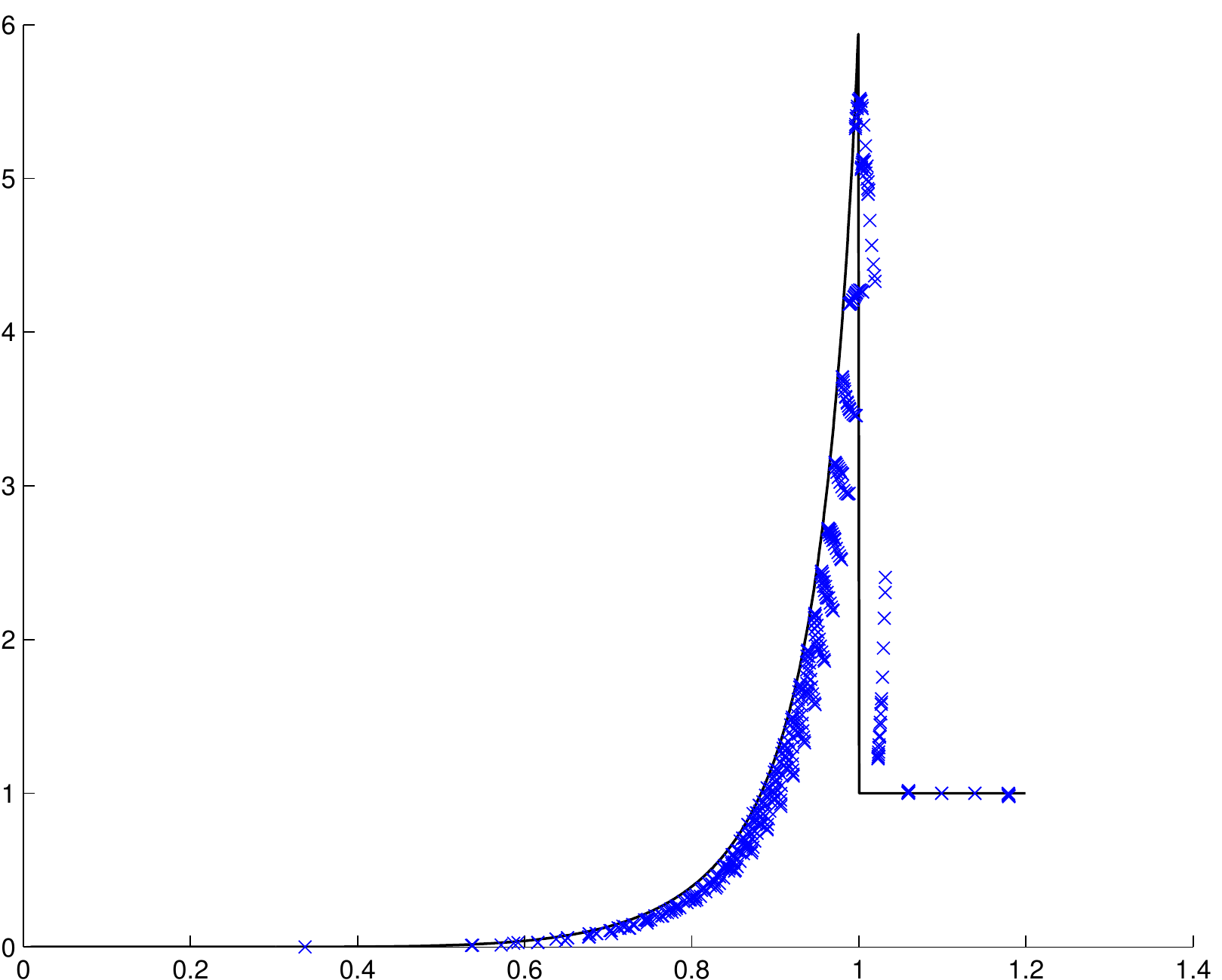} & 
\includegraphics[scale=0.35]{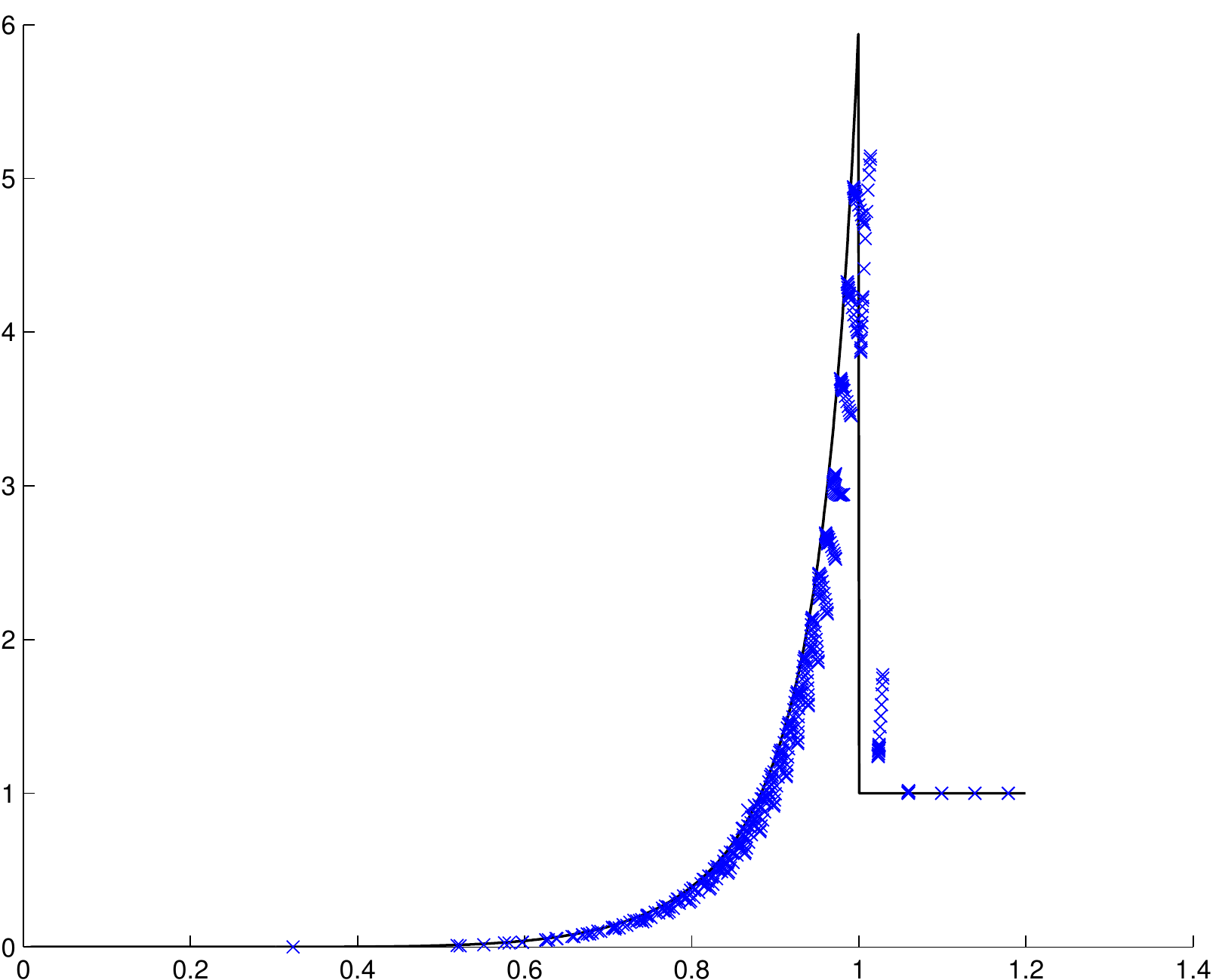} 
\end{tabular}
\caption{Sedov problem. From the top to the bottom: Interface positions, density maps, density profiles defined as a function of the cell center radius compared to the analytical solution at final time step for the pure Lagrangian computation using new scheme (left) and the original EUCCLHYD scheme (right). \label{sedov:2}}
\end{figure}

\begin{figure}[H]
\centering
\begin{tabular}{cll} 
{} &
 {~~~~~~~~New CCALE-MOF } & 
  {~~~~ EUCCLHYD CCALE-MOF } \\
\begin{sideways}{~~~~~ Interface positions } \end{sideways} &
\includegraphics[scale=0.25]{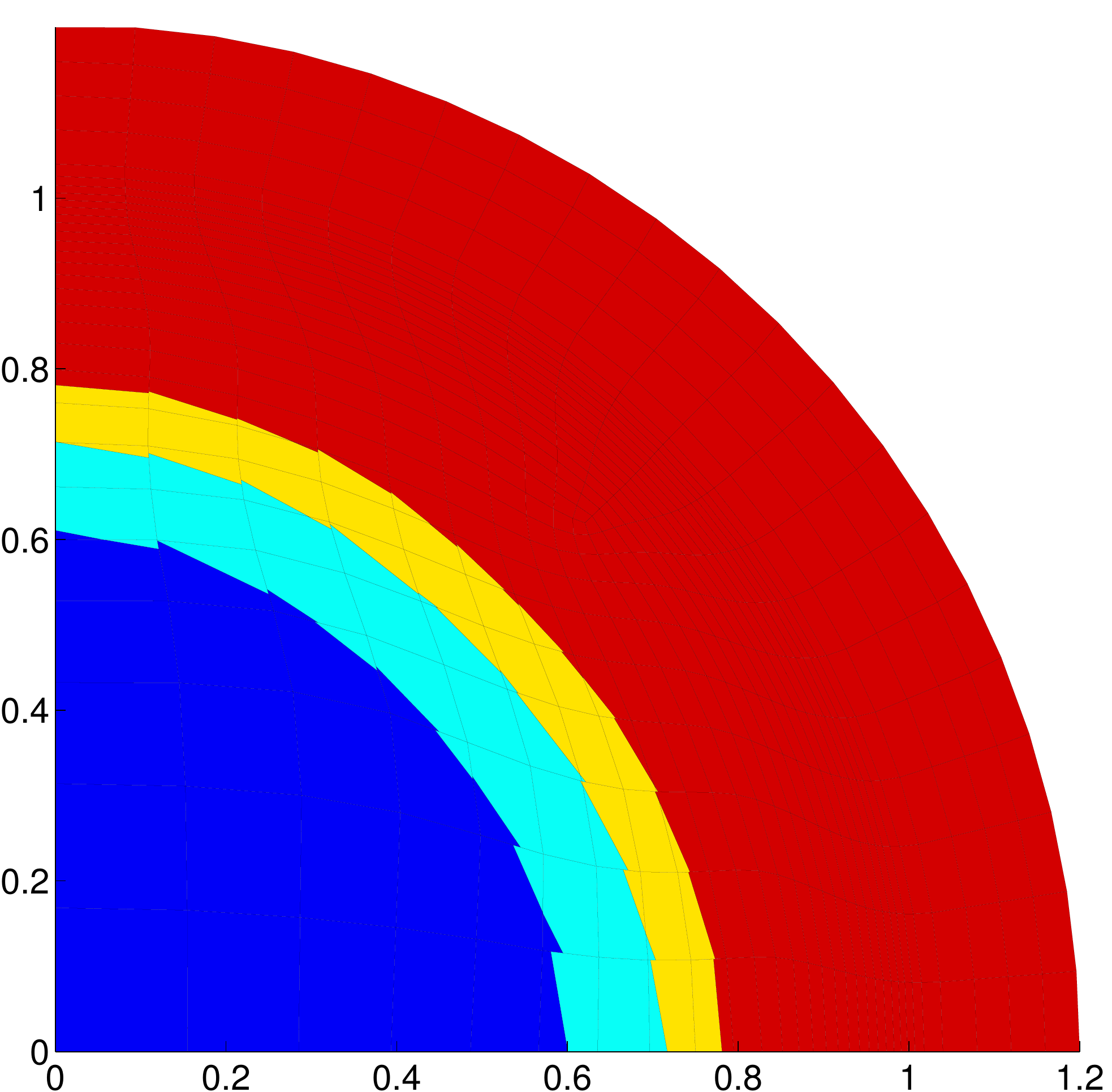} & 
\includegraphics[scale=0.25]{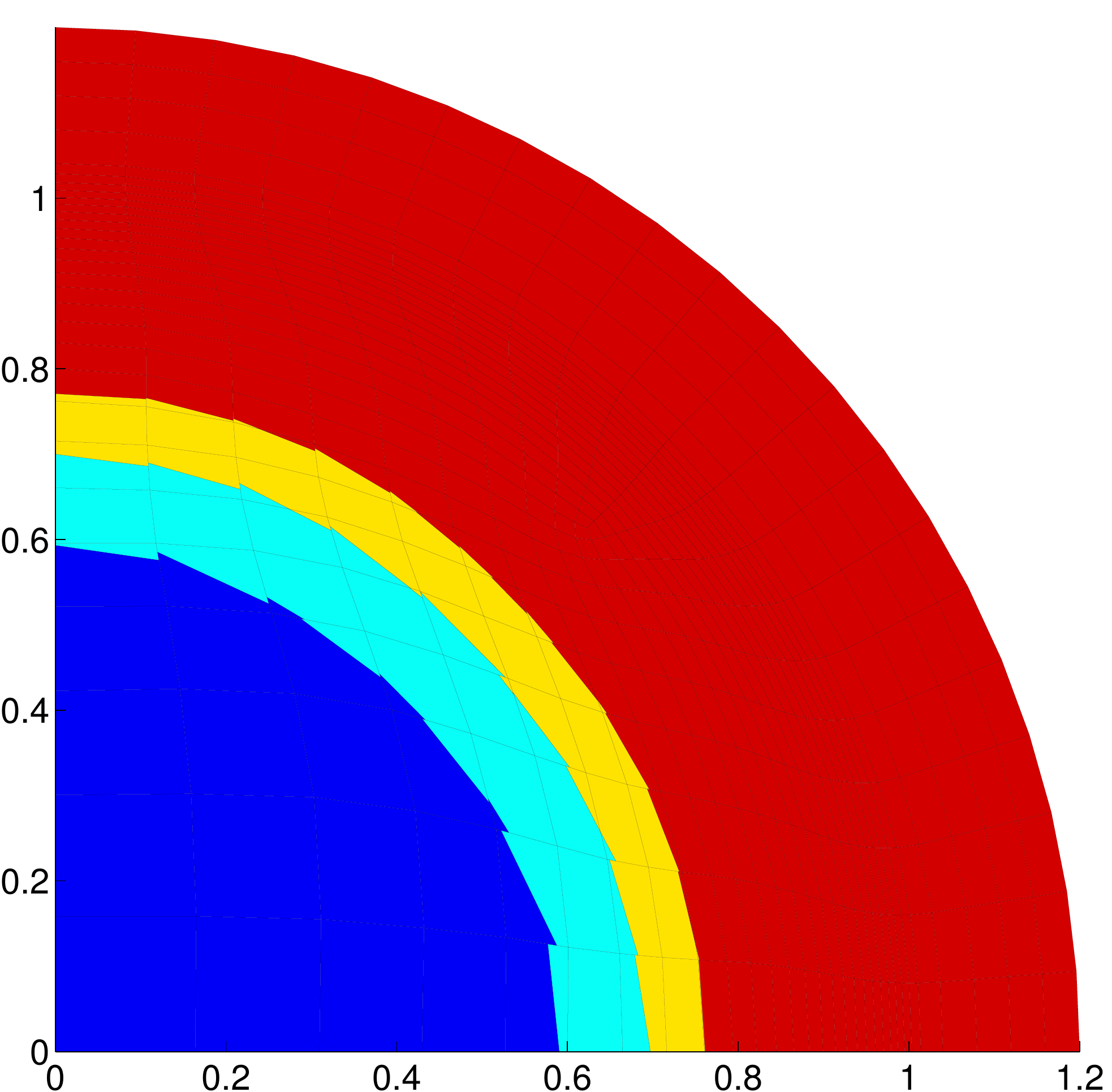} \\
\begin{sideways}{~~~~~~~~ Density} \end{sideways} &
\includegraphics[scale=0.38]{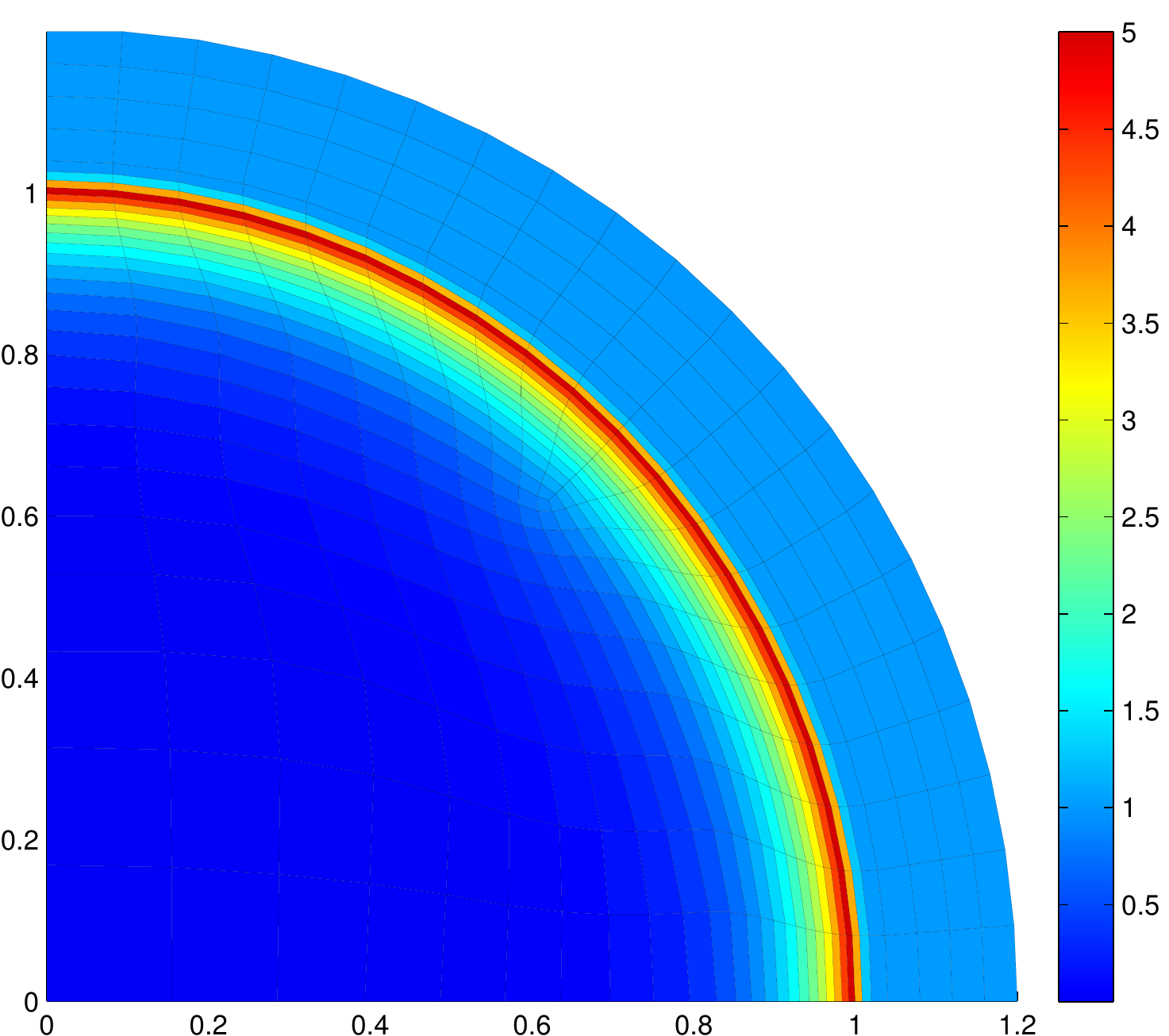} & 
\includegraphics[scale=0.38]{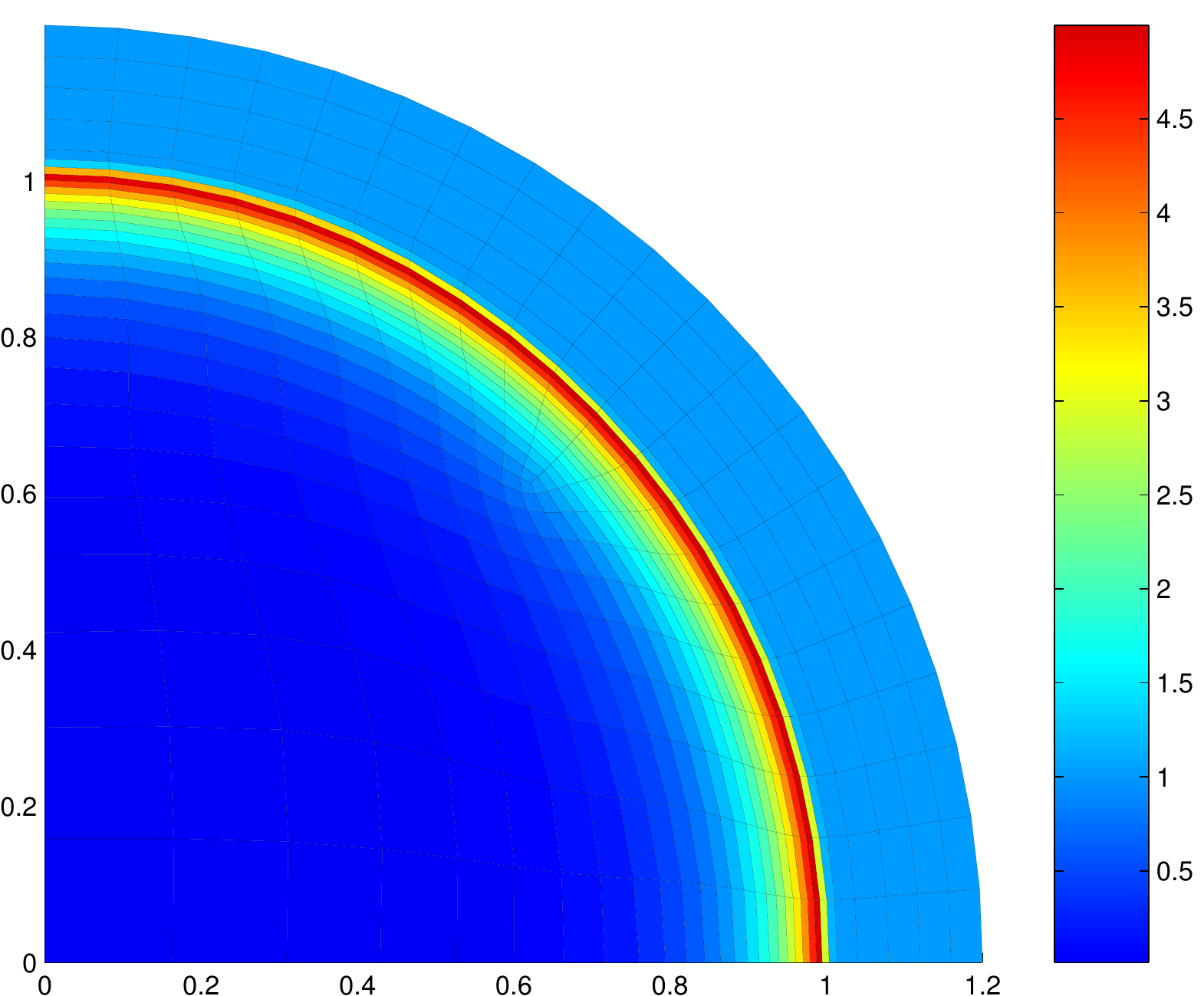} \\
\begin{sideways}{~~~ Density profile} \end{sideways} &
\includegraphics[scale=0.35]{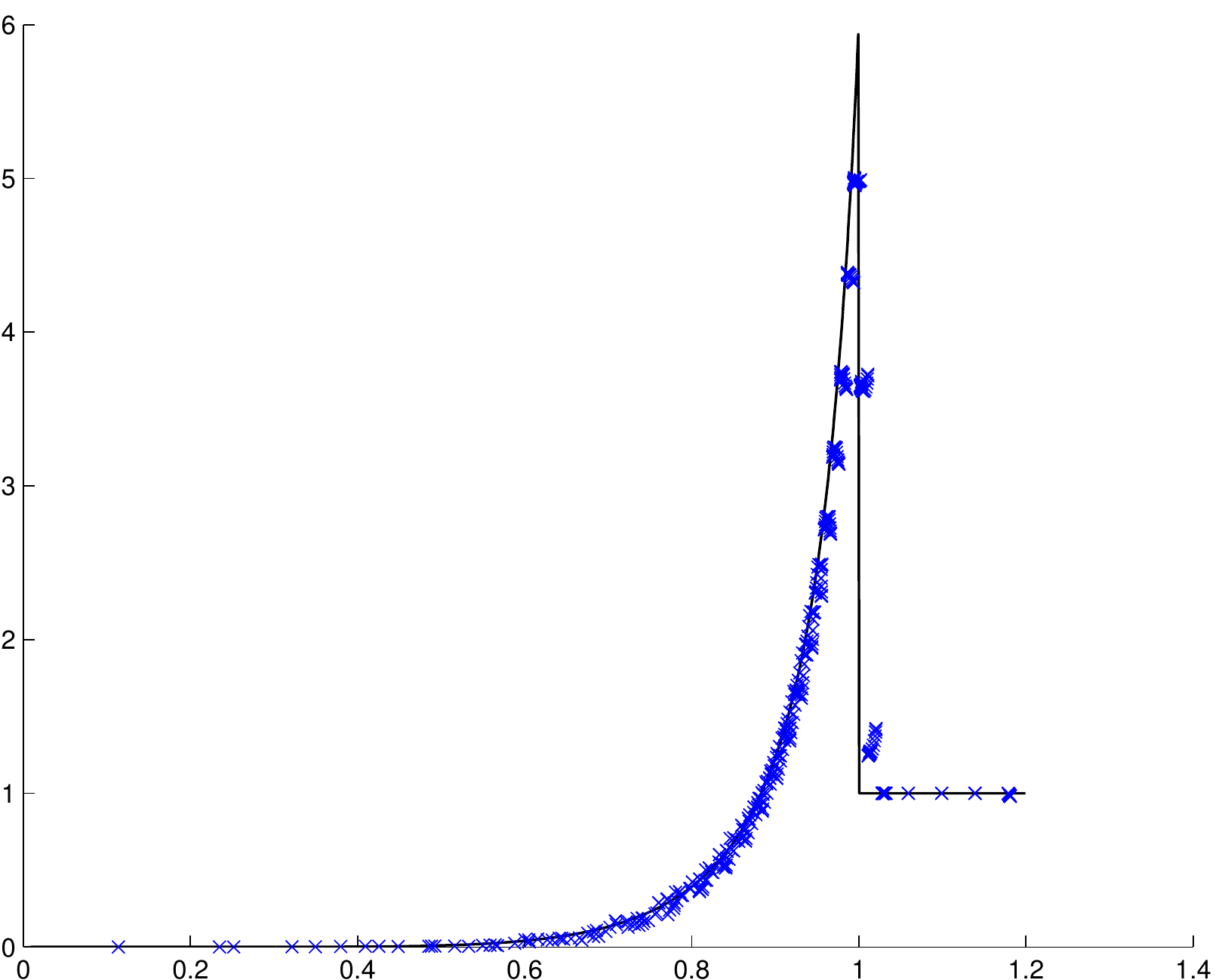} & 
\includegraphics[scale=0.35]{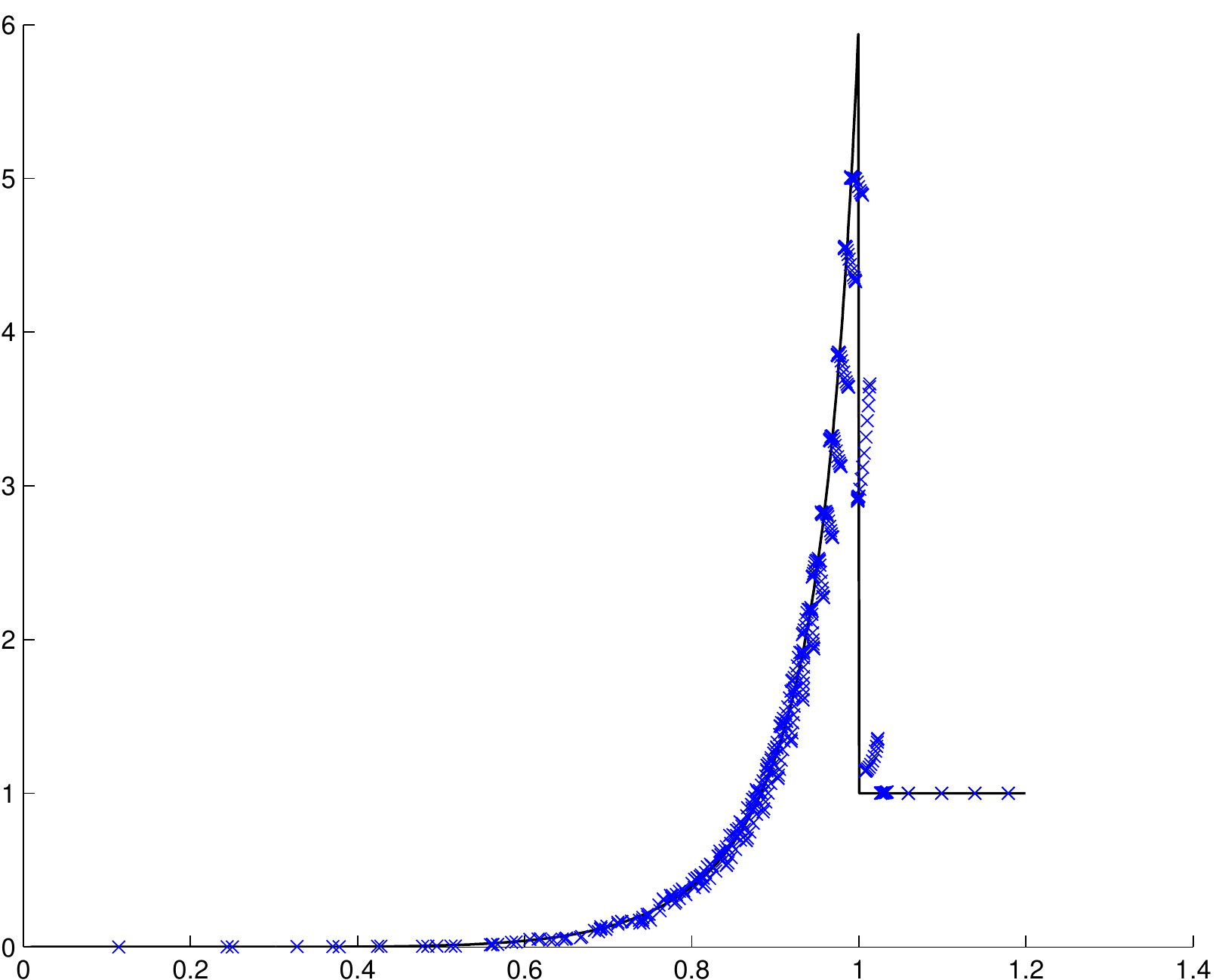} 
\end{tabular}
\caption{Sedov problem. From the top to the bottom: Interface positions, density maps, density profiles defined as a function of the cell center radius compared to the analytical solution at final time step for the new CCALE-MOF procedure (left) and the EUCCLHYD CCALE-MOF procedure (right). \label{sedov:3}}
\end{figure}

We present in this first section a Sedov problem for a point blast in a uniform medium with spherical symmetry. We use this test case to compare our new formulation with the original EUCCLHYD scheme in pure Lagrangian and coupled to the CCALE-MOF procedure. The initial conditions are  given by $(\rho^0,P^0,\VU^0) = (1, 10^{-6}, \boldsymbol{0})$ in a spherical domain of radius $1.2$ except in the cell at the origin $(0,0)$ where an initial delta-function  energy source is set through the pressure
$$
P_{or} = (\gamma-1)\rho_{or} \dfrac{{\cal E}_0}{V_{or}},
$$
with $V_{or}$ the volume of the origin cell and ${\cal E}_0=0.851072$ is the total amount of released energy. The fluid has its polytropic index $\gamma$ equal to $\dfrac{7}{5}$. Contrary to the original single material test case, we add here three artificial interfaces, to test our multi-material CCALE-MOF algorithm. These interfaces are initially located for a radius equals to $0.1$, $0.2$ and $0.3$ (see \figref \ref{sedov:1}).\\
Here we consider both Lagrangian and ALE computations for an initial unstructured mesh depicted on \figref \ref{sedov:1}. This grid is obtained after one rezoning step, with $\omega_p=1$ of an unstructured mesh initially paved with $500$ quadrangular cells. Numerical results are depicted on \figref \ref{sedov:2} and\figref \ref{sedov:3}  for a final time of $t_{end} = 1$ and compared to the analytical solution computed using self-similar arguments as done in \cite{Galera2}. It consists of a diverging shock wave whose front is exactly localized at radius $R=1$. As it is illustrated on \figref \ref{sedov:2}, the pure Lagrangian solutions are in good agreement with the  analytical one for both approaches. We can notice that the new formulation is less dissipative as we reach a higher density level in the shock region. Indeed for the Lagrangian method as for the CCALE-MOF one the shock location is well resolved without any spurious oscillation (\figref \ref{sedov:3}). In addition, this simple problem underlines the robustness (better mesh quality near the origin) and accuracy (shock location) of the axisymmetric CCALE-MOF approach especially when considering multi-material flows whose interfaces are well captured thanks to the MOF reconstruction (see \figref \ref{sedov:3}).

We point out that during the Lagrangian computation, non-convex cells appeared. This may lead to interface reconstruction 
failure when considering multi-material flows.  As illustrated by the previous numerical results, the proposed CCALE-MOF algorithm remains adapted to treat such configuration without any difficulty demonstrating once again its robustness.
\subsection{Axi-symmetric triple point problem}

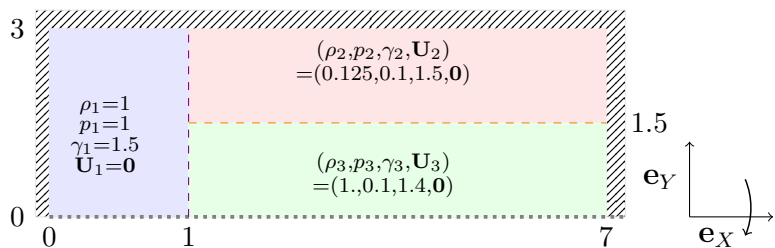
\begin{figure}[h!]
\centering
 \begin{tikzpicture}[xscale=0.37,yscale=0.5]
 \draw[gray!20,pattern= north east lines]  (-0.5,0) rectangle (20.7,5.5);
  \filldraw[fill=red!10,draw=red!10] (0,0) rectangle (20,5);
  \filldraw[fill=blue!10,draw=blue!10] (0,0) rectangle (5,5); 
  \filldraw[fill=green!10,draw=green!10]   (5,0) rectangle (20,2.5); 
  \draw[gray,dotted,line width=0.5mm]  (0,0)--(20.7,0);
  \draw[dashed,violet] (5,0) -- (5,5);
  \draw[dashed,orange] (5,2.5) -- (20,2.5);
  \draw (2,3.5) node[below]{$\substack{\rho_1 = 1\\ p_1 = 1 \\ \gamma_1 = 1.5 \\ {\VU}_1 = \boldsymbol{0} }$}; 
  \draw (12,5)node[below]{$\substack{(\rho_2,p_2,\gamma_2, {\VU}_2) \\ = (0.125,0.1,1.5,\boldsymbol{0})}$};
  \draw (12,2)node[below]{$\substack{(\rho_3,p_3,\gamma_3, {\VU}_3) \\ = (1.,0.1,1.4,\boldsymbol{0})}$};
  \draw (-0.5,5) node [left]{\small $3$};
  \draw (20,0) node [below]{\small $7$};
  \draw (20.5,2.5) node [right]{ \small $1.5$};
  \draw (5,0) node [below]{\small $1$};	
  \draw (0,0) node [below]{\small $0$};
  \draw (-0.5,0) node [left]{\small $0$};
  \draw[<-,line width=0.2mm] (25,-0.5) arc (-30:30:1.5);
  \node[left] at (23,1) {${\bf e}_Y$};
\node[below] at (24,0) {${\bf e}_X$};
\draw[->] (23,0)--(23,2); 
\draw[->] (23,0)--(26,0);
\end{tikzpicture}
 \caption{Axi-symmetric triple point problem : geometry and initial data.\label{fig:1}}
\end{figure}

We consider in this part a three-material problem that corresponds to a three-state Riemann problem in an axisymmetric geometry. This problem has been wisely studied in Cartesian geometry and here we propose new results for cylindrical geometry. The computational domain is rectangular and composed of three regions (blue, green, red) whose dimensions are depicted on \figref \ref{fig:1}. The top, left and right boundaries are closed thanks to walls. A symmetry condition is applied to the bottom boundary corresponding to the $X$-axis axi-symmetry. Initially, the blue region contains a fluid with  high pressure and density taken equal to $(\rho_1,p_1)= (1,1)$. The green region contains a low density and pressure fluid whose initial state is  $(\rho_2,p_2)= (0.1,0.125)$. The third fluid in the red region,  initially has a low pressure and an high density equal to  $(\rho_3,p_3)= (1,0.1)$.  At the beginning of the computation, all fluids are supposed to be at rest then $\VU_1=\VU_2=\VU_3 = \boldsymbol{0}$. The blue and green material have the same polytropic index $\gamma_1=\gamma_2 =1.5$, despite the red one has $\gamma_3=1.4$.\\

\begin{figure}[h!]
\centering
\includegraphics[scale=0.4]{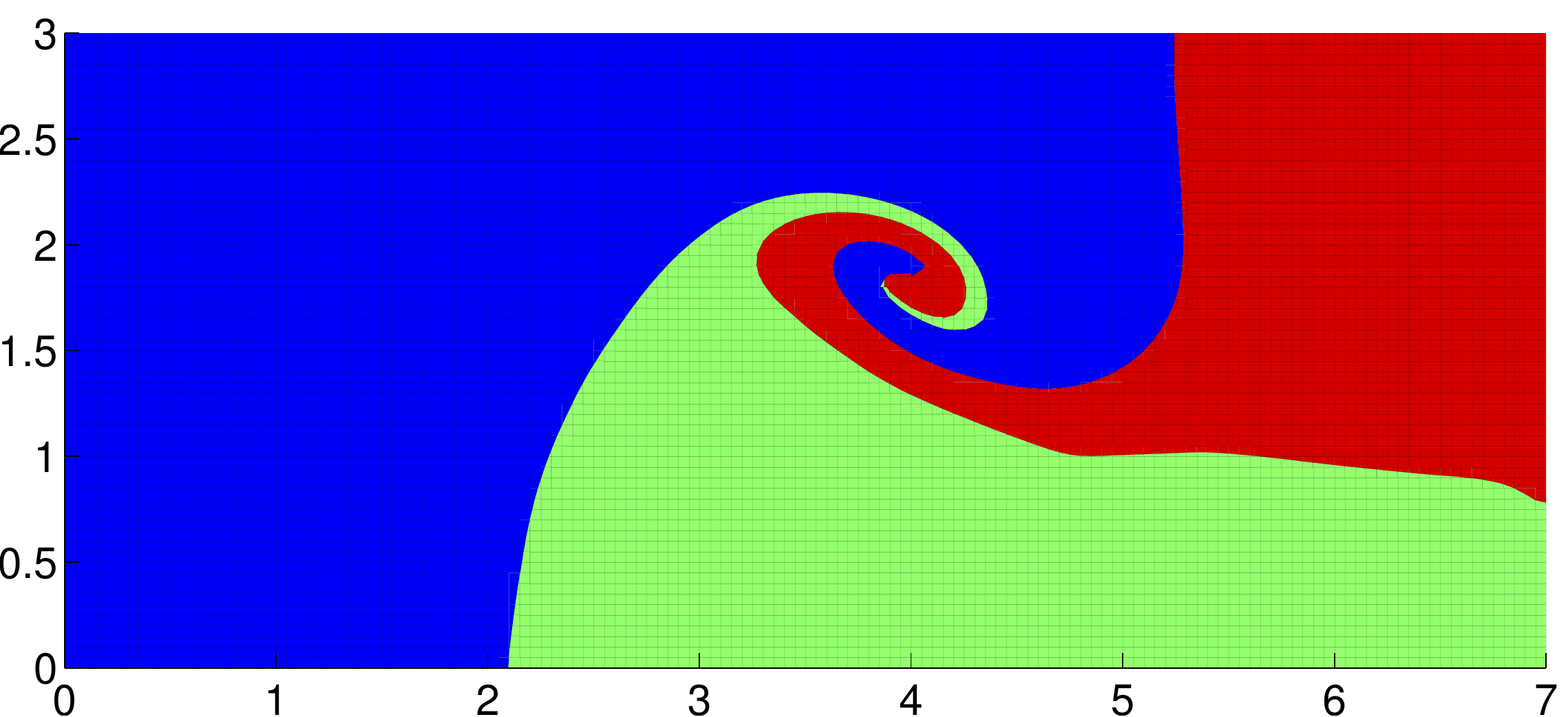} 
\caption{Axi-symmetric triple point problem. Mesh and material positions at $t=5$ for Eulerian computation. \label{pt:1}}
\end{figure}
\begin{figure}[h!]
\centering
\includegraphics[scale=0.4]{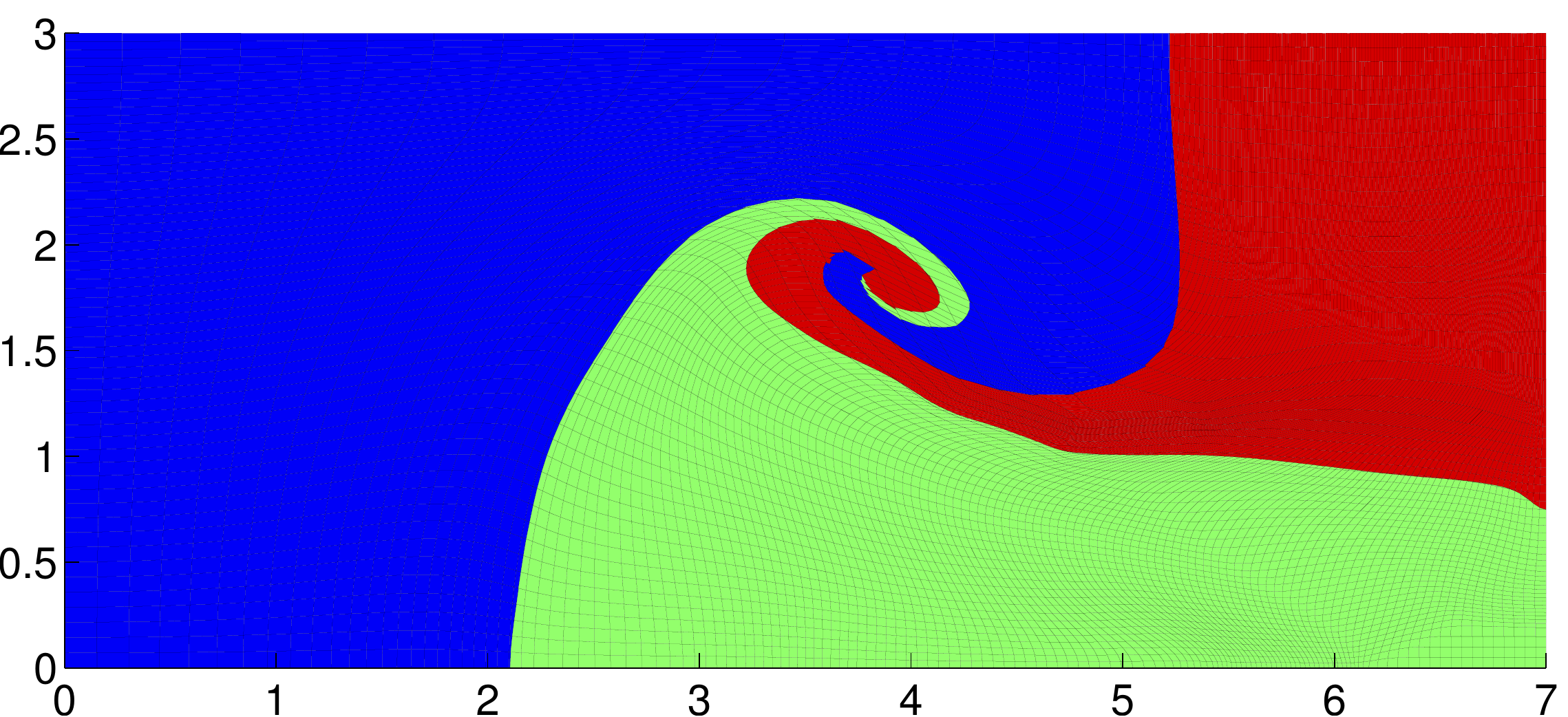} 
\caption{Axi-symmetric triple point problem. Mesh and material positions at $t=5$ for ALE computation. \label{pt:2}}
\end{figure}

The computation using the presented axisymmetric extension of the CCALE-MOF algorithm is made on a grid initially paved with $140\times60$ square cells until a final time $t_f=5$. For this simulation, comparison with a full Lagrangian computation can not be performed since its suffers from important mesh tangling as shown in \cite{Loubere1}. However comparison to full Eulerian simulations is done. In this case, nodes are moved to their initial positions during the rezoning step. Numerical results for both ALE and Eulerian methods representing interfaces and meshes are depicted on \figref \ref{pt:1}-\ref{pt:2}. As expected, since there is a shock wave with high speed that propagates from the heavy material (blue) to the light one (red), the interface is sheared at the triple point producing a Kelvin-Helmholtz like instability. 
Here, comparison to planar 2D computations \cite{Galera1} demonstrates that axisymmetric geometry particularly affects the vortex shape that is 3D. Although the global behavior of the solutions is very similar comparing ALE approach to the Eulerian one.

\subsection{Spherical Air-Helium shock/bubble interaction test}

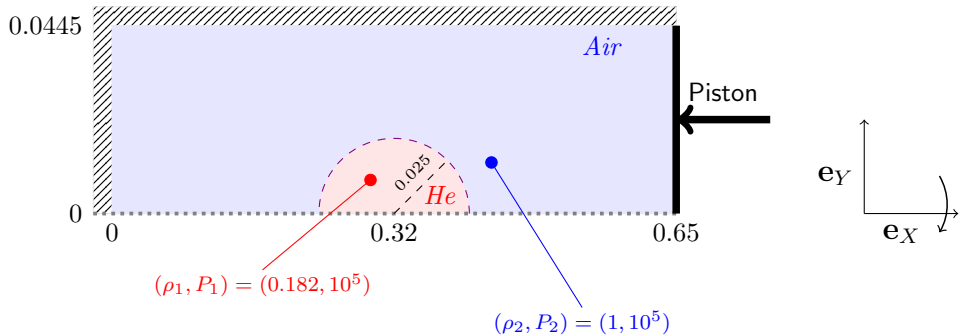
\begin{figure}[h!]
\centering
\begin{tikzpicture}[xscale=0.5,yscale=0.5]
\draw[gray!20,pattern= north east lines]  (-0.5,0) rectangle (15,5.5);
\filldraw[fill=blue!10,draw=blue!10] (0,0) rectangle (15,5);
\fill[red!10](9.5,0) arc (0:180:2);
\draw[gray,dotted,line width=0.5mm]  (-0.5,0)--(15,0);
\draw[violet,dashed](9.5,0) arc (0:180:2);
\draw (0,0) node [below]{\footnotesize  $0$};
\draw (7.5,0)  node[below]{\footnotesize  $0.32$};
\draw (-0.5,5)  node[left]{\footnotesize  $0.0445$};
\draw (-0.5,0)  node[left]{\footnotesize  $0$};
\draw (15,0)  node[below]{\footnotesize $0.65$};
\draw[dashed] (7.5,0)--(8.91,1.41) node[above left,rotate=45]{\tiny ~0.025};
\draw (13,5)node[below,blue]{\footnotesize \it Air};
\draw (8.7,0.5)node[red]{\footnotesize  \it He};
\draw[line width=1mm] (15,0)--(15,5);
\draw[line width=1mm,<-] (15,2.5)--(17.5,2.5);
\draw (16.25,2.7)node[above]{\footnotesize \sf Piston};
\draw[*-,draw=red,fill=red] (7,1)--(4,-1.5); 
\draw (4,-2.5) node [above]{\scriptsize \textcolor{red}{$(\rho_1,P_1)=(0.182,10^5) $}};
\draw[*-,draw=blue,fill=blue] (10,1.5)--(12.5,-2.5); 
\draw (12.5,-3.5) node [above]{\scriptsize \textcolor{blue}{$(\rho_2,P_2)=(1,10^5) $}};
\draw[<-,line width=0.2mm] (22,-0.5) arc (-30:30:1.5);
\node[left] at (20,1) {${\bf e}_Y$};
\node[below] at (21,0) {${\bf e}_X$};
\draw[->] (20,0)--(20,2.5); 
\draw[->] (20,0)--(22.5,0);
\end{tikzpicture}
\caption{Air-Helium shock/bubble interaction. Initial geometry and data. \label{fig:10}}
\end{figure}

\begin{figure}[H]
\centering
\small 
\begin{tabular}{c}
\includegraphics[scale=0.25, clip,trim=0cm 2.2cm 0cm 1cm ]{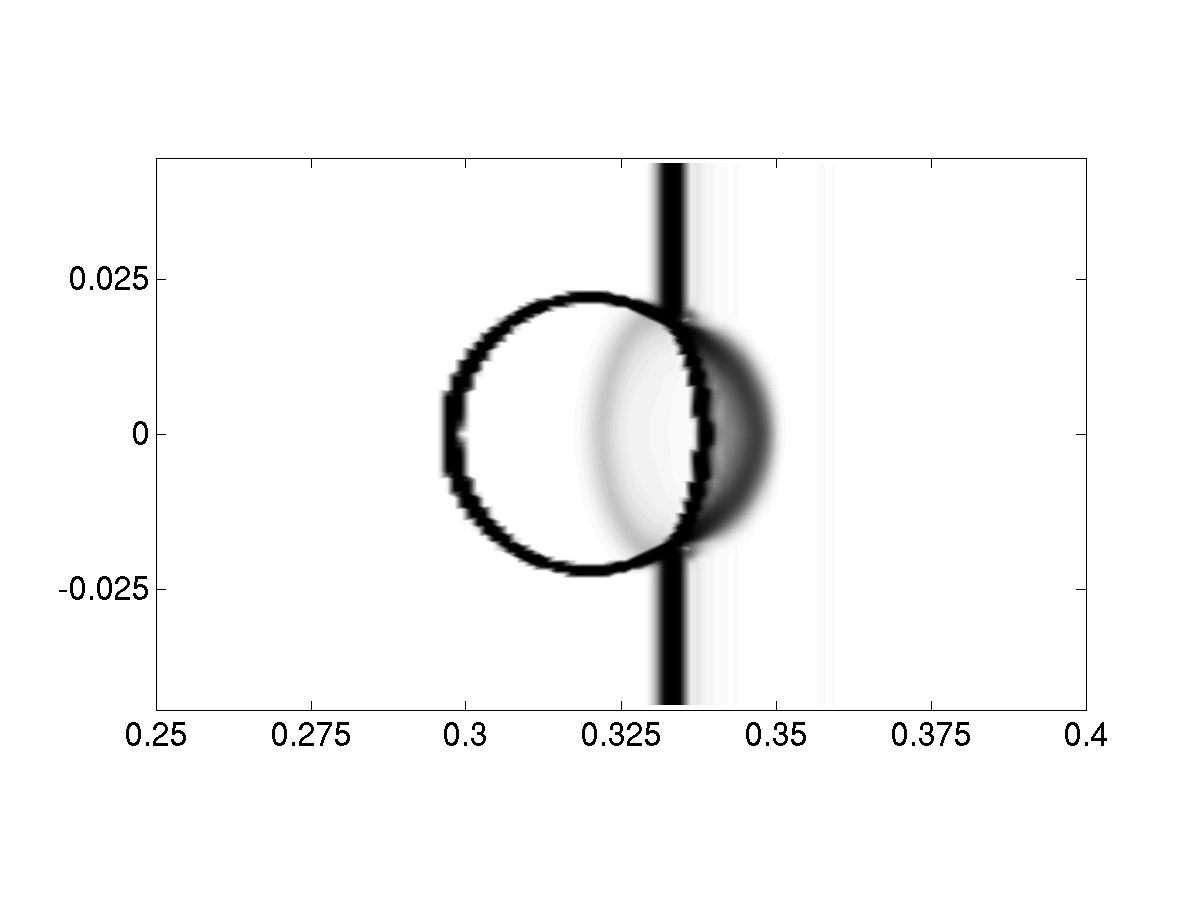}
\includegraphics[scale=0.35]{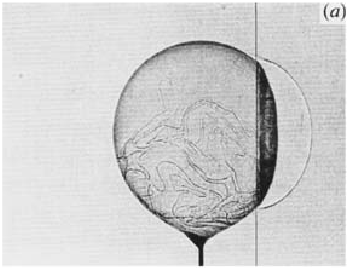} \\[-0.25cm]
$t=t_i+20\times10^{-6}$\\[-0.15cm]
\includegraphics[scale=0.25, clip,trim=0cm 2.2cm 0cm 1cm ]{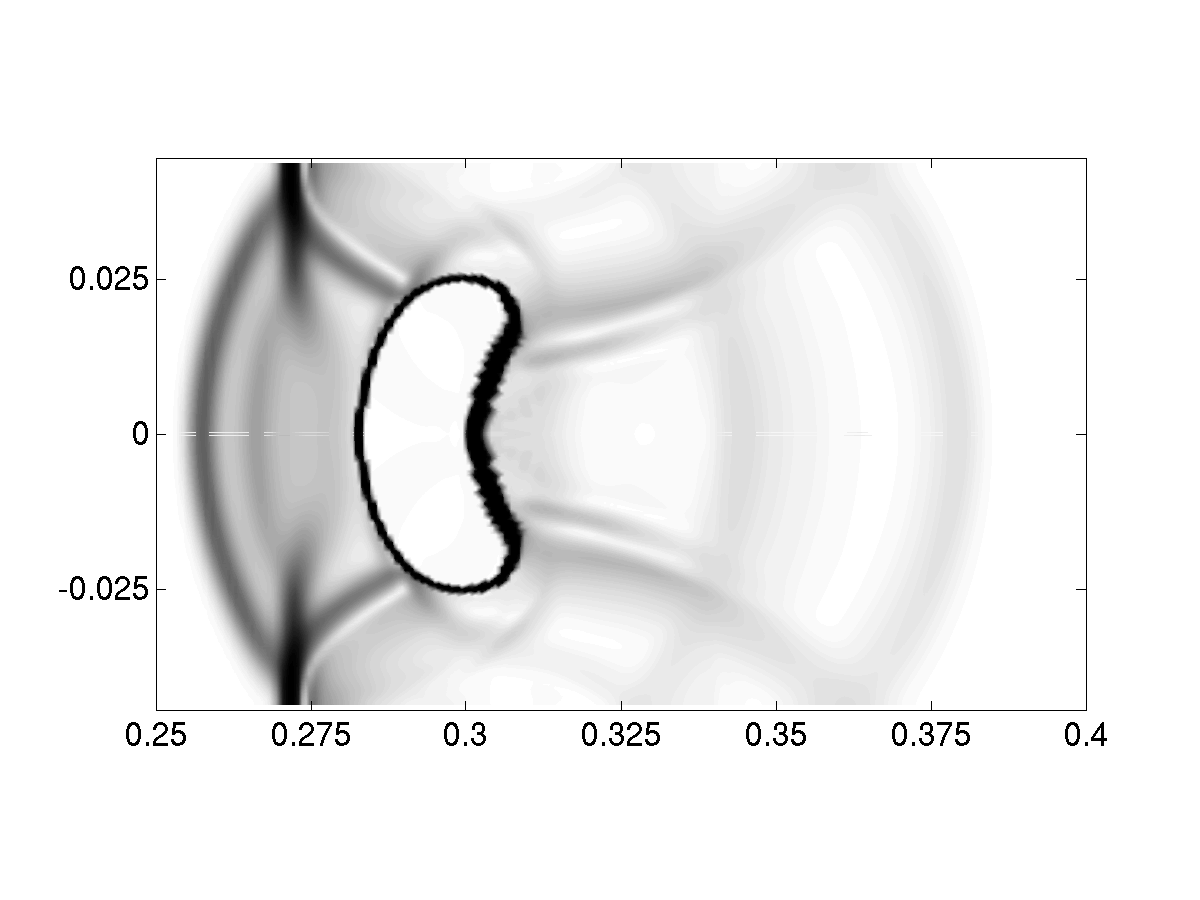}
\includegraphics[scale=0.35]{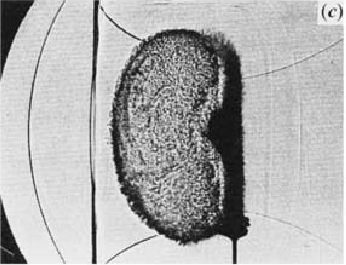} \\ [-0.25cm]
$t=t_i+145\times10^{-6}$\\[-0.15cm]
\includegraphics[scale=0.25, clip,trim=0cm 2.2cm 0cm 1cm ]{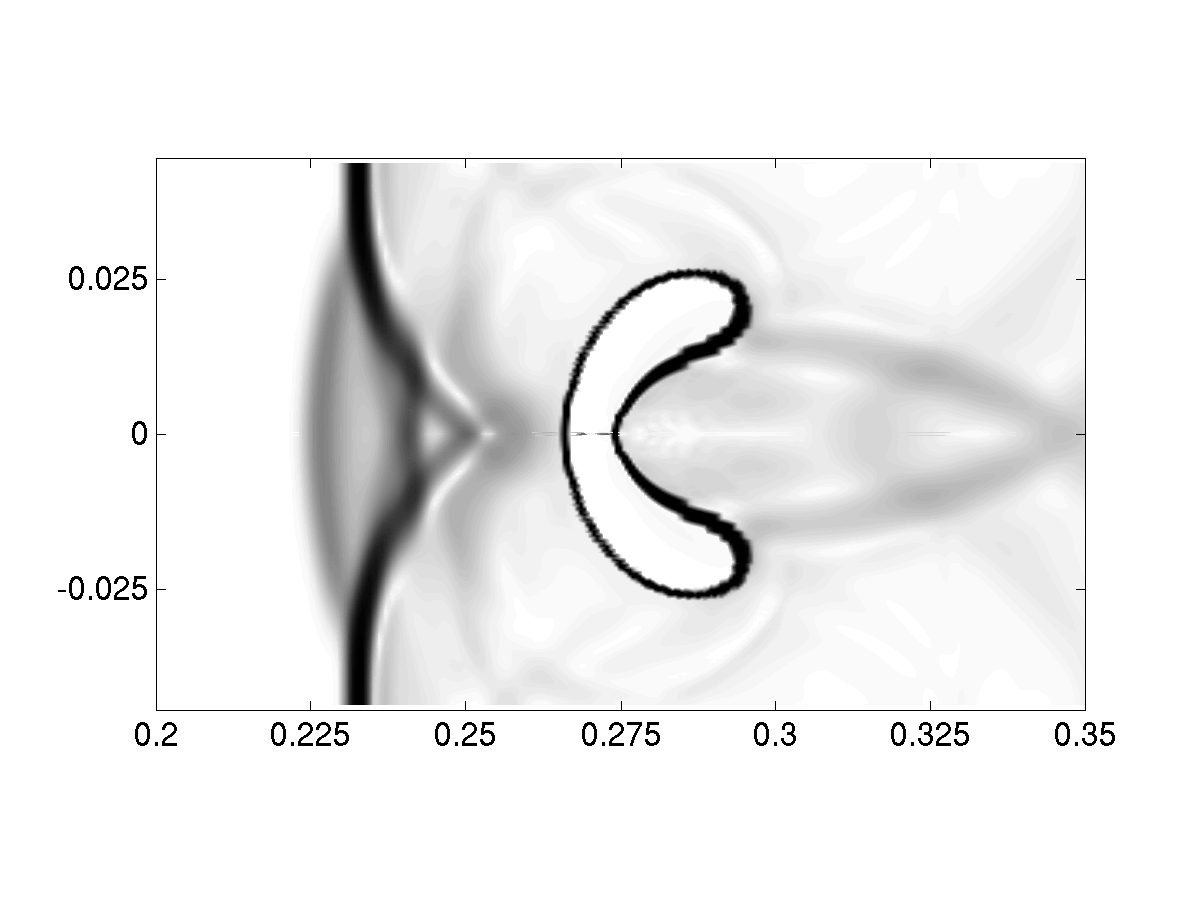}
\includegraphics[scale=0.35]{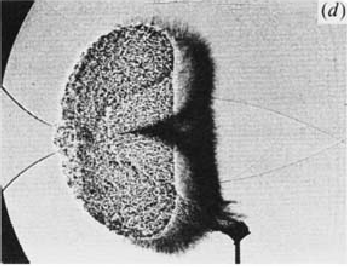} \\[-0.25cm]
$t=t_i+223\times10^{-6}$\\[-0.15cm]
\includegraphics[scale=0.25, clip,trim=0cm 2.2cm 0cm 1cm ]{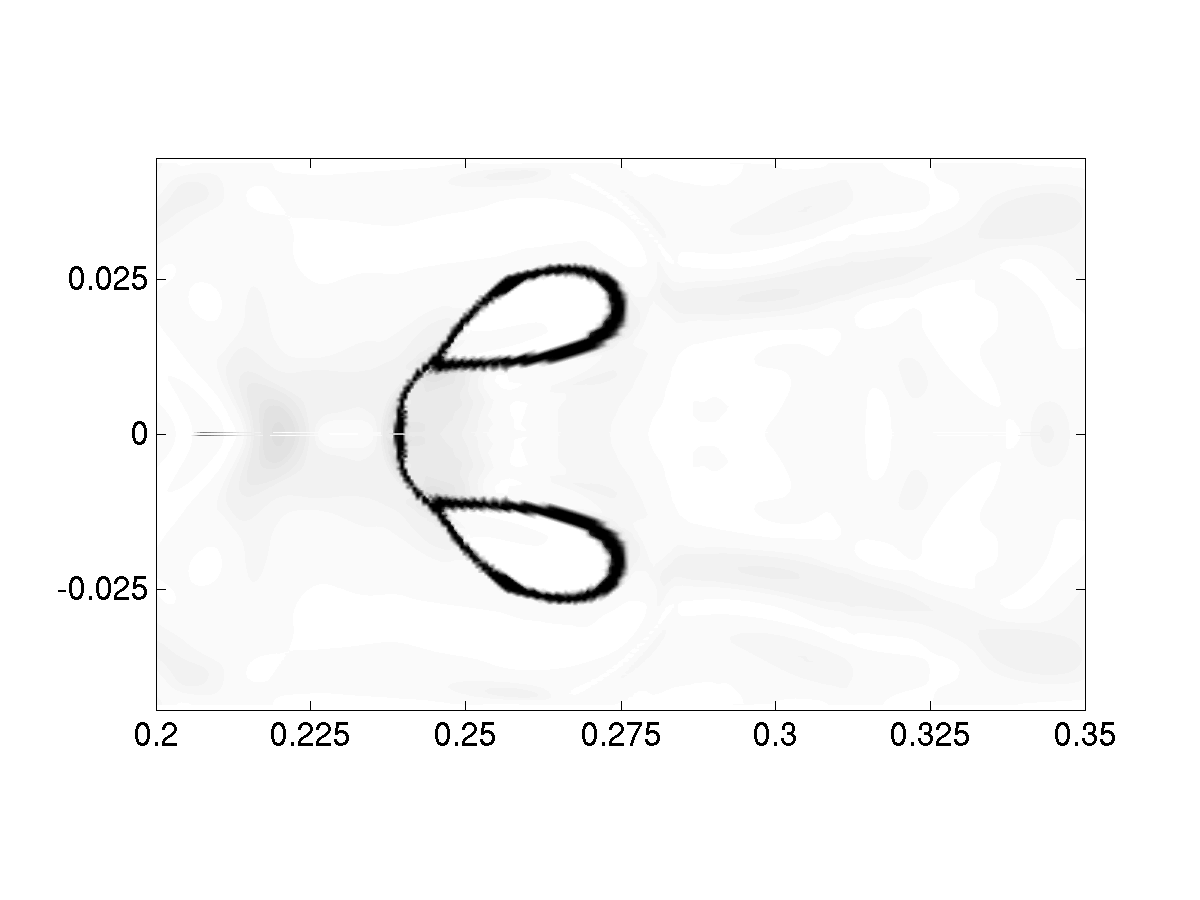}
\includegraphics[scale=0.35]{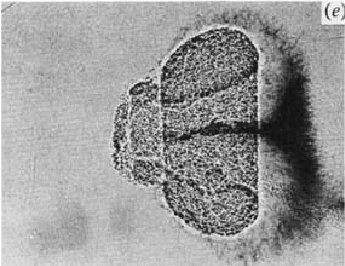} \\[-0.25cm]
$t=t_i+350\times10^{-6}$ \\[-0.15cm]
\includegraphics[scale=0.25, clip,trim=0cm 2.2cm 0cm 1cm ]{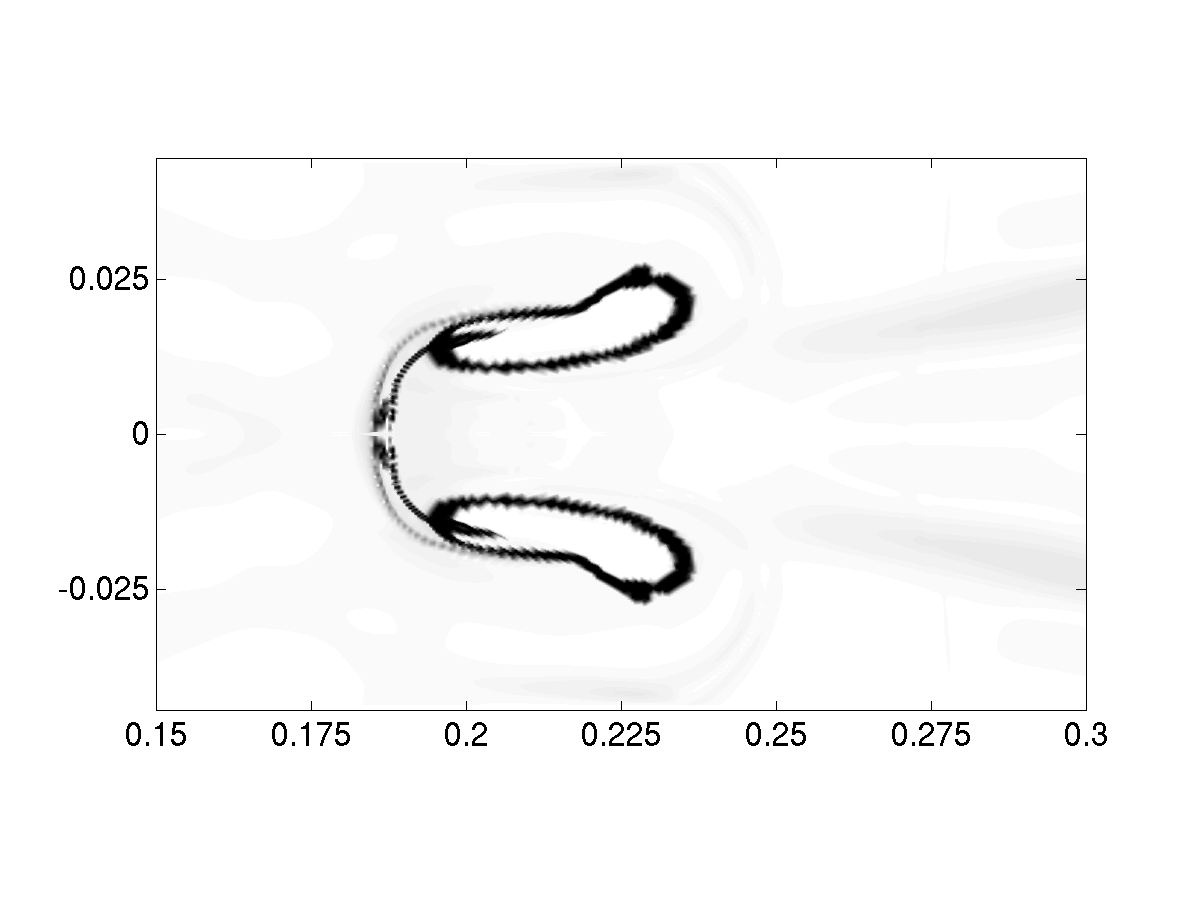}
\includegraphics[scale=0.35]{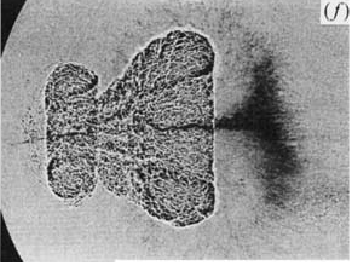} \\[-0.25cm]
$t=t_i+600\times10^{-6}$
\end{tabular}
\caption{Spherical Air-Helium shock/bubble interaction. Schlieren diagram of density. Axi-symmetric CCALE-MOF results (on the left) compared to experimental results (on the right) \cite{Haas1}. \label{bulle:1}}
\end{figure}

\begin{figure}[H]
\centering
\small
\begin{tabular}{cc}
\includegraphics[scale=0.32,clip,trim=0.cm 7cm 0cm 7cm]{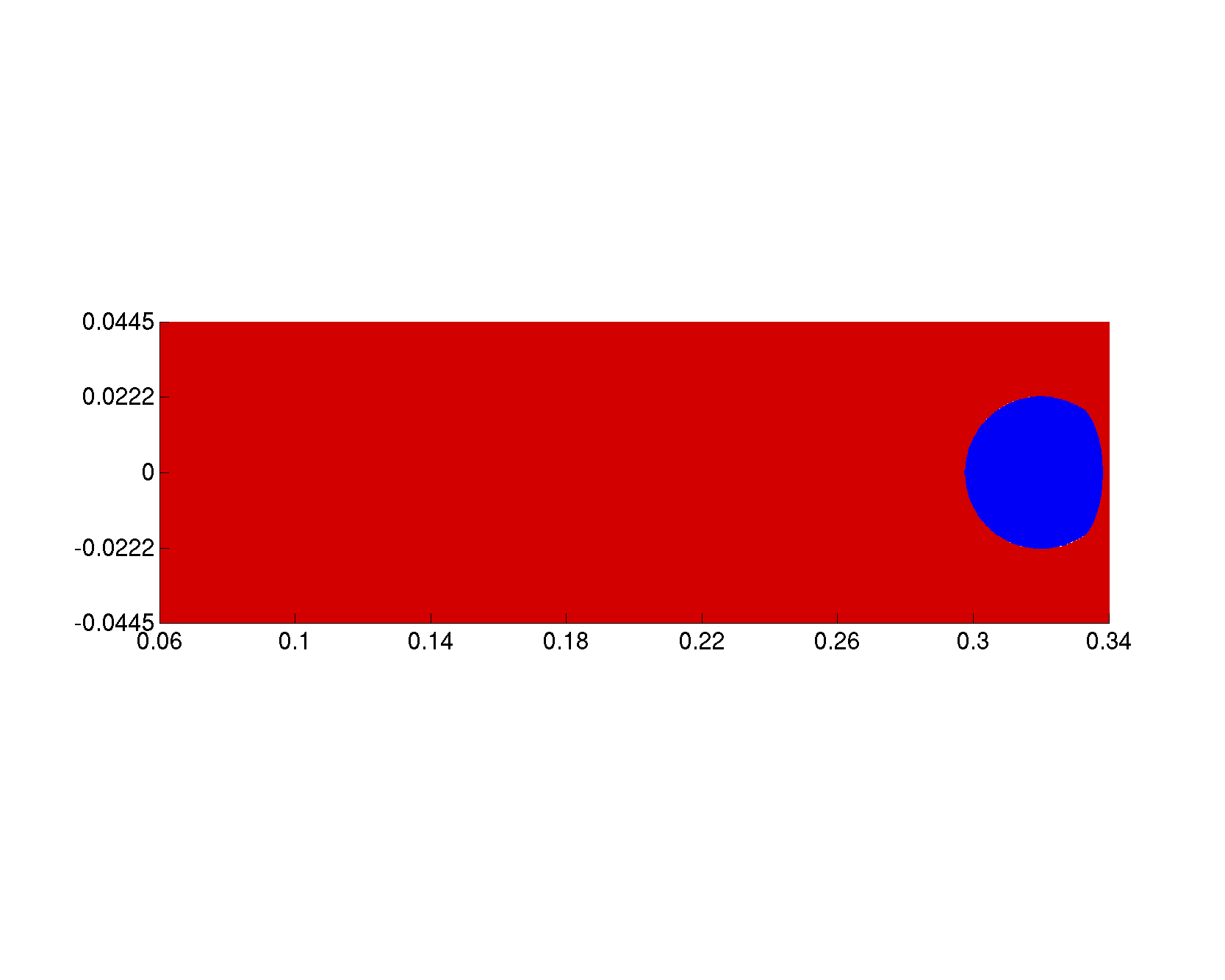} \\
$t=t_i+20\times10^{-6}$\\
\includegraphics[scale=0.32,clip,trim=0.cm 7cm 0cm 7cm]{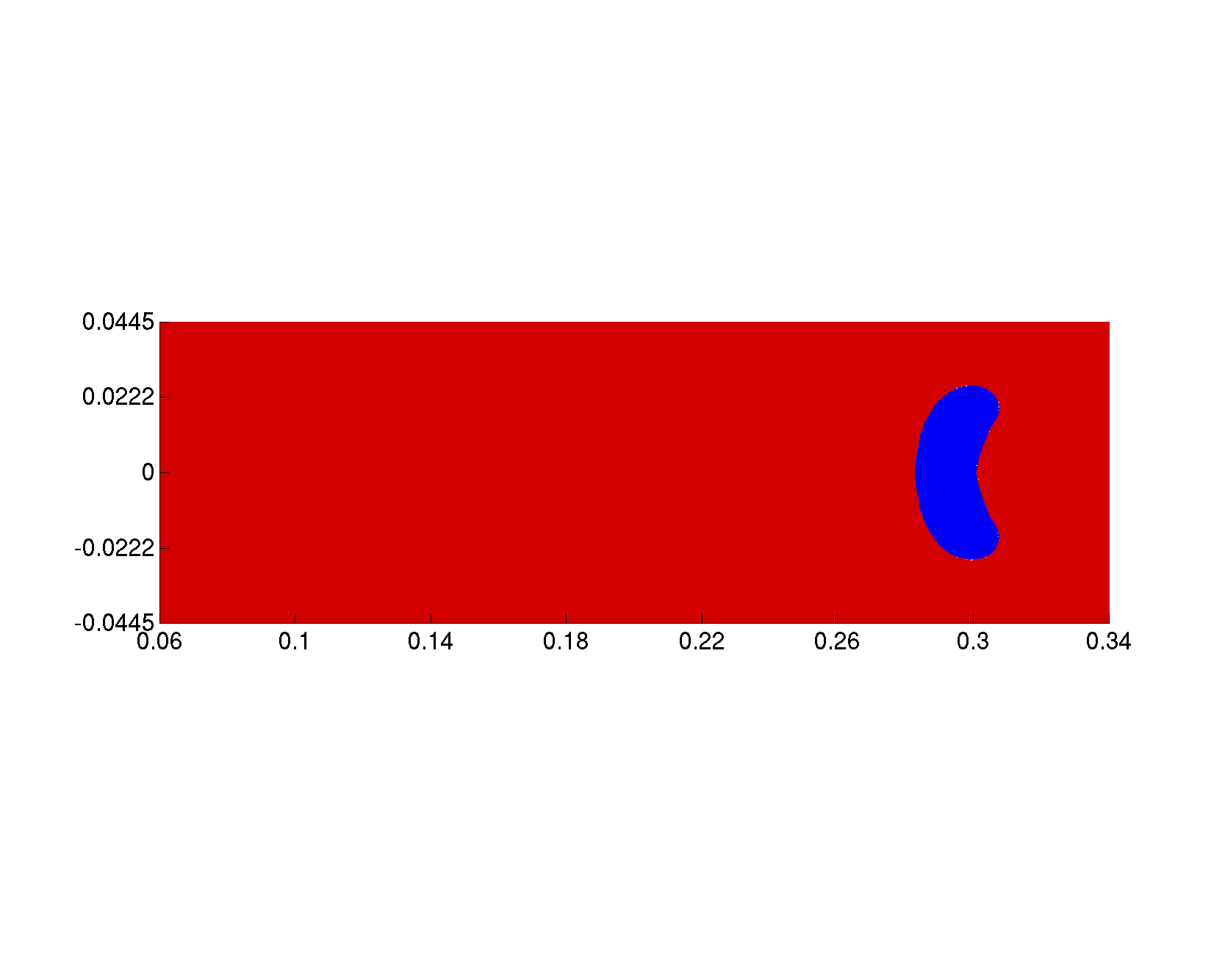} \\ 
$t=t_i+145\times10^{-6}$\\
\includegraphics[scale=0.32,clip,trim=0.cm 7cm 0cm 7cm]{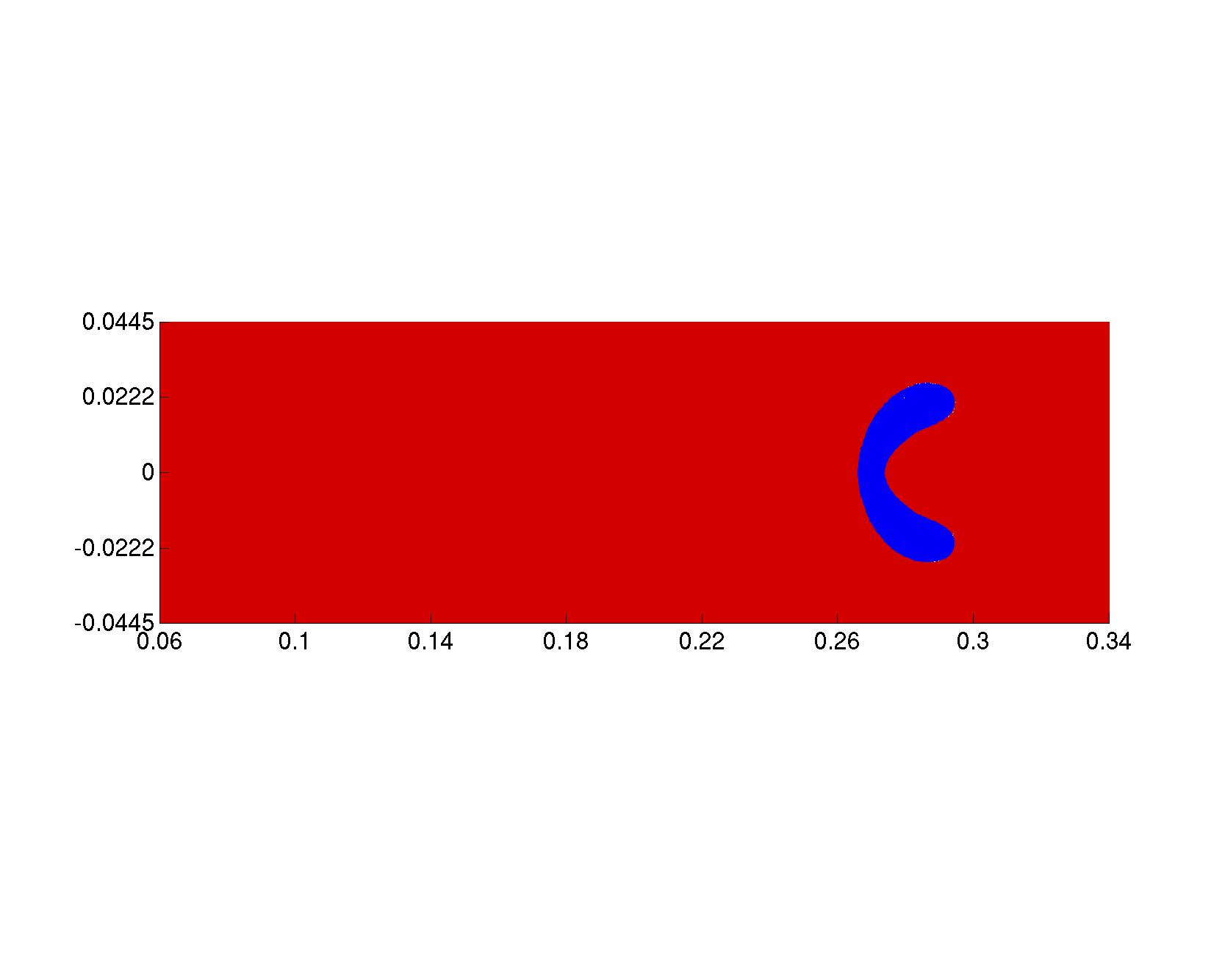}\\
$t=t_i+223\times10^{-6}$\\
\includegraphics[scale=0.32,clip,trim=0.cm 7cm 0cm 7cm]{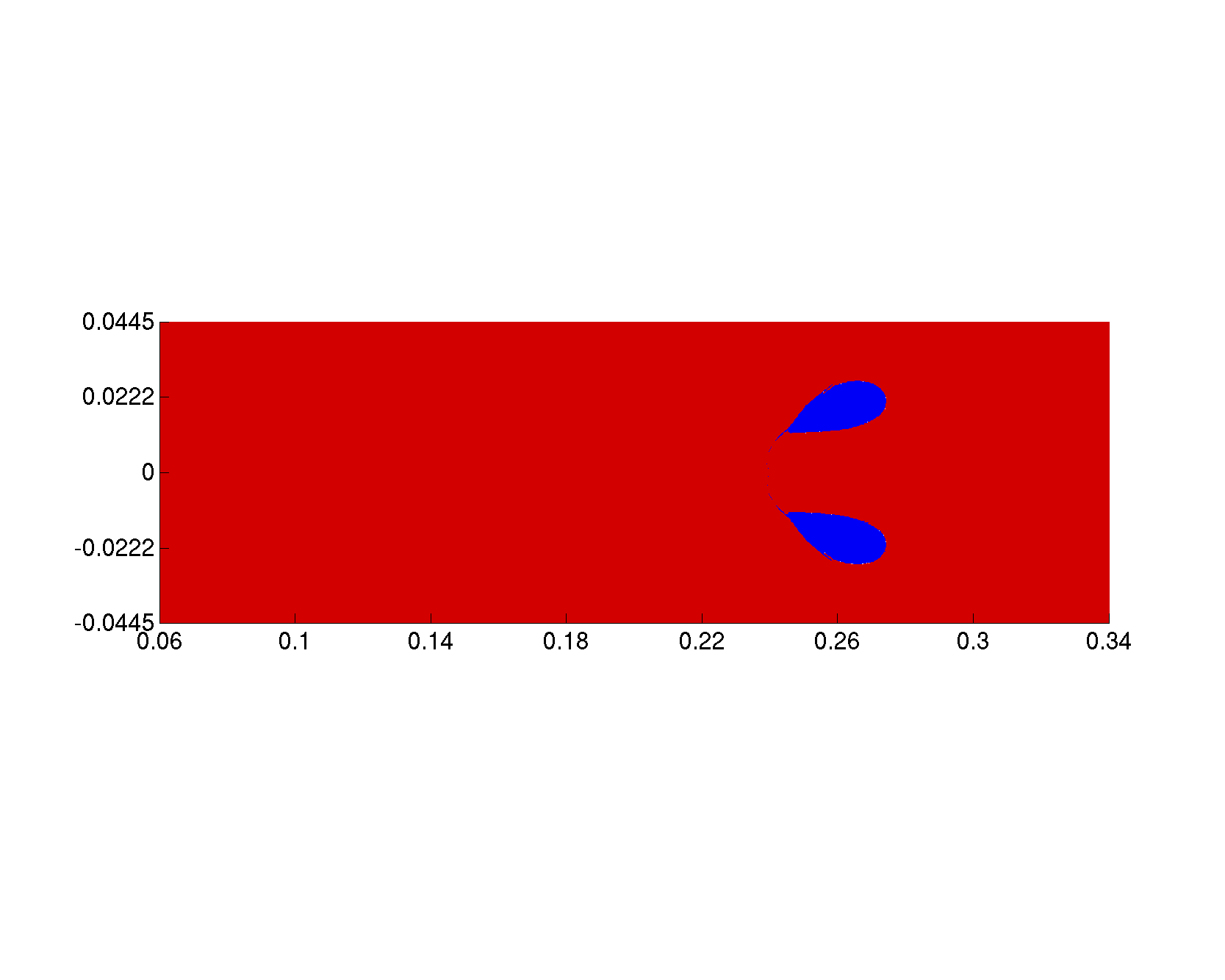} \\
$t=t_i+350\times10^{-6}$ \\
\includegraphics[scale=0.32,clip,trim=0.cm 7cm 0cm 7cm]{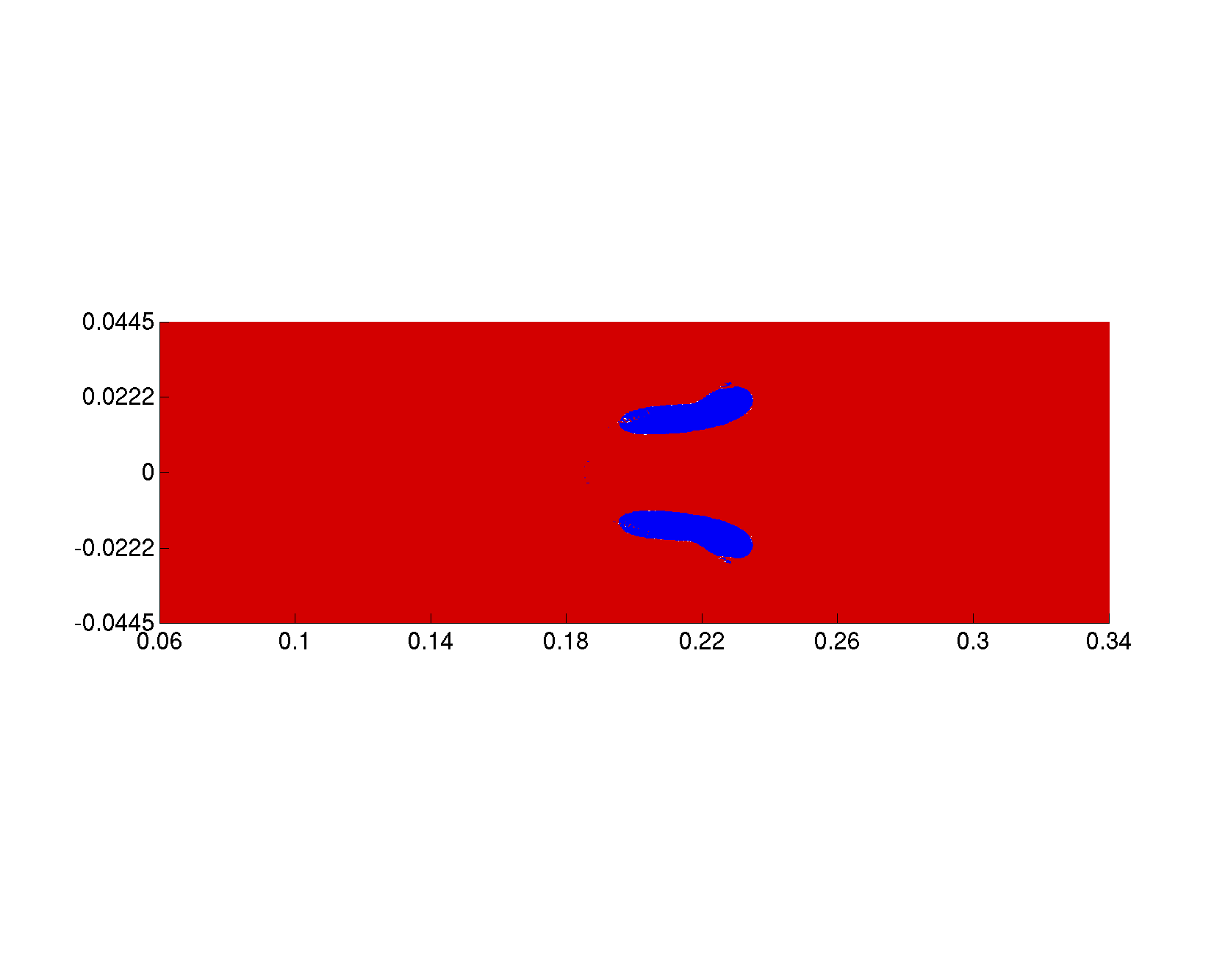} \\
$t=t_i+600\times10^{-6}$
\end{tabular}
\caption{Spherical Air-Helium shock/bubble interaction. Mesh and material interface evolution after the shock hits the bubble at time $t_i=657.463\times10^{-6}$. \label{bulle:2}}
\end{figure}

We deal in this part with the numerical simulation of the experiment of \cite{Haas1} concerning the impact of a Mach $1.25$ shock travelling through the air onto a spherical bubble of Helium. To this goal, let us consider a rectangular domain of dimensions $[0,0.65]\times[0,0.0445]$ initially full of Air of data $(\rho_1,P_1)=(0.182,10^5)$ except in an half disc centered in $(0,032)$ of radius $0.0225$ that contains Helium characterized by $(\rho_2,P_2)=(1,10^5) $  as depicted on \figref \ref{fig:10}.  Here, spherical geometry is obtain thanks to a rotation around the $X$-axis. EOS parameters for each fluids are stated on \tabref \ref{tab:1}. Wall boundary and  symmetry conditions are respectively chosen for the left, top boundaries. Despite, a piston-like condition is imposed to the right one for an incoming velocity equal to ${\VU}^* =(u^*,0)$. Here, the horizontal velocity $u^*$ is computed thanks to Rankine-Hugoniot conditions and is given by $u^*=-140.312$ corresponding to an incident shock moving at the velocity $  
D_c=-467.707$.\\
\begin{table}
\centering
\begin{tabular}{l|cc}
\hline
\hline
Fluid   & Polytropic index $\gamma$ & Molar mass $\mathscr{M}$ \\
\hline
Air     & $1.4$                       & $28.963$             \\
Helium  & $1.648$                     & $5.269\times10^{-3}$ \\
\hline
\end{tabular}
\caption{Air-Helium shock/bubble interaction: EOS parameters. \label{tab:1}}
\end{table}
The domain is initially paved with a structured cartesian grid composed of $520\times72$ cells.
Here, the bubble is directly initialized through the volume fraction on this mesh. Computations are done for the multi-material axisymmetric CCALE-MOF for a final time chosen equal to $t_f=t_i+600\times10^{-6}$ where $t_i=657.463\times10^{-6}$ corresponds to the time of the shock/bubble interaction. Here once again, simulations can not be achieved using pure Lagrangian framework due to the apparition of important mesh distortion. Numerical results associated to the Schlieren density profiles \cite{Hadjadj1} and interface positions deduced from the MOF reconstruction are respectively depicted on \figref \ref{bulle:1} and   \figref \ref{bulle:2}. Let us note that each pictures are obtained thanks to an axial symmetry  with respect to the $X$-axis. Comparisons between the Schlieren density profiles and the sadow-graphs of the experiment show a good agreement, especially when observing the bubble shape deformations. Moreover, waves generated by the initial shock  are well localized and illustrate multiple reflections and refractions especially on the bubble and the domain boundaries. These main points clearly demonstrate the accuracy and the robustness of the method and validate the axisymmetric CCALE-MOF approach when computing spherical test-cases coming from experiment. 

\subsection{Spherical implosion}

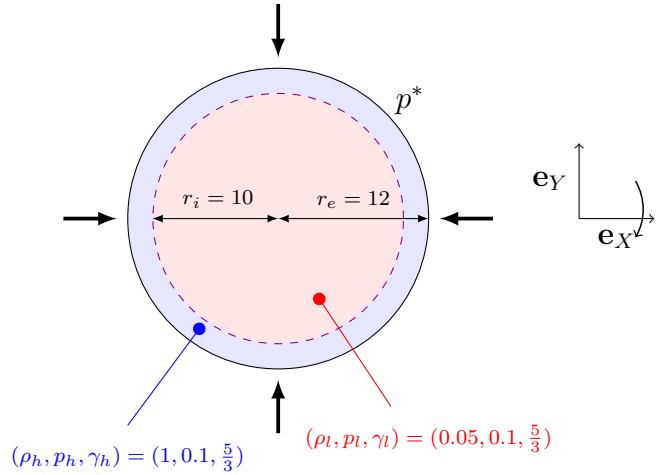
\begin{figure}[h!]
\centering
\begin{tikzpicture}[xscale=0.5,yscale=0.5]
\filldraw[fill=blue!10,draw=black] (5,5) circle (4);
\filldraw[fill=red!10,draw=violet,dashed] (5,5) circle (3.33);
\tikzstyle{fleche}=[->,>=latex,line width=1mm]
	\tikzstyle{fleche1}=[->,>=latex,line width=0.5mm]
	\tikzstyle{fleche2}=[<->,>=latex]
	\node (D1) at (1,5){};
	\node (G1) at (-1,5){};
	\node (D2) at (11,5){};
	\node (G2) at (9,5){};
	\node (D3) at (5,1){};
	\node (G3) at (5,-1){};
	\node (D4) at (5,11){};
	\node (G4) at (5,9){};
	\draw[fleche1] (G1)--(D1);
	\draw[fleche1] (D2)--(G2);
	\draw[fleche1] (G3)--(D3);
	\draw[fleche1] (D4)--(G4);
	\draw[fleche2] (5,5)--(9,5);
	\draw[fleche2, above] (5,5)--(1.65,5); 
	\draw (3.4,5) node [above]{\scriptsize $r_i=10$};
	\draw (7,5) node [above]{\scriptsize $r_e=12$};
	\draw (9,-1.5) node [above]{\scriptsize \textcolor{red}{$(\rho_l,p_l,\gamma_l)=(0.05,0.1,\frac{5}{3}) $}};
	\draw[-*,draw=red,fill=red] (8,0)--(6,3); 
	\draw (1,-2) node [above]{\scriptsize \textcolor{blue}{$(\rho_h,p_h,\gamma_h)=(1,0.1,\frac{5}{3}) $}};
	\draw[-*,draw=blue,fill=blue] (1,-0.5)--(3,2.2);
	\draw(8.5,7.5) node [above] {$p^*$};
	\draw[<-,line width=0.2mm] (14.5,4.5) arc (-30:30:1.5);
	\node[left] at (13,6) {${\bf e}_Y$};
	\node[below] at (14,5) {${\bf e}_X$};
	\draw[->] (13,5)--(13,7); 
	\draw[->] (13,5)--(15,5);
\end{tikzpicture}  
\caption{Multi-mode implosion in spherical geometry. Initial geometry and data. \label{fig:4}}
\end{figure}

The last test-case of this paper deals with the numerical computation of a spherical  implosion as initially treated in \cite{Youngs1}. The interest of this simulation is twofold. First, this is a realistic problem quite close to those encountered in Ignition Confinement Fusion (ICF)  simulation. Then, it allows to test the capability of the multi-material CCALE-MOF algorithm with hybrid rezoning.

\begin{figure}[h!]
\centering
\begin{tabular}{cc}
\includegraphics[scale=0.3]{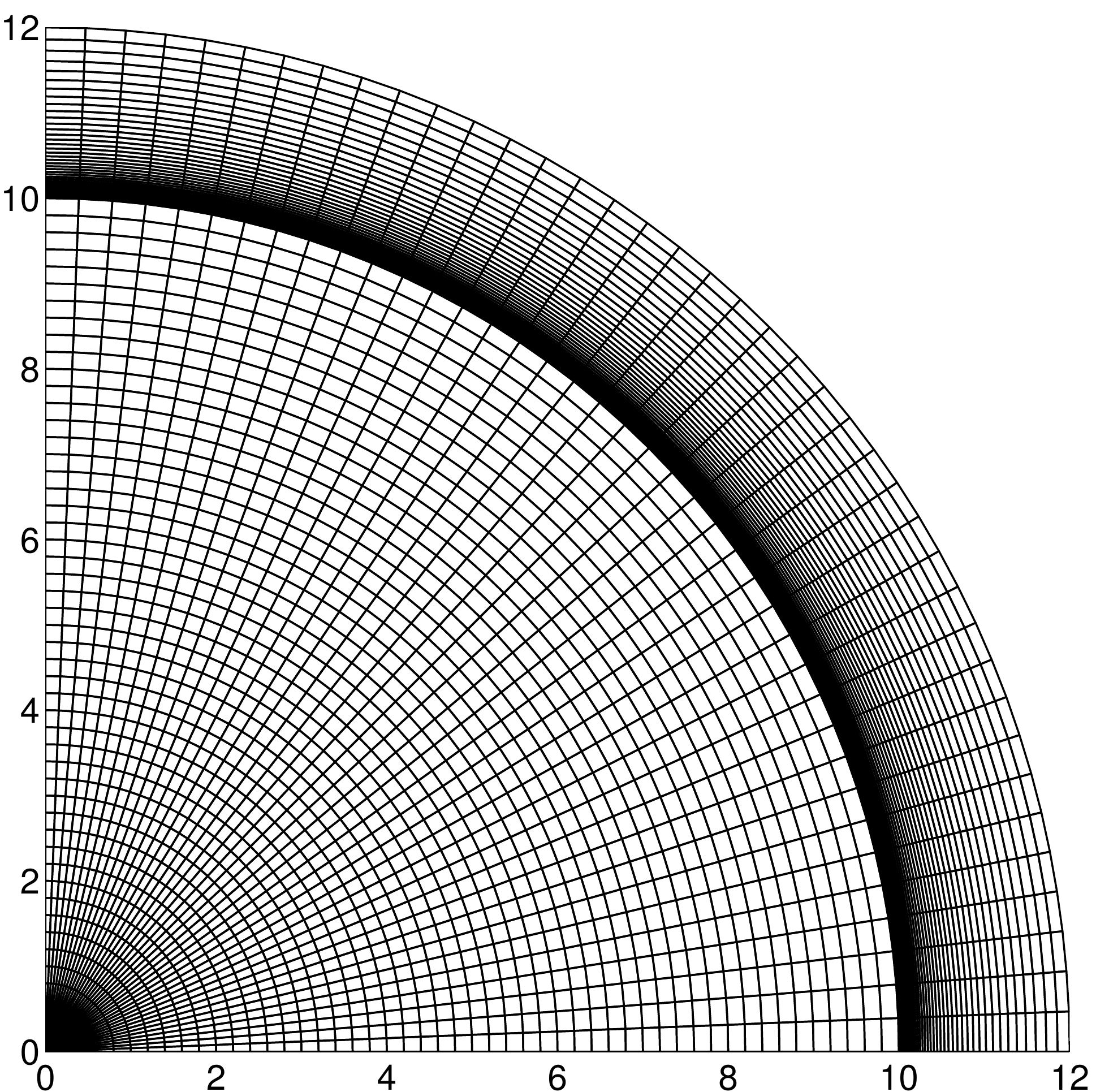} & 
\includegraphics[scale=0.3]{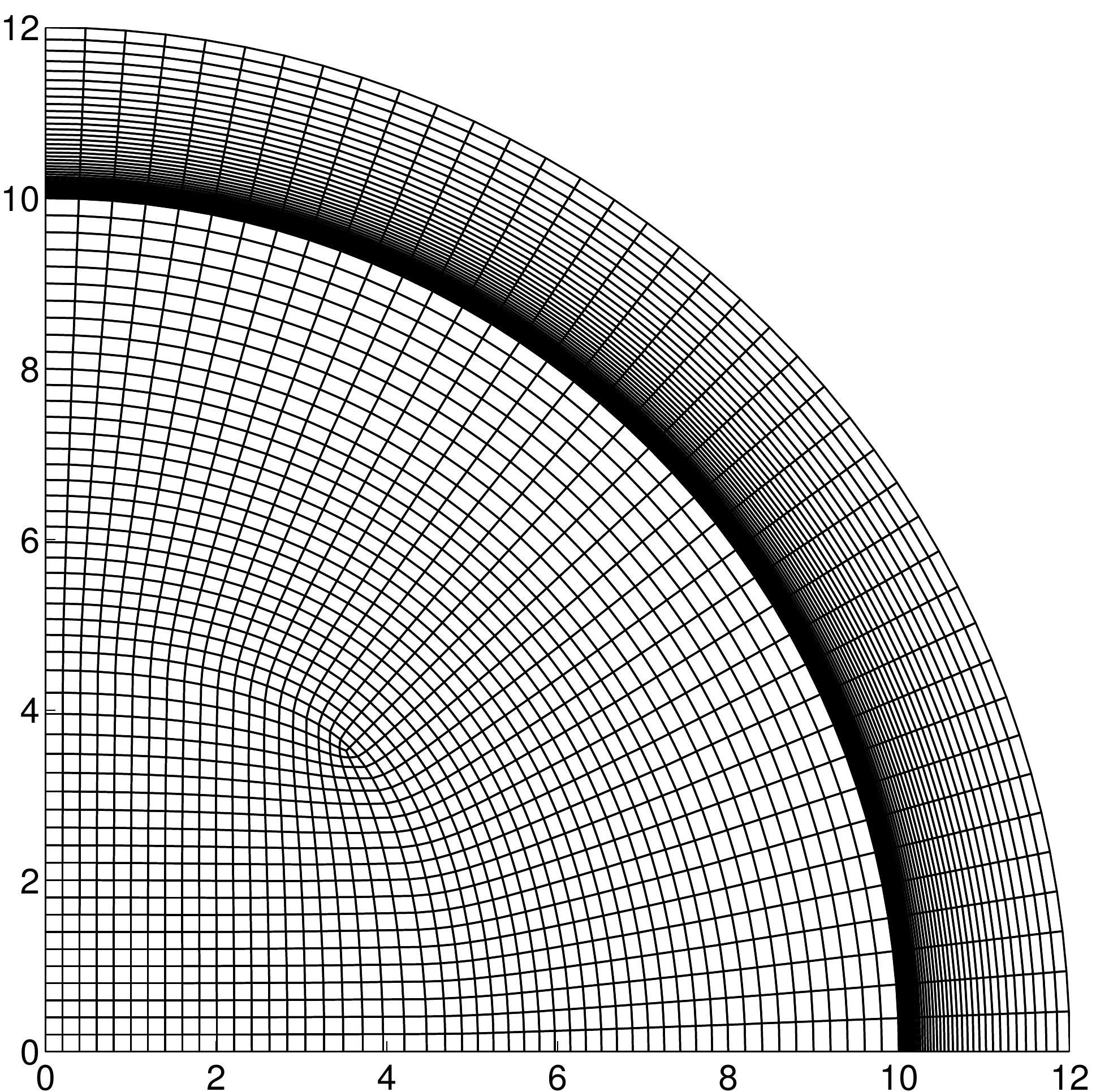}
\end{tabular}
\caption{Spherical implosion. Initial polar (left) and unstructured (right) grids. \label{youngs:1}}
\end{figure}

Here we focus on the treatment of perturbed interfaces where compressible Rayleigh-Taylor instabilities occur.\\
Let us consider a spherical ball of light fluid ($r\in [0,10]$) initially surrounded by a shell of heavy fluid ($R\in[10,12]$) as depicted on \figref \ref{fig:4}. For both fluid the polytropic index is the same $ \gamma_l=\gamma_h=\dfrac{5}{3}$. The initial pressures and densities are $(\rho_l,p_l) = (0.05,0.1)$ and $(\rho_h,p_h) =(1,0.1)$. The implosion is driven by imposing the following pressure law on the dense shell boundary
$$
p^*(t) = 
\left\{
\begin{array}{rl}
10      &  \text{ if } t\in [0,0.5], \\
12 - 4t &  \text{ if } t\in [0.5, 3].
\end{array}
\right.
$$
Finally, the interface between the light and the heavy fluids is initially perturbed according to the law
$$
r_p^{per} = r_p (1+a_0 {\cal D}(r_p) P_l(\cos(\theta_p))
$$
with the damping factor
$$
{\cal D}(r_p)=
\left\{
\begin{array}{rl}
1-\dfrac{r_p-r_i}{r_e-r_i}     &  \text{ if } r_p\in [r_i,r_e], \\
1-\dfrac{r_i-r_p}{r_i}         &  \text{ if } r_p \in[0,r_i].
\end{array}
\right.
$$
where $r_i^{per}$ denotes the perturbed radius and $a_0$ is the amplitude of the perturbation. Finally, $P_l$ is the $l^{th}$ Legendre polynomial. In the sequel $l=10$ and several values of $a_0$ are considered from the non-perturbed case $a_0=0$, to weakly and strongly perturbed one with respectively $a_0=2\times10^{-4}$ and $a_0=1\times10^{-3}$.\\
Computations are made for two different meshes until the final time $t_f=3$. The first one is a polar grid displayed on \figref \ref{youngs:1}-(left) composed of $90\times40$ cells. Size of cells in the radial direction have been chosen respecting a mass radial spacing deduced from the equivalent one-dimensional test case. The other grid, is obtained after an hybrid regularization for $\omega_p=1$ of an unstructured mesh initially paved with $3200$ quadrangular cells respecting the mass radial spacing (see \figref \ref{youngs:1}-(right)).

\begin{figure}[h!]
\centering
\begin{tabular}{cc}
\includegraphics[scale=0.45]{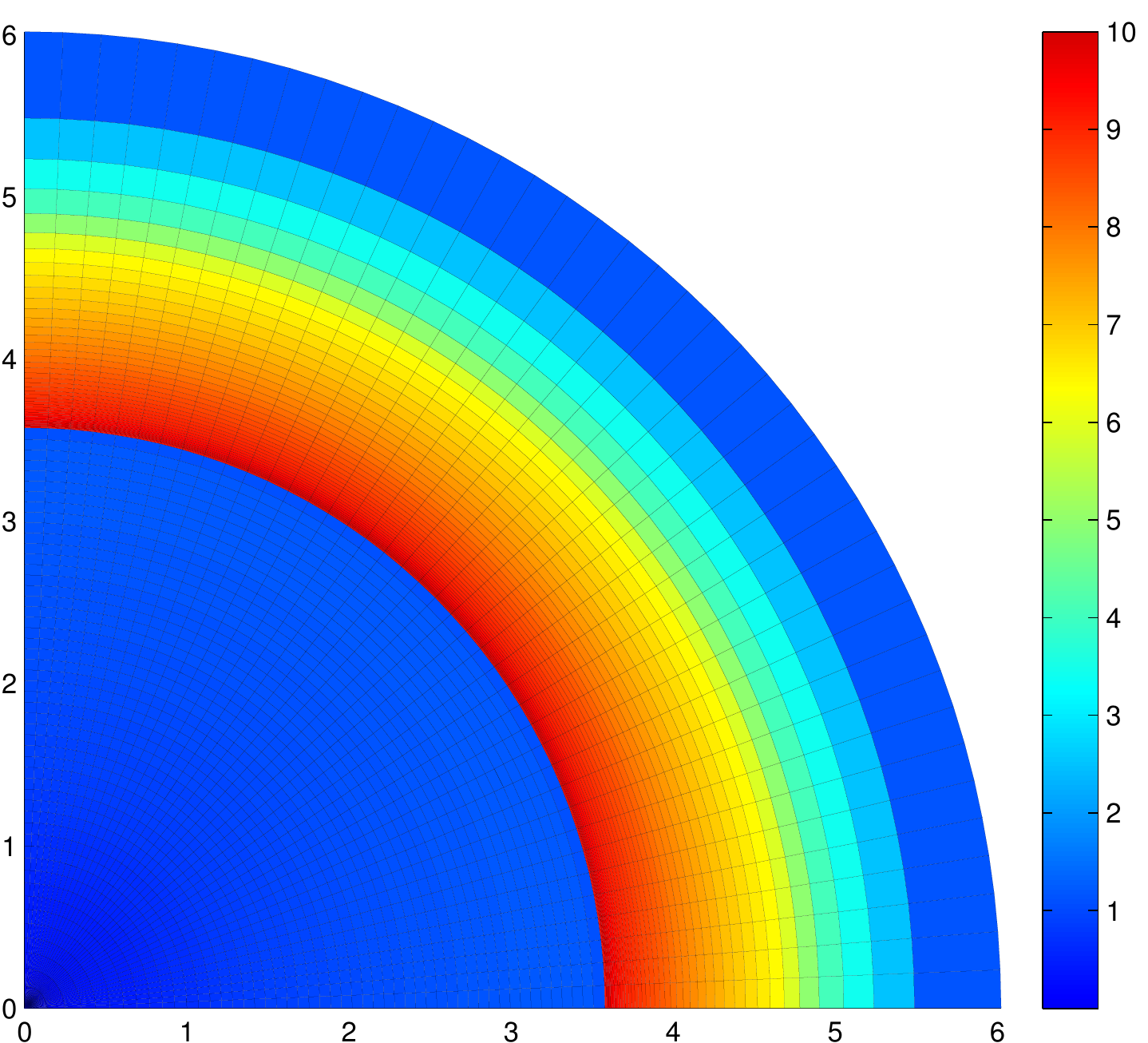} & 
\includegraphics[scale=0.45]{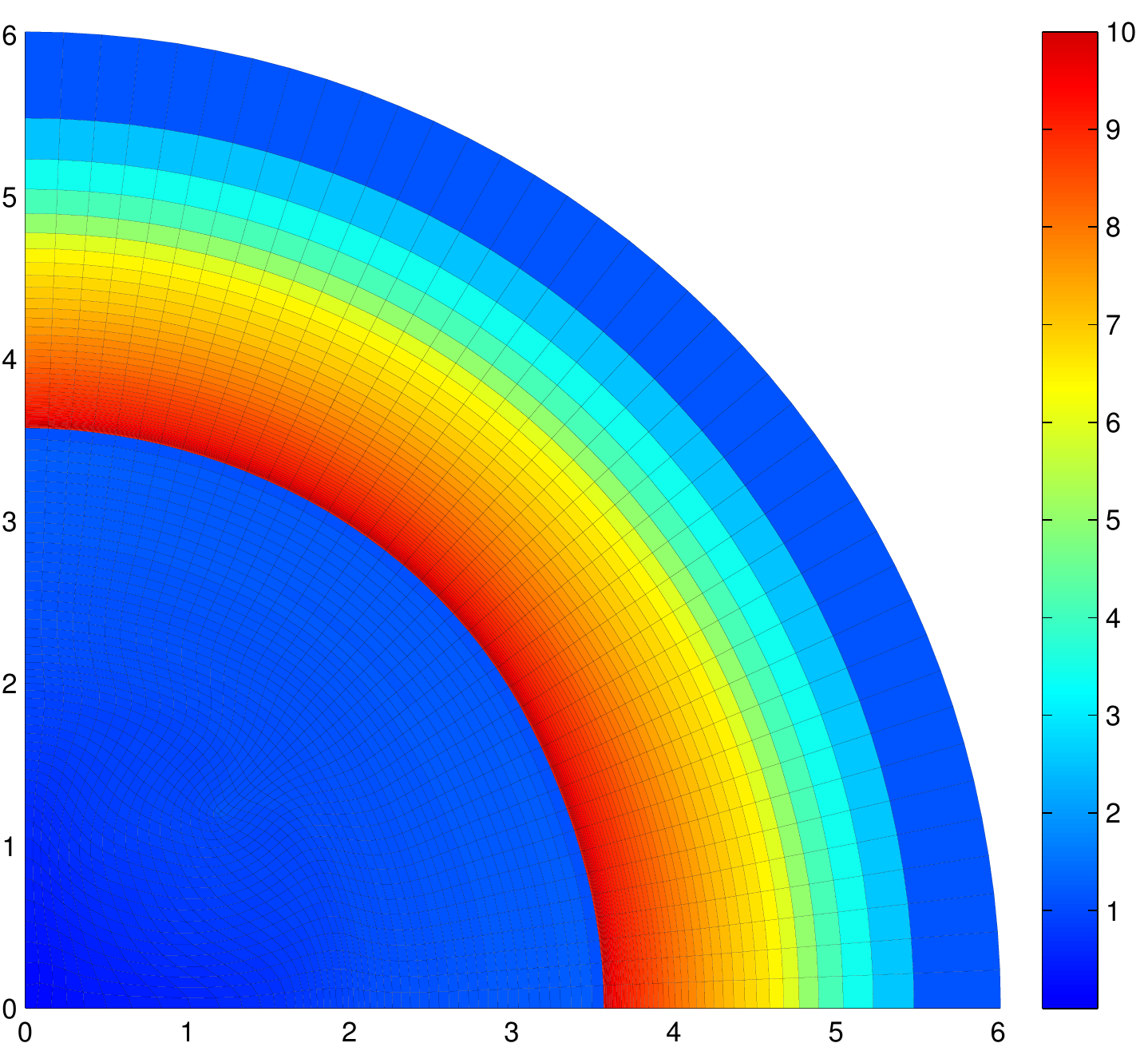}
\end{tabular}
\caption{Spherical implosion without deformation. Mesh and density for polar (left) and unstructured (right) grids at final time $t_f=3$. \label{youngs:2}}
\end{figure}

\paragraph{Non-perturbed case with $a_0=0$}
As a first study, we test the behavior of our algorithm in axisymmetric geometries in pure Lagrange computation for both meshes. As shown on \figref \ref{youngs:2}, numerical results for both meshes are similar. Nevertheless, the method remains faster on the unstructured mesh. Indeed, it has the advantage to not impose a drastic time step for computation due to triangular cells with high aspect ratio in the polar mesh as shown in \cite{Galera2}. 

\paragraph{Weakly perturbed case with  $a_0=2\times10^{-4}$}

\begin{figure}[h!]
\centering
\begin{tabular}{cc}
\includegraphics[scale=0.45]{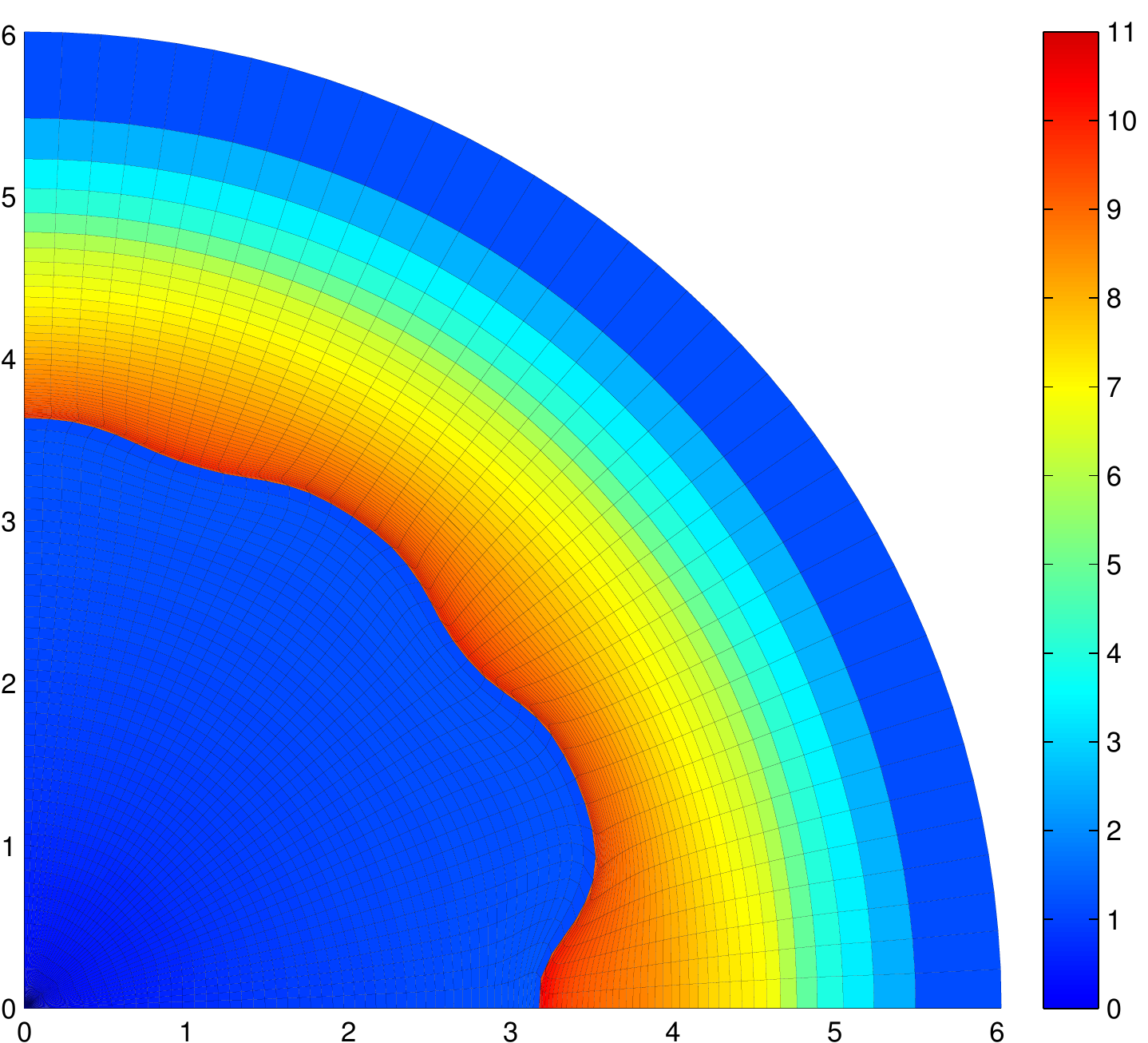} & 
\includegraphics[scale=0.45]{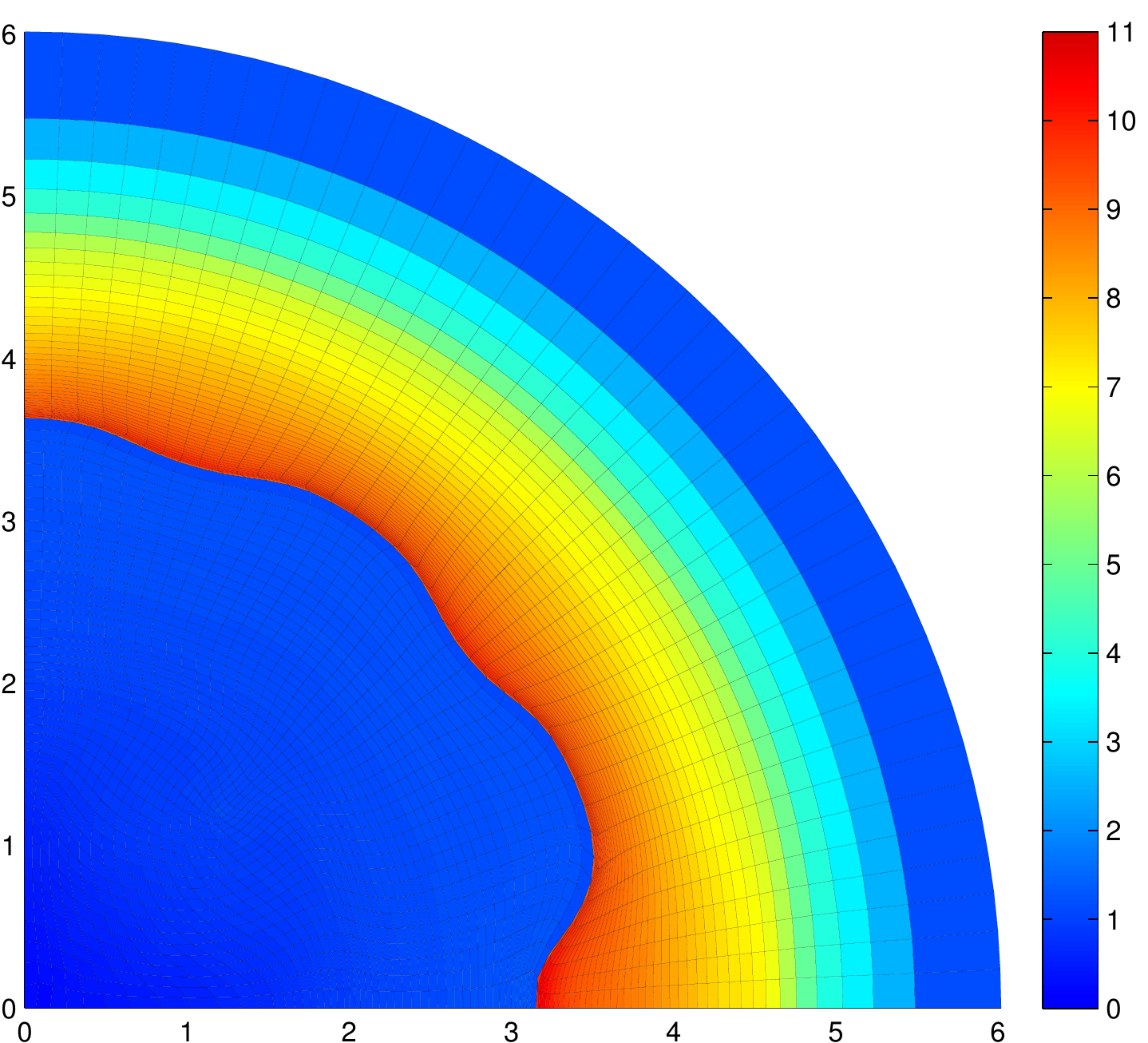}\\
\includegraphics[scale=0.32]{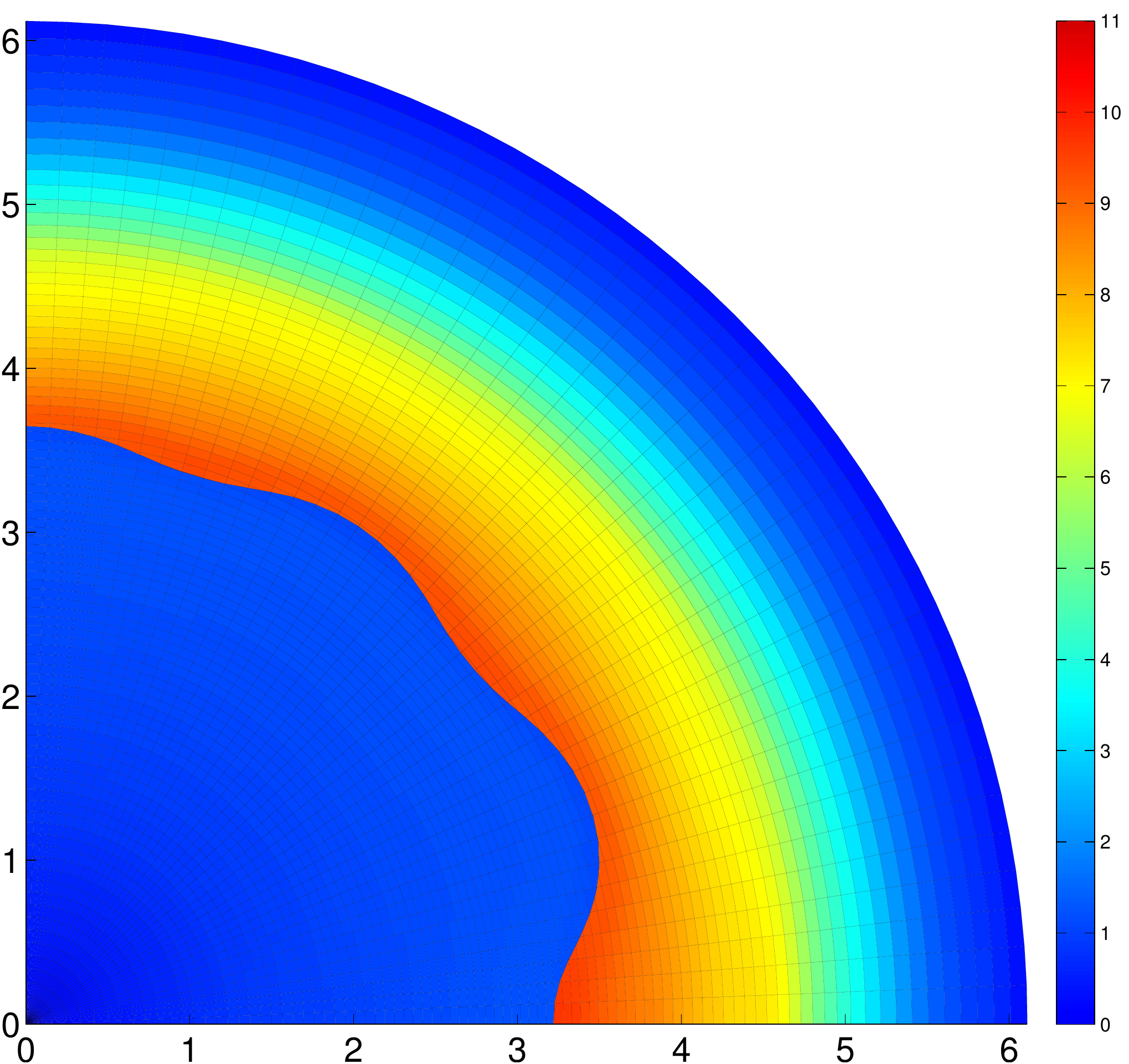} & 
\includegraphics[scale=0.32]{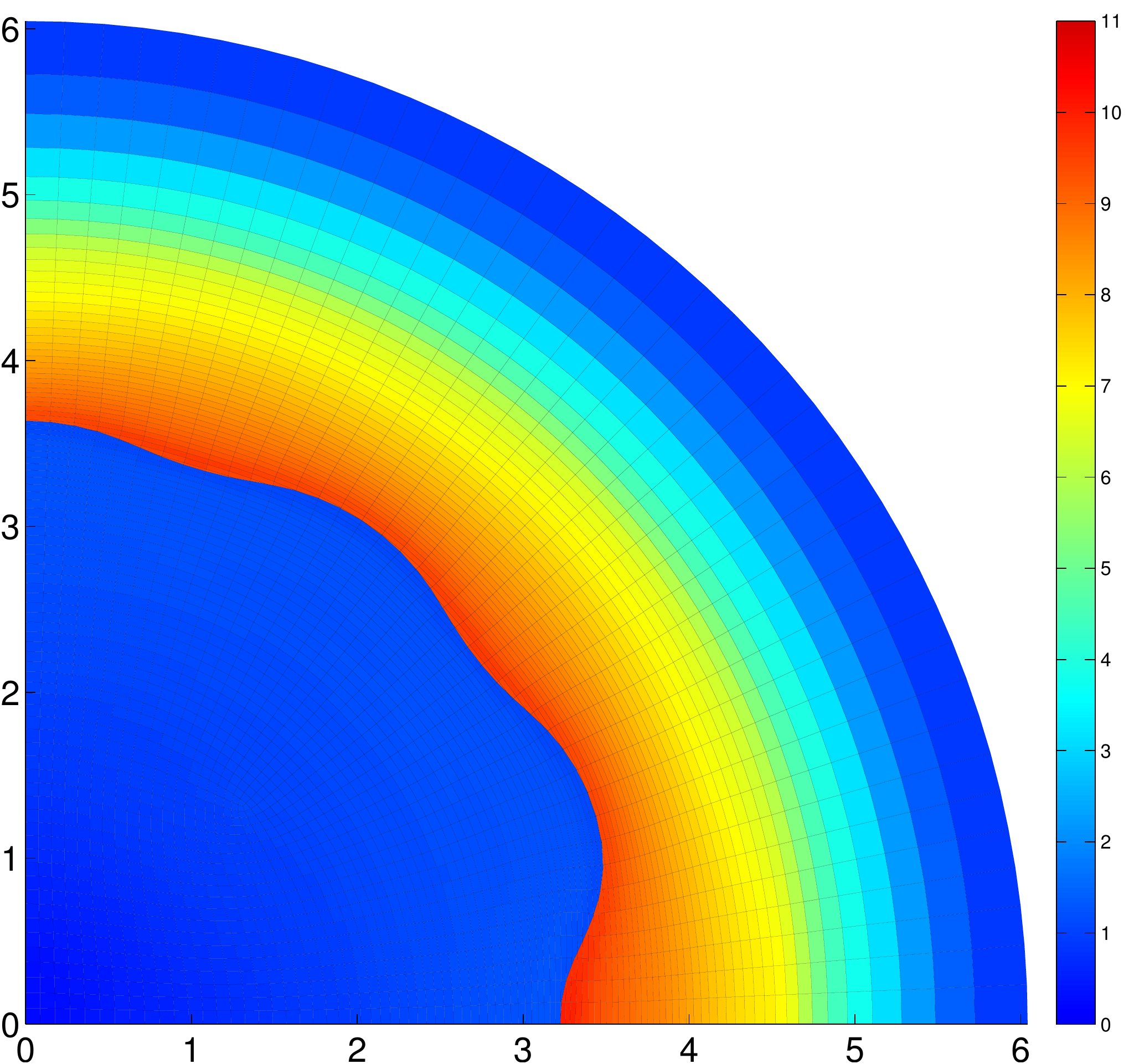}
\end{tabular}
\caption{Spherical implosion with small deformation. Mesh and density for Lagrangian (top) and ALE (bottom) computations at final time $t_f=3$ for both polar (left) and unstructured (right) grids. \label{youngs:3}}
\end{figure}

Now, we investigate the capability of our CCALE-MOF algorithm to treat perturbed interfaces on both non-structured and polar meshes. To this end, comparisons with pure Lagrangian results are first achieved for weakly perturbed interfaces imposing $a_0=2\times10^{-4}$. Here for both polar and hybrid meshes, the GCNS is used. As demonstrated on \figref \ref{youngs:3}, for the polar mesh as well as for the non-structured mesh, ALE results, especially concerning the interface deformation, are in very good agreement to thoses obtained thanks to pure Lagrangian computations. Furthemore, one should note that for the ALE  computation on polar grid the quality of the mesh is improved near the origin. Indeed, the central cells are  not systematically shifted to the origin contrary to computations achieved using CNS rezoning. 

\paragraph{Strongly perturbed case with  $a_0=1\times10^{-3}$}

\begin{figure}[h!]
\centering
\begin{tabular}{cc}
\includegraphics[scale=0.32]{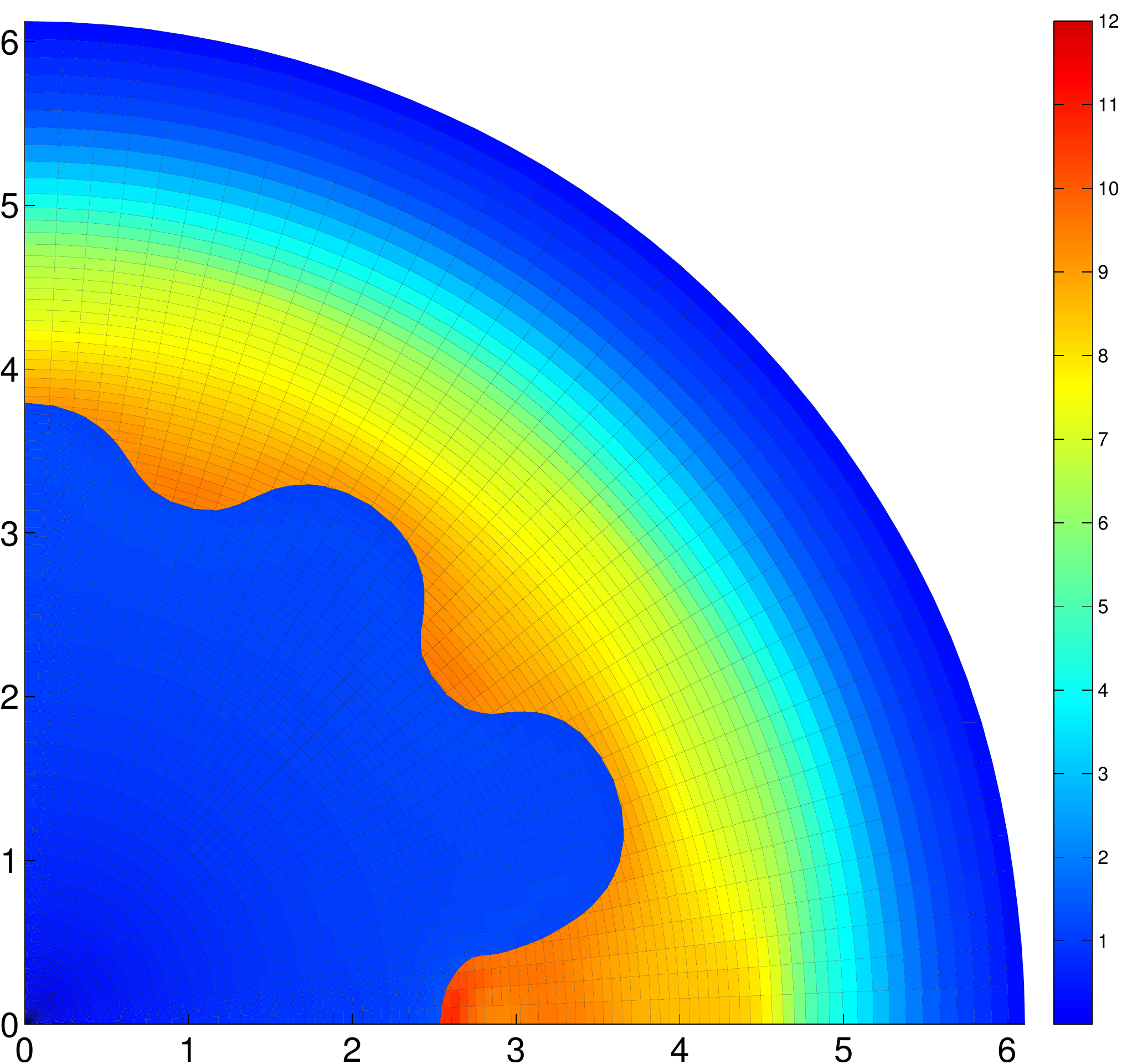} & 
\includegraphics[scale=0.32]{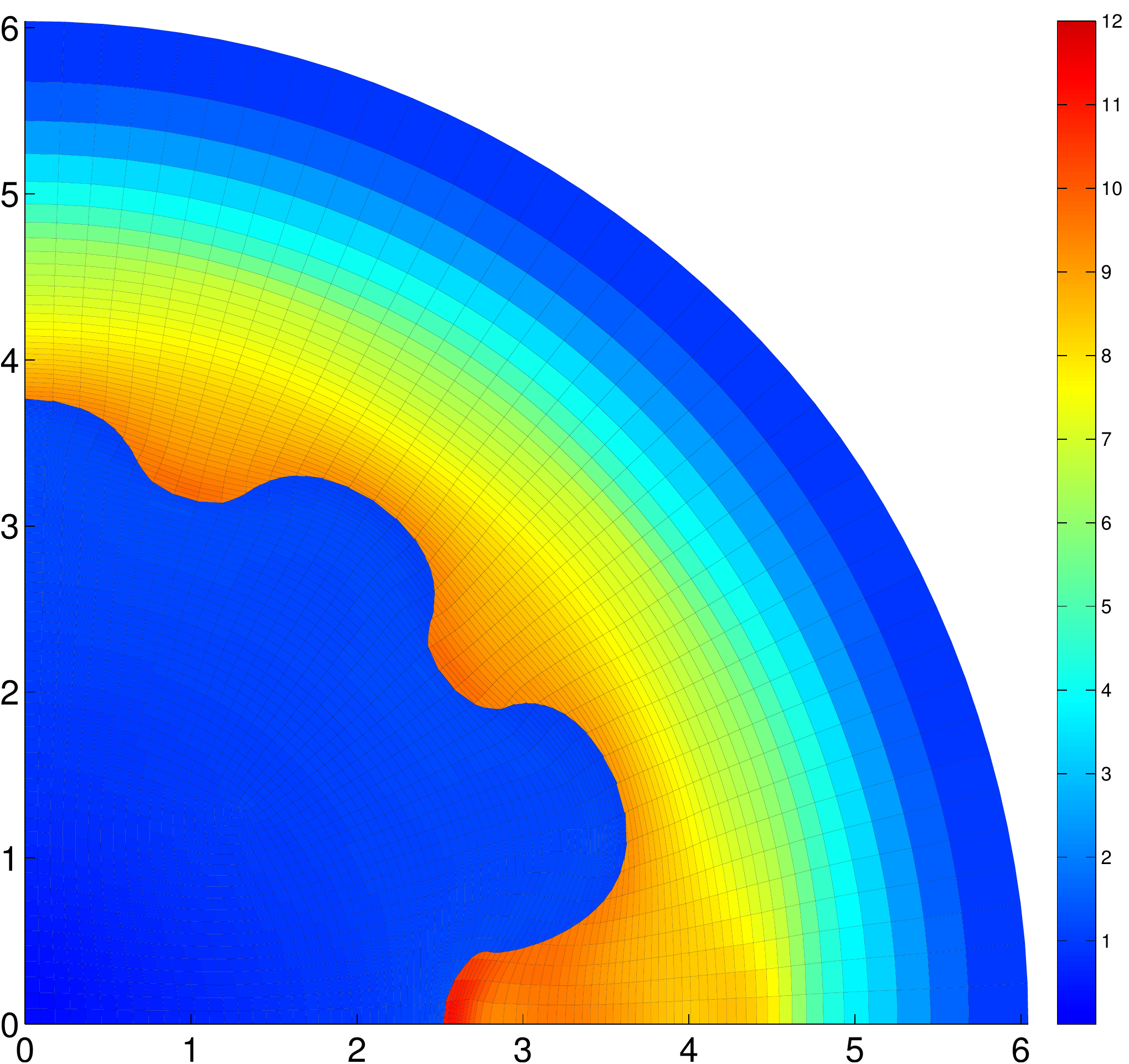}
\end{tabular}
\caption{Spherical implosion with important deformation. Mesh and density for ALE computation for both polar (left) and unstructured (right) grids at final time $t_f=3$. \label{youngs:4}}
\end{figure}

Finally, we perform a computation of this implosion for a more pertubated interface choosing $a_0$ five times greater than previously with  $a_0=1\times10^{-3}$. Due to mesh tangling, this is not possible to purchase such a test case using only Lagrangian method whose computation fails for $t>t_{fail}=2.6$.  Here, only results obtained thanks to our axisymmetric multi-material CCALE-MOF are presented. Contrary, to Lagrangian computations, the multi-material ALE simulations run without any difficulties thanks to specific rezoning. For both grids, final results (see \figref \ref{youngs:4}) are very close. In particular we note the Rayleigh-Taylor instability has grown in a same way leading to similar interface shape deformation at final time. 

\section{Conclusion and future work}

In this paper, we have presented several extensions concerning a Cell-Centered Arbitrary Lagrangian-Eulerian (CCALE) strategy using the Moment of Fluid (MOF) interface reconstruction devoted to the numerical simulation of multi-material compressible flows especially in axisymmetric geometry on both polar and Cartesian unstructured meshes. To this end, we have introduced a simple and unified formulation of the Lagrangian scheme relying on an area-weighted formulation, a multi-material MOF interface reconstruction, a new formulation of rezoning for both polar and Cartesian grids and finally a general hybrid remap procedure for both axisymmetric and Cartesian geometry. As demonstrated on several academical as well as ICF-like test cases, the proposed method remains accurate and robust.\\
As future work, we plan to incorporate the proposed method in the multi-physic code CHIC dedicated to the simulation of ICF experiment. The main goal is to treat eventually more general configurations notably coupling realistic EOS, laser energy deposition, with multi-material hydrodynamics in the lines of \cite{Breil1}.


\newpage
\begin{thebibliography}{99.}
\bibitem{Ahn1} H.T. Ahn, M.J. Shashkov: {\it Multi-material interface reconstruction on generalized polyhedral meshes}, J. Comput. Phys., 226(2):2096-2132, 2007.  
\bibitem{Anbarlooei1} H.R. Anbarlooei, K. Mazaheri: {\it 'Moment of fluid' interface reconstruction method in axisymmetric coordinates}, Int. J. Numer. Meth. Biomed. Engng., 27(10):1640-1651, 2011.
\bibitem{Barlow1} A.J. Barlow, P.L. Roe: {\it  A cell centred Lagrangian Godunov scheme for shock hydrodynamics}, Comput. Fluids, 46(1):133-136, 2011.
\bibitem{Berndt1} M. Berndt, J. Breil, S. Galera, M. Kucharik, P.-H. Maire, M. Shashkov: {\it Two-step hybrid conservative remapping for multi-material arbitrary Lagrangian-Eulerian methods}, J. Comput. Phys., 230(17):6664-6687, 2011.
\bibitem{Breil1} J. Breil, S. Galera, P.-H. Maire: {\it Multi-material ALE computation in Inertial Confinement Fusion CHIC}, Comput. Fluids, 46(1):161-167, 2011.
\bibitem{Breil2}  J. Breil, S. Galera, P.-H. Maire: {\it A two-dimensional VOF interface reconstruction in a multi-material cell-centered ALE scheme}, Int. J. Numer. Meth. Fluids, 65(11-12):1351-1364, 2011.
\bibitem{Despres1} G. Carré, S. Del Pino, B. Després: {\it A cell-centered Lagrangian hydrodynamics scheme on general unstructured meshes in arbitrary dimension}, J. Comput. Phys., 228(14):5160-5183, 2009.
\bibitem{Dyadechko1} V. Dyadechko, M. Shashkov: {\it Reconstruction of Multi-material Interfaces from Moment Data}, J. Comput. Phys., 227(11):5361-5384, 2008.
\bibitem{Dukowicz1} J. K. Dukowicz: {\it A general, non-iterative Riemann solver for Godunov's method}, J. Comput. Phys., 61(1):119-137, 1985.
\bibitem{Galera1} S. Galera, J. Breil and P.-H. Maire: {\it A 2D unstructured multi-material Cell-Centered Arbitrary Lagrangian-Eulerian (CCALE) scheme using MOF interface reconstruction}, Comput. Fluids, 46(1):237-244, 2011.
\bibitem{Galera2} S. Galera, P.-H. Maire and J. Breil: {\it A two-dimensional unstructured cell-centered multi-material ALE scheme using VOF interface reconstruction},  J. Comput. Phys., 229(16):5755-5787, 2010.
\bibitem{Haas1} J.-F. Haas, B. Sturtevant: {\it Interaction of weak shock wave with cylindrical and spherical gas inhomogeneities}, J. Fluid. Mech., 181:41-76, 1987.
\bibitem{Hadjadj1} A. Hadjadj,  A. Kudryavtsev  {\it Computation and flow visualization in high-speed aerodynamics}, Journal of Turbulence, 6(16), 2005.
\bibitem{Hirt1} C.W. Hirt, A. Amsden, and J.L. Cook: {\it An arbitrary Lagrangian-Eulerian computing method for all flow speeds}, J. Comput. Phys., 14:227-253, 1974.
\bibitem{Knupp1} P. Knupp: {\it Achieving finite element mesh quality via optimization of the Jacobian matrix norm and associated quantities.
Part I-- a framework for surface mesh optimization}, Int. J. Numer. Meth. Engng, 48:401-420, 2000.
\bibitem{Kucharik1} M. Kucharik, J. Breil, S. Galera, P.-H. Maire, M. Berndt, M. Shashkov: {\it Hybrid remap for multi-material ALE}, Comput.  Fluids, 46(1):293-297, 2011.
\bibitem{Kucharik2} M. Kucharik, R.V. Garimella, S.P. Schofield, M.J. Shashkov: {\it A comparative study of interface reconstruction methods for multi-material ALE simulations}, J. Comput. Phys., 229(7):2432:2452, 2010.
\bibitem{Loubere1} R. Loub\`ere, P.-H. Maire, M. Shashkov, J. Breil, S. Galera: {\it ReALE: A reconnection-based arbitrary-Lagrangian-Eulerian method}, J. Comput. Phys., 229(12):4724-4761, 2010.
\bibitem{Maire1} P.-H. Maire: {\it A high-order cell-centered Lagrangian scheme for compressible fluid flows in two-dimensional cylindrical geometry}, J. Comput. Phys., 228(18):6882-6915, 2009.
\bibitem{Maire2} P.-H. Maire: {\it A high-order cell-centered Lagrangian scheme for two-dimensional compressible fluid flows on unstructured meshes}, J. Comput. Phys., 228(7):2391-2425, 2009.
\bibitem{Maire3} P.-H. Maire, R. Abgrall, J. Breil, J. Ovadia: {\it A cell-centered Lagrangian scheme for two-dimensional compressible flow problems}, SIAM Journal of Scientific Computing, 29(4):1781-1824, 2007.
\bibitem{Maire4}
P.-H. Maire:{\it Contribution to the numerical modeling of Inertial Confinement Fusion},
Habilitation \`a Diriger des Recherches, Bordeaux University, 2011;
Available at: \url{http://tel.archives-ouvertes.fr/docs/00/58/97/58/PDF/hdr_main.pdf}.
\bibitem{Margolin1} L. G. Margolin, M. Shashkov: {\it Second-order sign-preserving conservative interpolation (remapping) on general grids},  J. Comput. Phys., 184(1):266-298, 2003.
\bibitem{Shashkov1} M.Shashkov: {\it Closure models for multidimensional cells in arbitrary Lagrangian-Eulerian hydrocodes},  Int. J. Numer. Meth. Fluids 56:1497-1504, 2008.
\bibitem{Vachal1} P. Vachal, P.-H. Maire: {\it Discretizations for weighted condition number smoothing on general unstructured meshes}, Comput. Fluids, 46(1):479-485, 2011.
\bibitem{Vachal2} P. Vachal, R.V. Garimella, M.J. Shashkov: {\it Untangling of 2D meshes in ALE simulations}, J. Comput. Phys., 196:627-644, 2004.
\bibitem{Youngs1} D. L. Youngs: {\it Multi-mode implosion in cylindrical 3D geometry}, 11th International Workshop on the Physics of Compressible Turbulent Mixing (IWPCTM11), Santa Fe, 2008.
\end{thebibliography}
\end{document}